\documentclass[12pt, a4paper]{report} 
\usepackage[english]{babel}
\usepackage{theorem}
\usepackage{amssymb}
\usepackage{amsfonts}
\usepackage{amsmath}
\usepackage[all]{xy}
\usepackage{geometry}

\newtheorem{thm}{Theorem}[section]
\newtheorem{prop}[thm]{Proposition}
\newtheorem{lem}[thm]{Lemma}
\newtheorem{cor}[thm]{Corollary}

\makeatletter
\begingroup
\gdef\th@upshape{\normalfont
  \def\@begintheorem##1##2{%
        \item[\hskip\labelsep \theorem@headerfont ##1\ ##2.]}%
\def\@opargbegintheorem##1##2##3{%
   \item[\hskip\labelsep \theorem@headerfont ##1\ ##2\ (##3).]}}
\endgroup
\makeatother

\theoremstyle{upshape}

\newtheorem{rem}[thm]{Remark}
\newtheorem{defn}[thm]{Definition}

\def\bg{\begin}
\def\ed{\end}

\def\a{\alpha}
\def\b{\beta}
\def\g{\gamma}
\def\d{\delta}
\def\e{\varepsilon}
\def\f{\varphi}
\def\k{\kappa}
\def\l{\lambda}
\def\ps{\psi}
\def\r{\rho}
\def\s{\sigma}
\def\S{\Sigma}
\def\ta{\theta}
\def\t{\tau}
\def\w{\omega}

\def\Om{\Omega}
\def\z{\zeta}
\def\del{\partial}
\def\DEL{\Delta}
\def\G{\Gamma}

\def\LL{\Lambda}

\def\bC{{\mathbf C}}
\def\bN{{\mathbf N}}
\def\bQ{{\mathbf Q}}
\def\bR{{\mathbf R}}
\def\bZ{{\mathbf Z}}
\def\bH{{\mathbf H}}

\def\bK{\mathbf{K}}
\def\bT{\mathbf{T}}
\def\nullset{\emptyset}

\def\A{{\cal A}}
\def\B{{\cal B}}
\def\C{{\cal C}}
\def\D{{\cal D}}
\def\E{{\cal E}}
\def\F{{\cal F}}
\def\H{{\cal H}}

\def\L{{\cal L}}
\def\cl{\ell}
\def\M{{\cal M}}
\def\N{{\cal N}}
\def\Q{{\cal Q}}
\def\O{{\cal O}}
\def\cP{{\cal P}}
\def\T{{\cal T}}
\def\X{{\cal X}}

\def\fg{{\mathfrak g}}

\def\fT{{\mathfrak T}}

\def\ae{\mathrm{a.e.}}
\def\Ad{\mathrm{Ad}}
\def\Aut{\mathrm{Aut}}

\def\ch{\mathrm{ch}}
\def\Cl{{\it Cl}}
\def\cl{{\it cl}}
\def\cvtu{\noindent \emph{\textbf{Proof}.\ }~}
\def\demo{\noindent \emph{\textbf{Proof}.\ }~}

\def\Der{\mathrm{Der}}
\def\dim{\mathrm{dim}}
\def\dist{\mathrm{dist}}
\def\dom{\mathrm{dom}}
\def\End{\mathrm{End}}

\def\GL{\mathrm{GL\,}}
\def\Ham{\mathrm{Ham}}
\def\id{\mathrm{id}}
\def\im{\mathrm{im}}
\def\Inv{\mathrm{Inv}}
\def\invt{\mathrm{inv}}
\def\ker{\mathrm{ker}}
\def\op{\mathrm{op}}
\def\Perm{\mathrm{Perm}}
\def\Pol{\mathrm{Pol}}
\def\ssa{\mathrm{sa}}
\def\sgn{\mathrm{sgn}}
\def\spin{\mathrm{spin}}

\def\th{\mathrm{th}}
\def\Tor{\mathrm{Tor}}
\def\th{^\mathrm{th}}
\def\tr{\mathrm{tr}}

\def\inv{^{-1}}
\def\cantrip{(\csm,L^2(S),\Dirac)}

\def\csm{C^{\infty}(M)}
\def\cx{C(X)}
\def\by{\times}
\def\oby{\otimes}
\def\obys{\wh{\otimes}}
\def\obyss{\,\overline{\otimes}}
\def\omd{(\Om,d)}

\newcommand{\Dirac}{{D \! \! \! \! / \;}}
\def\wed{\wedge}
\def\sseq{\subseteq}
\def\wt{\widetilde}
\def\wh{\widehat}
\def\ol{\overline}
\def\la{\left\langle}
\def\<{\left\langle}
\def\ra{\right\rangle}
\def\>{\right\rangle}
\def\bs{\backslash}
\def\mto{\mapsto}

\def\lfrom{\longleftarrow}

\def\lto{\longrightarrow}

\def\qed{\hfill\ensuremath{\square}\par\medskip}

\def\alg{algebra~}
\def\algn{algebra}
\def\algs{algebras~}
\def\algsn{algebras}

\def\calg{$C^{*}$-algebra~}
\def\calgs{$C^{*}$-algebras~}
\def\calgsn{$C^{*}$-algebras}
\def\calgn{$C^{*}$-algebra}

\def\cdns{conditions~}

\def\cts{continuous~}
\def\cpt{compact~}
\def\cptn{compact}
\def\cqms{\cpt quantum metric space~}
\def\cqmsn{\cpt quantum metric space}
\def\cqg{\cpt quantum group~}
\def\cqgn{\cpt quantum group}
\def\cqsg{\cpt quantum semigroup~}
\def\cqsgn{\cpt quantum semigroup}
\def\ctsn{continuous}
\def\cty{continuity~}

\def\cvg{convergence~}

\def\fn{function~}
\def\fnn{function}
\def\fns{functions~}

\def\iff{if, and only if,~}

\def\lcs{locally convex space~}
\def\lcsn{locally convex space}

\def\ncg{noncommutative geometry~}
\def\ncgn{noncommutative geometry}
\def\nc{noncommutative~}
\def\ncn{noncommutative}

\def\nbd{neighbourhood~}
\def\nbdn{neighbourhood}
\def\nbds{neighbourhoods~}
\def\nbdsn{neighbourhoods}

\def\pve{positive~}
\def\qft{quantum field theory~}
\def\qftn{quantum field theory}

\def\qgn{quantum group}
\def\qghd{quantum Gromov--Hausdorff distance~}
\def\qghdn{quantum Gromov--Hausdorff distance}

\def\qgh{quantum Gromov--Hausdorff~}

\def\sa{self-adjoint~}
\def\san{self-adjoint}
\def\ss{\subset}
\def\st{such that~}
\def\stn{such that}

\def\tsp{topological space~}

\def\tvs{topological vector space~}
\def\tvsn{topological vector space}
\def\vbd{vector bundle~}
\def\vbdn{vector bundle}
\def\vbds{vector bundles~}
\def\vbdsn{vector bundles}
\def\wrt{with respect to~}
\begin{document}
\thispagestyle{empty}

\begin{center}
\mbox{}
\vspace{3cm}

{\LARGE A Presentation of Certain New Trends in \\[3mm]
Noncommutative Geometry}

\vspace{8mm}

{\large R\'{e}amonn \'{O} Buachalla}
\vspace{15mm}

Thesis submitted for MSc by Research \\
National University of Ireland, Cork\\
2006

\vspace{3mm}

Department of Mathematics\\
Faculty of Science

\vspace{5mm}

\vspace{3mm}

Head of Department: Prof.\ J. Berndt\\
Supervisors: Prof.\ G. J. Murphy, Dr S. Wills
\end{center}

\pagebreak
\pagenumbering{arabic}

\tableofcontents
\chapter*{Introduction}

{\em Einstein was always rather hostile to quantum mechanics. How
can one understand this? I think it is very easy to understand,
because Einstein had proceeded on different lines, lines of pure
geometry. He had been developing geometrical theories and had
achieved enormous success. It is only natural that he should think
that further problems of physics should be solved by further
development of geometrical ideas. How to have $a \times b$ not equal
to $b \by a$ is something that does not fit in very well with
geometric ideas; hence his hostility to it.} ~~~~~~~~~~~~P. A. M.
Dirac

\bigskip

The development of quantum mechanics in the first half of the
twentieth century completely revolutionized classical physics. In
retrospect, its effect on many areas of mathematics has been no less
profound. The mathematical formalism of quantum mechanics was
constructed immediately after its birth by John von Neumann. This
formulation in turn gave birth to the theory of operator \algsn. In
the following years the work of mathematicians such as Israil
Gelfand and  Mark Naimark showed that by concentrating on the
function \alg of a space, rather than on the space itself, many
familiar mathematical objects in topology and measure theory could
be understood as `commutative versions' of operator
\algn-structures. The study of the \nc versions of these structures
then became loosely known as {\em \nc topology} and {\em \nc measure
theory} respectively. In the first chapter of this thesis we shall
present the work of von Neumann, Gelfand, Naimark and others in this
area.

From the middle of the twentieth century on, geometers like
Grothendieck, Atiyah, Hirzebruch and Bott  began to have great
success using algebraic formulations of geometric concepts that were
expressed in terms of the function \algs of spaces. $K$-theory,
which will be discussed in Chapter $3$, is one important example.
Later, it was realised that these algebraic formulations are well
defined for any \algn, and that the structures involved had
important roles to play in the operator theory. We cite algebraic
$K$-theory, and Brown, Fillmore, and Douglas' K-homology as primary
examples. However, the first person to make a serious attempt to
interpret all of this work as a \nc version of differential geometry
was a French mathematician named Alain Connes. Connes first came to
prominence in the $1970$s and he was awarded the Fields medal in
$1982$ for his work on von Neumann \algn s. Connes advanced the
mathematics that existed on the borderlines between differential
geometry and operator \algs further than anyone had previously
imagined possible. In $1994$ he published a book called {\em
Noncommutative Geometry} \cite{CON} that gave an expository account
of his work up to that time. It was an expanded translation of a
book he had written in French some years earlier called {\em
G\'{e}om\'{e}trie Noncommutative}. Noncommutative Geometry was
hailed as a `milestone for mathematics' by Connes' fellow Fields
medalist Vaughan Jones.

Connes developed \ncg to deal with certain spaces called singular
spaces. These arise naturally in many problems in classical
mathematics (usually as quotient spaces) but they are badly behaved
from the point of view of the classical tools of mathematics, such
as measure theory, topology and differential geometry. For example,
the space, as a topological space, may not be Hausdorff, or its
natural topology may even be the coarse one; consequently, the tools
of topology are effectively useless for the study of the problem at
hand. The idea of noncommutative geometry is to replace such a space
by a canonically corresponding noncommutative
${\text{$C^*$-algebra}}$ and to tackle the problem by means of the
formidable tools available in \ncgn. This approach has had enormous
success in the last two decades.

We give the following simple example: Let $X$ be a compact Hausdorff
space and let $G$ be a discrete group acting upon it. If the action
is sufficiently complicated then the quotient topology on $X/G$ can
fail to separate orbits, and in extreme cases $X/G$ can even have
the indiscrete topology. Thus, the traditional tools of mathematics
will be of little use in its examination. However, the \nc algebra
$C(X) \ltimes G$, the crossed product of $G$ with $C(X)$, is a
powerful tool for the study of $X/G$.

In the third and fourth chapters of this thesis we shall present the
fundamentals  of Connes' work.

\bigskip

In Connes' wake, \ncg has become an extremely active area of
mathematics. Applications have been found in fields as diverse as
particle physics and the study of the Riemann hypothesis. At the
start of the decade a group of European universities and research
institutes across seven countries formed an alliance for the
purposes of co-operation and collaboration in research in the theory
of operator algebras and noncommutative geometry; the alliance was
called the {\em European Union Operator Algebras Network}. The
primary aim of this thesis is to examine some of the areas currently
being explored by the Irish based members of the network. The Irish
institutions in question are the National University of Ireland,
Cork, and the Dublin Institute for Advanced Studies (DIAS).

In Chapters $3$ and $4$ we shall present some of the work being
pursued in Cork. Specifically, we shall present the work being done
on the interaction between Connes' \ncg and \cqgn s. Compact quantum
groups are a \nc generalisation of \cpt groups formulated by the
Polish mathematician and physicist Stanislaw Lech Woronowicz. The
relationship between these two theories is troublesome and ill
understood. However, as a research area it is showing great promise
and is at present the subject of very active investigation.

In the fifth chapter we shall present an overview  of some of the
work being done in Dublin. We shall focus on the efforts being made
there to use \nc methods in the renormalisation of quantum field
theory. This area is also a very active area of research; it is
known as {\em fuzzy physics}. We shall not, however, present a
detailed account of the DIAS work. Its highly physical nature is not
well suited to a pure mathematical treatment. Instead, an exposition
of Marc Rieffel's \cqms theory will be given in the sixth chapter.
Compact quantum metric spaces were formulated by Rieffel with the
specific intention of providing a sound mathematical framework in
which to discuss fuzzy physics. Indeed, the papers of DIAS members
are often cited by Rieffel when he discusses \cqmsn s.


\chapter{Motivation}

As explained in the introduction, \ncg is based upon the fact that
there exist a number of correspondences between  basic mathematical
structures and the commutative versions of certain operator \algn ic
structures. In this chapter we shall present the three prototypical
examples of these correspondences: the Gelfand--Naimark Theorem, the
Serre--Swan Theorem, and the characterisation of commutative von
Neumann Algebras. We shall also show how these results naturally
lead us to define `\nc versions' of the mathematical structures in
question.

\section{The Gelfand--Naimark Theorem}

The Gelfand--Naimark Theorem is often regarded as the founding
theorem of \ncgn. Of all the results that we shall present in this
chapter it was the first to be formulated, and, from our point of
view, it is the most important. It is the theorem that motivated
mathematicians to consider the idea of a `\nc generalisation'  of
locally \cpt Hausdorff spaces.

The theorem, and indeed most of \ncgn, is expressed in the langauge
of \calgsn. Thus, we shall begin with an exposition of the basics of
these \algsn, and then progress to a proof of the theorem itself.

\subsection{\calgs}
The motivating example of a \calg is  $B(H)$ the normed \alg of
bounded linear operators on a Hilbert space $H$. The algebra
operations of $B(H)$ are defined pointwise and its norm,  which is
called the {\em operator norm}, is defined by
${\|A\|=\sup\{\|A(x)\|:\|x\|\leq 1\}}$,  $A \in B(H)$. With
respect to these definitions, $B(H)$ is a {\em Banach algebra};
that is, an \alg with a complete submultiplicative norm.
Completeness of the norm is a standard result in Hilbert space
theory, and submultiplicativity follows from the fact that
\[
\|AB\|=\sup_{\|x\| \leq 1} \|AB(x)\| \leq \sup_{\|x\| \leq 1}
(\|A\|\|B(x)\|) =\|A\|\|B\|,
\]
for all $A,B \in B(H)$.

A central feature of Hilbert space operator theory is the fact that
for every ~~~~${A \in B(H)}$, there exists an operator $A^* \in
B(H)$, called the {\em adjoint} of $A$, \st
\[
\<Ax,y\>=\<x,A^*y\>,\text{~~~~~for all~}x,y \in H.
\]
As is well known, and easily verified, $(A + \l
B)^*=A^*+\ol{\l}B^*$, $A^{**}=A$, and $(AB)^*=B^*A^*$, for all $A,B
\in B(H)$, $\l \in \bC$. Another well known, and very important,
equation that involves an operator and its adjoint is
\[
\|A^*A\|=\|A\|^2, \text{~~~~~for all~}A \in B(H).
\]

\bigskip

Now that we have reviewed the relevant features of Hilbert space
operator theory, we are ready to begin generalising. A {\em
$*$-\alg} is an \alg $A$, together with a mapping
\[
*:A \to A, ~~~~ a\mapsto a^{*},
\]
called an {\em \alg involution}, such that for all $a,b \in A,~\l
\in \bC$;
\begin{enumerate}
\item  $(a+\l b)^{*}=a^{*}+\bar{\l}b^{*}$ ({\em conjugate-linearity}),%
\item  $~a^{**}=a$ ({\em involutivity}), %
\item  $(ab)^{*}=b^{*}a^{*}$ ({\em anti-multiplicativity}).
\end{enumerate}

If $A$ and $B$ are two $*$-\algs and $\f$ is an \alg homomorphism
from $A$ to $B$, then $\f$ is called a {\em $*$-\alg homomorphism}
if $\f(a^*)=\f(a)^*$, for all $a \in A$. If the homomorphism is
bijective, then it is called a {\em $*$-isomorphism}. If $a \in A$
\st $a^*=a$, then $a$ is called {\em self-adjoint}; we denote the
subset of self-adjoint elements of $A$ by $A_{\ssa}$.  If $X$ is a
subset of $A$ \st $X^*=\{x^*:x \in X\} = X$, then $X$ is called a
{\em self-adjoint subset} of $A$.

\bigskip

\bg{defn} A {\em \calg } $\A$ is a $*$-\alg that is also a Banach
\algn, \st
\[
\|a^*a\| = \|a\|^2,~~~~ \text{for all } a \in \A.
\]
\ed{defn}

From our comments above we see that $B(H)$ is a \calg for every
Hilbert space $H$. The simplest example of a \calg is the complex
numbers $\bC$. Any closed self-adjoint sub\alg of a \calg is clearly
also a \calgn.

One of the most important examples is $C_{0}(X)$, the \alg of \cts
\fns that vanish at infinity on a locally compact Hausdorff space
$X$. (We recall that for a \fn $f \in C_{0}(X)$, to {\em vanish at
infinity} means that for each $\e>0$ there exists a compact subset
$K \sseq X$ such that $|f(x)| < \e$, for all $x\notin K$.) If we
endow $C_0(X)$ with the standard supremum norm and define
\[
f^*(x)=\overline{f(x)}, \qquad \text{ for all } x \in X,
\]
then the conditions for a \calg are clearly fulfilled.

This \alg is one of the standard examples of a \calg that is not
necessarily unital. Notice that if $X$ is non-compact, then $1
\notin C_{0}(X)$. On the other hand, if $X$ is \cptn, then
$C_0(X)=C(X)$, the \alg of \cts \fns on $X$, and so it is unital.

We should also notice that $C_0(X)$ is a commutative \calgn. In
fact, as we shall see later, all commutative \calgs are of this
form.

\subsubsection{The Universal Representation}

As we have seen,  \calgs  generalise the \algs of bounded linear
operators on Hilbert spaces. A {\em representation} of a \calg $\A$
is a pair $(U,H)$ where $H$ is a Hilbert space and $U:\A \to B(H)$
is a $*$-\alg homomorphism. If $U$ is injective, then the
representation is called {\em faithful}. What is interesting to note
is that for any \calg $\A$ there exists a distinguished faithful
representation of $\A$ called its {\em universal representation}.
Thus, every \calg is isometrically \mbox{$*$-isomorphic} to a closed
$*$-sub\alg of $B(H)$, for some Hilbert space $H$. This result is
due to Gelfand and Naimark, for further details see \cite{MUR}.

\subsubsection{Unitisation}

If $\A$ is  a non-unital \calgn, then it often proves useful to
embed it into a unital \calg using a process known as {\em
unitisation}: One starts with the linear space $\A \oplus \bC$ and
defines a multiplication on it by setting
\begin{equation}\label{1:eqn:unitmult}
(a, \l).(b,\mu) = (ab + \l b + \mu a, \l \mu),
\end{equation}
and an involution by setting
\[
(a, \l)^*=(a^*,\ol{\l}).
\]
One then defines a norm on it by setting
\[
\|(a,\l)\|=\sup\{\|ab+\l b\|:b \in \A,\|b\| \leq 1\}.
\]
It is not too hard to show that this norm makes it a \calgn. (As we
shall see later, this norm is necessarily unique.) The element
$(0,1)$ clearly acts as a unit. We denote this new unital \calg by
$\wt{\A}$, and embed $\A$ into it in the canonical manner. It is
easily seen to be an isometric embedding.

\subsubsection{Spectrum of an Element}

Let $\A$ be a unital \calg and $a \in \A$. We define the {\em
spectrum} of $a$ to be the set
\[
\s(a) = \{ \l \in \bC : (\l 1 -a) \notin \Inv(\A) \},
\]
where $\Inv(\A)$ is the set of invertible elements of $\A$. If $a$
is an element of a non-unital \alg then we define its {\em spectrum}
to be the set
\[
\s(a) = \{ \l \in \bC : (\l 1 -a) \notin \Inv(\wt{\A}) \}.
\]
Three basic facts about the spectrum of an element are that: $\s(a)
\neq \nullset $, for any $a \in \A$; if $a=a^*$, then $\s(a) \sseq
\bR$; and if $\l \in \s(a)$, then $|\l| \leq \|a\|$. (These results
hold for the unital or non-unital definitions of the spectrum.)

We define the \textit{spectral radius} of an element $a$ to be
\[
r(a)=\sup\{|\l|:\l \in \s(a)\};
\]
note that this is a purely algebraic definition and takes no account
of the algebra's norm. It is a well known result of Beurling that,
for any \calg $\A$,
\[
r(a)= \lim_{n \to \infty}\|a^{n}\|^{\frac{1}{n}},\text{~~~~~for
all~}a \in \A.
\]
Now,~if $a \in \A$ is \san, then $\|a^2\|=\|a^*a\|=\|a\|^2$. Thus,
\begin{equation} \label{1:eqn:saradius}
r(a)=\lim_{n \to \infty}\|a^n\|^{\frac{1}{n}}=\lim_{n \to
\infty}\|a^{2^n}\|^{\frac{1}{2^n}}=\lim_{n \to \infty}\|a\|=\|a\|.
\end{equation}
This result will be of use to us in our proof of the
Gelfand--Naimark Theorem.

As a more immediate application, we can use it to show that there is
at most one norm on $*$-\alg making it a \calgn. Suppose that $\|
\cdot \|_1$ and $\| \cdot \|_2$ are two norms on a $*$-\alg $\A$
making a \calgn. Then they must be equal since
\[
\|a\|_1^2=\|aa^*\|_1=r(aa^*)=\|aa^*\|_2=\|a\|^2_2,
\]
for all $a \in \A$. A useful consequence of this fact is that any
$*$-isomorphism is an isometric mapping.

\subsection{The Gelfand Transform}

A {\em character} $\f$ on an algebra $A$ is a non-zero algebra
homomorphism from $\A$ to $\bC$; that is, a non-zero linear
functional $\f$ such that
\[
\f(ab)=\f(a)\f(b), \text{~~~~~for all~} a,b \in A.
\]
We denote the set of all characters on $A$ by $\Om(A)$.

Let $\A$ be a unital \calg and $\f \in \Om(\A)$, then
$\f(a)=\f(a1)=\f(a)\f(1)$, for all $a \in \A$, and so $\f(1)=1$. If
$a \in \Inv (\A)$, then $1=\f(a a \inv)=\f(a)\f(a \inv)$. Thus, if
$\f(a)=0$, $a$ cannot be invertible. The fact that $\f(\f(a)-a)=0$
then implies that $(\f(a)-a) \notin \Inv(\A)$. Thus, $\f(a) \in
\s(a)$, for all $a \in \A$, and for all $\f \in \Om(\A)$.

If we now recall that $|\l| \leq \|a\|$, for all $\l \in \s(a)$,
then we see that
\begin{equation} \label{eqn1:characterradius}
\|\f\|= \sup\{|\f(a)|: \|a\| \leq 1,~a \in \A\} \leq 1.
\end{equation}
(In fact, in the unital case $\|\f\|=1$, since $\f(1)=1$.)  It
follows that each $\f \in \Om(\A)$ is norm \ctsn, and that $\Om(\A)$
is contained in $\A^*_1$; where $\A^*$ is the space of bounded
linear functionals on $\A$, and $\A^*_1$ is the closed unit ball of
$\A^*$.

Another fact about characters is that they are {\em Hermitian},
that is, $\f(a^*)=\ol{\f(a)}$, for all $\f \in \Om(\A)$. This is
shown as follows: every element $a$ of $\A$ can be written
uniquely in the form $a=a_{1}+ia_{2}$, where $a_{1}$ and $a_{2}$
are the two self-adjoint elements
\[
a_{1}=\frac{1}{2}(a+a^{*}), \quad a_{2}=\frac{1}{2i}(a-a^{*}).
\]
Clearly,
\[
\f(a^{*})=\f((a_{1}+ia_{2})^{*})=\f(a_{1}-ia_{2})=\f(a_{1})-i\f(a_{2}).
\]
Then, since  $\f(a_{i}) \in \s(a_{i})$, and $\s(a_{i}) \subset \bR$,
for $i=1,2$, we have that
\[
\f(a_{1})-i\f(a_{2})=\overline{\f (a_{1})+i\f
(a_{2})}=\overline{\f(a)}.
\]

\subsubsection{Gelfand Topology}
We now find it convenient to endow $\Om(\A)$ with a topology
different from the norm topology. We do this by putting the weak$^*$
topology on $\A^*$, and restricting it to $\Om(\A)$. Recall that the
{\em weak$^*$ topology} on $\A^*$ is the weakest topology with
respect to which all maps of the form
\begin{equation}\label{1:eqn:weakstar}
\wh{a}:\A^* \to \bC, ~~~~ \f \mto \f(a),~~~~a \in \A,
\end{equation}
are \ctsn. Of course, we may alternatively describe it as the
weakest topology \wrt which a net $\{\f_{\l}\}_{\l}$ in $\A^*$
converges to $\f \in \A^*$ \iff $~\wh{a}(\f_{\l}) \to \wh a(\f)$; or
equivalently, \iff $~\f_{\l}(a) \to \f(a)$, for all $a \in \A$.

When $\Om(\A)$ is endowed with this topology we call it the {\em
spectrum} of $\A$. We call the topology itself the {\em Gelfand
topology}. We should note that the spectrum is a Hausdorff space
since the weak$^*$ topology is a Hausdorff topology.

With respect to the  weak$^*$ topology $\Om(\A) \cup \{0\}$ is
closed in $\A^*$. To see this take a net $\{\f_{\l}\}_{\l}$  in
$\Om(\A) \cup \{0\}$ that converges to $\f \in \A^*$. Now,~ $\f_\l$
is bounded for each $\l$, and so
\[
\f(ab)=\lim_{\l}\f_{\l}(ab)=\lim_{\l}\f_{\l}(a)\f_{\l}(b)=\f(a)\f(b),
\]
for all $a,b \in \A$. Thus, if $\f$ is non-zero, then it is a
character. Either way however, $\f$ is contained in $\Om(\A) \cup
\{0\}$, and so $\Om(\A) \cup \{0\}$ is closed.

If $\A$ is unital, then taking $\{\f_{\l}\}_{\l}$ and $\f$ as above,
we have that
\[
\f(1)= \lim_{\l}\f_{\l}(1)=1 \neq 0.
\]
Hence the zero functional lies outside the closure of $\Om(\A)$, and
so $\Om(\A)$ is closed.

We should now note that  equation (\ref{eqn1:characterradius})
implies that $\Om(\A) \cup \{0\}$ is contained in $\A^{*}_1$, which
is weak$^*$ compact by the Banach--Alaoglu Theorem. Hence, $\Om(\A)$
is compact in the unital case, and locally compact in the non-unital
case.

\subsubsection{Gelfand Transform}
Define a mapping $\G$, called the {\em Gelfand transform}, by
setting
\[
\G:\A \to C_0(\Om(\A)), ~~~~ a \mto \wh a;
\]
where by $\wh a$ we now mean the mapping defined in
(\ref{1:eqn:weakstar}) with a domain restricted to $\Om(\A)$. The
image of $\A$ under the Gelfand transform is contained in
$C_0(\Om(\A))$. To see this, choose an arbitrary $\e>0$ and
consider the set
\[
\Om_{\wh a,\e}=\{\f \in \Om(\A) : |\widehat{a}(\f)| \geq \e \}.
\]
Using an argument similar to the one above, we can show that this
set is weak$^*$ closed in $\A_1^{*}$. Thus, by the Banach--Alaoglu
Theorem, it is compact. Now,~for $\f \notin \Om_{\wh a , \e}$ we
have $|\widehat{a}(\f)| < \e $, and so $\widehat{a} \in
C_{0}(\Om(\A))$. However, if $\A$ is unital, then this tells us
nothing new since $\Om(\A)$ will be compact and $C_0(\Om(\A))$ will
be equal to $C(\Om(\A))$.

\subsubsection{Gelfand Transform and the Spectral Radius}

We saw earlier that if $\f \in \Om(\A)$, then $\f(a) \in \s(a)$,
for all $a \in \A$. If $\A$ is assumed to be commutative and
unital, then it can be shown that all elements of $\s(a)$ are of
this form; that is,
\[
\s(a)=\{\f(a)| \f \in \Om(A)\}.
\]
In the non-unital case we almost have the same result: if $\A$ is
a non-unital commutative \calgn, then
\[
\s(a)=\{\f(a)|\f \in \Om(A)\} \cup \{0\}.
\]
These two results are important because they give us the equation
\[
r(a)=\max\{|\l| :\l \in \s(a)\}=\max\{|\l| :\l \in
\widehat{a}(\Om(A))\}=\|\wh{a}\|.
\]
Hence, for any commutative \calg $\A$, we have that
\begin{equation} \label{1:eqn:radiusgelfand}
r(a)=\|\widehat{a}\|, ~~~~ \text{for all } a \in \A.
\end{equation}
It is instructive to note that this is the point at which the
requirement of commutativity (which is needed for the
Gelfand--Naimark theorem to hold) enters our discussion of \calgsn.

\subsection{Gelfand--Naimark Theorem}

The following result is of great importance. It allows us to
completely characterise commutative $C^{*}$-algebras. It first
appeared in a paper \cite{GN} of Gelfand and Naimark in $1943$,
and it has since become the principle theorem motivating
noncommutative geometry.

\begin{thm}
Let $\A$ be an abelian $C^{*}$-algebra. The Gelfand transform
\[
\G:\A \to C_{0}(\Om (\A)),~~~~ a \mapsto \widehat{a},
\]
is an isometric $*$-isomorphism.
\end{thm}
\demo That $\G$ is a homomorphism is clear from
\[
\G(a + \l b)(\f)=\wh{a + \l b}(\f)=\f(a)+\l \f(b)=\G(a)(\f)+\l
\G(b)(\f),
\]
and
\[
\G(ab)(\f)=\wh{ab}(\f)=\f(ab)=\f(a)\f(b)=(\G(a)\G(b))(\f).
\]
Hence, $\G(\A)$ forms a sub\alg of $C_0(\Om(\A))$.

If  $\f,\ps \in \Om(\A)$ and $\f \neq \ps$, then there exists an $a
\in \A$ such that $\f(a) \neq \ps(a)$. Hence, $\widehat{a}(\f)\neq
\widehat{a}(\ps)$, and so $\G(\A)$ separates points of $\Om(\A)$.

If $\f \in \Om(\A)$, then by definition $\f$ is non-zero. Therefore,
there exists an $a \in \A$ such that $\f(a) \neq 0$, and so $
\widehat{a}(\f) \neq 0$. Hence, for each point of the spectrum there
exists a function in $ \G(\A)$ that does not vanish there.

It is seen that $\G(\A)$ is closed under conjugation from
\[
\widehat{a}^*(\f)=\overline{\widehat{a}(\f)}=\overline{\f(a)}=\f(a^{*})=\widehat{a^{*}}(\f).
\]

We shall now show that $\G$ is an isometry. Firstly, we note that
since $aa^*$ is \sa
\[
\|aa^*\|=r(aa^*)=\|\widehat{aa^{*}}\|.
\]
Then  since $\G(\A)$ is a subset of the \calg $C_0(\Om(\A))$, we
have that ${\|\wh{a}\|^2=\|\wh{a}^*\wh{a}\|}$. It follows  that
\[
\|\widehat{a}\|^{2}=\|{\widehat{a}}\widehat{a}^*\|=\|\widehat{aa^*}\|=\|aa^*\|=\|a\|^{2}.
\]
Therefore, $\|a\|=\|\widehat{a}\|$ and $\G$ is an isometry. This
implies that $\G(\A)$ is complete, and therefore closed in
$C_0(\Om(\A))$. Since $\Om(\A)$ is locally compact, we can invoke
the Stone--Weierstrass Theorem  and conclude that
$\G(\A)=C_{0}(\Om(\A))$. \qed

\bg{cor} For $\A$ a unital abelian $C^{*}$-algebra, the Gelfand
transform
\[
\A \to C(\Om (\A)), ~~~~ a \mapsto \widehat{a},
\]
is an isometric $*$-isomorphism
\end{cor}

\subsection{The Algebra-Space Correspondence}

Following the definition of a $C^{*}$-algebra, we saw that if $X$
is a locally compact Hausdorff space, then $C_{0}(X)$ is a
commutative $C^{*}$-algebra. Therefore, every locally compact
Hausdorff space is naturally associated with an abelian \calgn.
Conversely, the Gelfand--Naimark Theorem showed us that every
abelian \calg $\A$ is naturally associated with a locally compact
Hausdorff space, namely $\Om(\A)$. If we could show that these
associations were inverse to each other, then we would have
established a one-to-one correspondence between spaces and
\algsn. In fact, this is a direct consequence of the following
theorem.

\begin{thm} \label{1:thm:cptsurjective} If $X$ is a locally compact
Hausdorff space, then $X$ is homeomorphic to $\Om(C_{0}(X))$.
\end{thm}
\demo We shall prove this result in the compact case only and we
refer the interested reader to the first chapter \cite{ZHU}.

Let us define the mapping
\[
\F:X \to \Om(C(X)),~~~~ x \mapsto \F_{x},
\]
by setting $\F_{x}(f)=f(x)$, for all $f \in C(X)$. It is clear that
$\F_x$ is a character, for all $x \in X$.

If $x_{1} \neq x_{2}$ in $X$, then it follows from Urysohn's Lemma
that there exists an $f$ in $C(X)$ such that $f(x_{1})=0$ and
$f(x_{2})=1$. This shows that  $\F$ is injective.

To show that $\F$ is surjective, take any $\ps \in \Om(C(X))$ and
consider the ideal
\[
I=\ker(\ps)=\{f \in C(X):\ps(f)=0\}.
\]
We shall show that there exists an $x_{0} \in X$ such that
$f(x_{0})=0$, for all $f \in I$. If this is not the case, then for
each $x \in X$, there is an $f_{x} \in I$ such that $f_{x}(x) \neq
0$. The continuity of $f$ implies that each $x$ has an open
neighborhood $U_{x}$ on which $f_{x}$ is non-vanishing. By the
compactness of $X$, there exist $x_{1},...,x_{n}$ in $X$ such that
$X=\bigcup_{k=1}^{n}U_{x_{k}}$. Let
\[
f(x)=\sum_{k=1}^{n}|f_{x_{k}}(x)|^{2}, \quad x \in X.
\]
Clearly $f$ is non-vanishing on $X$, and so it is invertible in
$C(X)$. This implies that $\ps(f) \neq 0$. On the other hand, since
$\ps$ is multiplicative,
\[
\ps(f)=\sum_{k=1}^{n}\ps(f_{x_{k}})\ps(\overline{f_{x_{k}}})=0
\]
This contradiction shows that there must exist some $x_{0} \in X$
such that {$f(x_{0})=0$}, for all $f \in I$.

Now,~if $f$ is an arbitrary element in $C(X)$, then $f-\ps(f)$ is in
$I$. Thus,
\[
(f-\ps(f))(x_0)=f(x_{0})-\ps(f)=0,
\]
or equivalently $\ps(f)=f(x_{0})$. Therefore, $\ps=\F(x_{0})$, and
so $\F$ is surjective.

If $(x_{\l})_{\l}$ is a net in $X$ that converges to $x$, then
$f(x_{\l}) \to f(x)$, for every $f \in C(X)$. This is equivalent to
saying that, $\ps_{x_{\l}}(f) \to \ps_x(f)$, for all $f \in C(X)$;
or that $\wh{f}(\ps_{x_{\l}}) \to \wh{f}(\ps_x)$, for all $f \in
C(X)$. Hence, $\ps_{x_{\l}} \to \ps_x$ \wrt the weak$^*$ topology,
and so $\F$ is a \cts mapping.

We have now shown that $\F$ is bijective continuous function from
the compact Hausdorff space $X$ to the compact Hausdorff space
$\Om(C(X))$. Therefore, we can conclude that it is a homeomorphism.
\qed


\section{Noncommutative Topology}
We shall now stop and reflect on what we have established: we have
shown that there is a one-to-one correspondence between
commutative $C^{*}$-algebras and locally compact Hausdorff spaces.
Thus, the \fn \alg of a space contains all the information about
that space. This means that nothing would be lost if we were to
study the \alg alone and `forget' about the space. This approach
is common in other areas of mathematics, most notably algebraic
geometry. The idea behind noncommutative topology is to take it
one step further: loosely speaking, noncommutative topology views
\cpt Hausdorff spaces as special commutative examples of general
\calgn s, and studies them in this context. The subject can be
described as the investigation of those \calgn ic structures that
correspond to topological structures when the \alg is commutative.

As often happens when one works in greater generality, results that
were previously complex or technical become quite straightforward
when this approach is used. It has allowed formerly `unsolvable'
problems in topology to be solved. Also, the application of our
topological intuition to \nc \calgs has helped in the discovery of
new algebraic results. Quite often, the algebraic structures that
generalise topological structures are much richer and have features
with no classical counterparts.

\bigskip

However, in all of this one does notices a lack of duality. While
algebraic properties of commutative \calgs correspond to geometric
properties, no such correspondence exists for \nc \algsn. This
prompts us to consider the possibility that, in a very loose
sense, every \nc \calg could be viewed as the {\em \fn \alg of
some type of `\ncn'} or {\em `quantum space'}. This is the basic
heuristic principle upon which most of the vocabulary of
noncommutative topology is based.

Quantum spaces are imagined to be a type of generalised set, and
those spaces that have a classical point set representation are
considered to be special case. Noncommutative topology is then
thought of as the investigation of these quantum spaces through
their  \fn \algsn.

It must be stressed, however, that quantum spaces only exist as an
intuitive tool based upon an analogy. They do not have any kind of
proper definition. It is only when one ventures into the literature
of physics that the concept gains any concrete form.

\subsection{Some Noncommutative Generalisations}

Now that we have presented the general philosophy behind  \nc
topology we can move on and explore some more concrete aspects of
the subject. In this section we shall present some simple examples
of  generalisations of topological properties and structures to the
\nc setting.

\begin{enumerate}
\item As we saw above, points in a space $X$ are in one-to-one
correspondence with characters on the \calg $C_0(X)$. Therefore,
we view characters as the appropriate generalisation of points to
the \nc case:
\[
\bg{tabular}{lll} points & $\lto$ & characters.
\end{tabular}
\]

However, it must be noted that the non-emptiness of the spectrum of
a \nc \calg is not guaranteed. The easiest example of such a \calg
is $M_n(\bC)$, for $n \geq 2$. It is well known $M_n(\bC)$ has no
proper ideals for $n \geq 2$. However, if $\f$ is a character on any
\calg $\A$, then its kernel is clearly a proper two-sided proper
ideal of $\A$ (in fact, it is a maximal ideal). Therefore, the
spectrum of $M_n(\bC)$ must be empty.

For this reason the notion of a point will be of little use in the
noncommutative world. In fact, it is often written that `quantum
spaces are pointless spaces'.

\item

Let $\F$ be a homeomorphism from a locally \cpt Hausdorff space $X$
to itself. Consider the map
\[
\Phi:\A \to \A,~~~~f \mto f \circ \F,
\]
where $\A=C(X)$. A little thought will verify that it is a
$*$-isomorphism. Thus, we can associate a $*$-isomorphism to each
homeomorphism.

On the other hand, let $\Phi$ be an \alg $*$-isomorphism from  $\A$
to itself, and consider the mapping
\[
\F:\Om(\A) \to \Om(\A),~~~~\f \mto \f \circ \Phi.
\]
With the aim of establishing the \cty of $\F$, we shall consider a
net $(\f_\l)_\l$ in $\Om(\A)$ that converges to $\f$. Since
$\Om(\A)$ is endowed with the weak$^*$ topology, $\f_\l(a) \to
\f(a)$, for all $a \in \A$. Now,~ $\F(\f_\l)(a)=\f_\l(\Phi(a))$, and
since $\f_\l(\Phi(a)) \to \f(\Phi(a))$ and $\f(\Phi(a))=\F(\f)(a)$,
 $\F$ must be \ctsn.

Again, it is straightforward to show that $\F$ is bijective. Using
an argument similar to that above, we can also show that  $\F
\inv$ is \ctsn, and so $\F$ is a homeomorphism.

Because of the Gelfand--Naimark Theorem $\F$ can also be
considered as a homeomorphism  from $X$ to $X$. Thus, to each
$*$-isomorphism we can associate a homeomorphism.

A little extra work will verify that these two associations are
inverse to each other. Hence, we have established a one-to-one
correspondence between homeomorphisms and $*$-isomorphisms. This
motivates our next generalisation:

\[
\bg{tabular}{lll}
homeomorphisms & $\lto$ & $*$-isomorphisms. \\
\end{tabular}
\]

\item If $f \in C_0(X)$, for some \cpt Hausdorff space $X$, then
clearly $\l \in \s(f)$ \iff $f(x)=\l$, for some $x \in X$. (In the
locally \cpt case $\s(f)=\im(f) \cup \{0\}$.) Thus, the spectrum of
an element of a \calg is a generalisation of the image of a \fnn:
\[
\bg{tabular}{lll}
image of a \fn & $\lto$ & spectrum of a element. \\
\end{tabular}
\]

This motivates us to define a {\em \pve element} of a \calg $\A$
to be a self-adjoint element with positive spectrum; we denote
that $a \in \A$ is \pve by writing $a \geq 0$, and we denote the
set of \pve elements of $\A$ by $\A_+$. The self-adjointness
requirement is necessary because there may exist non-self-adjoint
elements of $\A$ with positive spectrum. One merely needs to look
in $M_2(\bC)$ for an example.

It is pleasing to note that, just as in the classical case, every
\pve element is of the form $a a^*$, for some $a \in \A$. Note
that this implies that when $\A=B(H)$, for some Hilbert space $H$,
then the \pve elements of $\A$ will coincide with the positive
operators on $H$.

We now record this generalisation:
\[ \bg{tabular}{lll}
\pve \fn & $\lto$ & \pve element. \\
\end{tabular}
\]

\item
Let $X$ be a \cpt Hausdorff space, and let $\mu$ be a regular
complex Borel measure on $X$. The scalar-valued function $\l$ on
$C_0(X)$, defined by $\l(f)=\int_{X}f d\mu$, is clearly bounded and
linear by the properties of the integral. Therefore, it is an
element of $C_0(X)^*$. The following well known theorem shows us
that every bounded linear functional arises in this way.

\begin{thm}[Riesz Representation Theorem]
Let $X$ be a locally compact Hausdorff space. For all $\l \in C_0
(X)^*$ there exists a unique regular complex Borel measure $\mu$ on
$X$ \st
\[
\l(f)=\int_{X}f d\mu, ~~~~ f \in C(X).
\]
\end{thm}

Thus, there exists a one-to-one  correspondence between the elements
of the dual space of $C_0 (X)$ and the regular complex Borel
measures on $X$. This motivates our next generalisation:
\[
\bg{tabular}{lll} regular complex Borel measures  & $\lto $ &
bounded linear functionals.
\end{tabular}
\]

\item As we noted earlier, if $X$ is a \cpt Hausdorff space, then
$C_0(X)=C(X)$ is a unital \algn; and if $X$ is a non-\cpt Hausdorff
space, then $C_0(X)$ is a non-unital \algn.

As is well known, we may compactify a non-\cpt space by adding to it
a point at infinity. We denote this new space by $X_{\infty}$. The
\alg $C_0(X_{\infty})=C(X_\infty)$ is then a unital \algn. If we
unitise $C_{0}(X)$, then we also get a unital \alg $\wt{C_0(X)}$.
What is interesting about this is that $C(X_\infty)$ is
isometrically $*$-isomorphic to $\wt{C_0(X)}$. An obvious
isomorphism is
\[
C(X_\infty) \to \wt{C_0(X)},~~~~  f \mto
((f-f(\infty)1)|_X,f(\infty)).
\]
These observations motivate the following generalisation:

\[
\bg{tabular}{lll}
compact spaces & $\lto$ & unital \calgn s, \\
one point compactifaction & $\lto$ & unitisation.
\end{tabular}
\]

\item Each closed subset $K$ of a compact Hausdorff space $X$ is
a compact Hausdorff space. Through our \algn-space correspondence,
$X$ and $K$ are associated to the \calg of  \cts \fns defined upon
them. With the aim of generalising closed sets to the \alg setting,
we shall examine the relationship between these two \calgsn.

Let $K$ be a closed subset of $X$ and write
\begin{equation} \label{eqn:1:ideal}
I=\{f \in C(X):f|_{K}=0\}.
\end{equation}
Clearly $I$ is a closed  {\em $*$-ideal} of $C(X)$; that is, a
closed self-adjoint ideal. This implies that $C(X)/I$ is well
defined as a $*$-\algn. We can define a norm on it by
\[
\|f+I\|=\inf_{h \in I}\|f + h\|_{\infty}.
\]
It is a standard result that when $C(X)/I$ is endowed with this
norm it is a Banach \algn. In fact, as can be routinely verified,
it is a \calgn.

Let us now consider the mapping
\[
\pi:C(X) \to C(K),~~~~f \mto f|_K.
\]
For any $g \in C(K)$, the Tietze extension theorem implies that
there exists an $f \in C(X)$, \st $f$ extends $g$. Thus, $\pi$ is
surjective. Since $\pi$ is  clearly a $*$-\alg homomorphism with
kernel equal to $I$, it induces a $*$-isomorphism from $C(X)/I$ to
$C(K)$.

If we could show that all the closed $*$-ideals of $C(X)$ were of
the same form as (\ref{eqn:1:ideal}), then we would have a
one-to-one correspondence between quotient \algs and  closed
subsets. The following lemma shows exactly this (for a proof see the
first chapter of \cite{ZHU}).

\begin{lem} Let $X$ be a compact Hausdorff space and let $I$ be a
closed $*$-ideal of $C(X)$. Then there exists a closed subset $K$ of
$X$ such that
\[
I=\{f \in C(X):f|_{K}=0\}.
\]
\end{lem}

\bigskip

Using a similar line of argument one can also show that there is a
bijective correspondence between the open sets of $X$ and the
ideals in $C(X)$.

This gives us our next generalisations:
\[
\bg{tabular}{lll}
closed sets of a \cpt space & $\lto $ & quotients of  unital \calgn s, \\
open sets of  a \cpt space & $\lto $ & ideals of unital \calgsn.
\end{tabular}
\]

\item With the aid of a simple definition and a standard result we
can generalise connectedness.

A {\em projection} in a \calg $\A$ is an element $a \in \A$ such
that
\[
a^{2}=a^{*}=a.
\]
An $*$-\alg $A$ is called {\em projectionless} if the only
projections it contains are $0$, and $1$ if $A$ is unital.

Let $X$ be a connected space. If $p \in C(X)$ is a projection,
then for all $x \in X$,  $(p(x))^2=p(x)$, so $p(x)$ is equal to
$0$ or $1$. Since $X$ is connected, $X=p \inv \{0\} \cup p \inv
\{1\}$ cannot be a disconnection of $X$, therefore $p$ is equal to
$0$ or $1$.

Conversely, if $C(X)$ is projectionless, then $X$ must be
connected because a non-trivial projection can easily be defined
on an unconnected space. This motivates the following
generalisation:
\[
\bg{tabular}{lll} connected \cpt space& $\lto$ & projectionless
unital \calgn.
\end{tabular}
\]

\item

\begin{thm}
If a compact Hausdorff space $X$ is metrisable, then $C(X)$ is
separable.
\end{thm}

\demo Let $X$ be metrisable with a metric $d$, and denote
\[
{B_r(x)=\{y: d(x,y)<r,~y \in X\}}, ~~~~ \text{for }r >0.
\]
This means that $\C_r=\{B_r(x):x \in X\}$ is an open cover of $X$.
Since $X$ is compact, $\C_r$ has a finite subcover $\C'_r$, for all
$r \geq 0$. The family $\C=\bigcup_{n=1}^{\infty}\C'_{\frac{1}{n}}$
will then form a countable base for the topology of $X$.

Let $x$ be an element of $B_1 \in \C$. Since $X$ is a compact
Hausdorff space, it is easy to see that there exists a $B_2 \in \C$,
\st
\begin{equation}
\label{eqn1:topbasis} x \in B_2 \ss \ol{B_2} \ss B_1.
\end{equation}
Thus, the countable family
\[
\D=\{(B_1,B_2) \in \C \by \C : \ol{B_2} \ss B_1 \},
\]
is non-empty. For each $(B_1,B_2) \in \D$, Urysohn's Lemma
guarantees the existence of a \fn $f_{B_1,B_2} \in \cx$ satisfying
\[
\bg{tabular}{ccc}
$f_{B_1,B_2}(X\bs B_1)=\{0\}$, & and & $f_{B_1,B_2}(\ol{B_2}) = \{1\}$.\\
\ed{tabular}
\]
We write
\[
\F=\{f_{B_1,B_2}:(B_1,B_2) \in \D\}.
\]

If $x \neq y$ in $X$, then there clearly exists $(B_1,B_2) \in \D$
\st $x \in B_2,$ and $y \in X \bs B_1$. Since
\[
f_{B_1,B_2}(x)=1 \neq 0 = f_{B_1,B_2}(y),
\]
 $\F$ must separate the points of $X$. Also, it is {\em non-vanishing}; that is, for any $x
\in X$, there must exist an $f \in \F$ \st $f(x) \neq 0$.

Let $\A$ be the smallest \alg that contains $\F \cup \F^*$. It is
clearly a non-vanishing self-adjoint \alg that separates points.
Hence, by the Stone-Weierstrass Theorem, it is dense in $C(X)$.
Clearly, there exists a countable subset of $\A$ that is dense in
$\A$. Take, for example, the smallest \alg over $\bQ$ that contains
$\F \cup \F^*$. This countable subset is then also dense in $C(X)$.
\qed

The following theorem establishes the result in the opposite
direction.

\begin{thm}
If $X$ is a compact Hausdorff topological space, then the
separability of  $C(X)$ implies the metrisability of $X$.
\end{thm}

\demo Since $C(X)$ is separable, it contains a dense sequence of
continuous functions $\{f_{n}\}_{n \in \bN}$. Define $g_n=
\frac{f_{n}}{1+\|f_{n}\|}$, this guarantees that $\|g_n\| < 1$.

If $x \neq y$ in $X$, then by Urysohn's Lemma there exists an $f \in
C(X)$ \st $f(x) \neq f(y)$. Therefore, there must exist a $g_n$ \st
$g_n(x) \neq g_n(y)$, and so $\{g_{n}\}_{n \in \bN}$ separates the
points of $X$.

Let us define
\[
d(x,y)=\sup_{n \in \bN}2^{-n}|g_{n}(x)-g_{n}(y)|.
\]
Since $\{g_{n}\}_{n \in \bN}$ separates the points of $X$,
$d(x,y)=0$ implies that $x=y$. Therefore, $d$ is a metric on $X$.

We shall now examine the open balls of $d$. Let $x \in X,~0 < \e <
1$, and choose $N \in \bN$ \st $2^{-N}<\e$. We write
\[
U_n = g_n \inv \{z \in \bC : |g_n(x)-z|< \e\},
\]
and
\[
U=\bigcap_{n=0}^{N} U_n.
\]
Each $U_n$ is open \wrt the original topology on $X$ and, as a
result, $U$ is also open. Let $y \in U$. If $n \leq N$, then
\[
2^{-n}|g_n(x)-g_n(y)| \leq 2^{-n} \e \leq \e.
\]
If $n > N$, then, since $\|g_n\| < 1$,
\[
2^{-n}|g_n(x)-g_n(y)| < 2^{-n}2  < \e.
\]
It follows that
\[
d(x,y)=\sup_{n \in \bN}2^{-n}|g_n(x)-g_n(y)|< \e.
\]
Thus, $U \subset B_{\e}(x)$, and so the open balls of the metric are
open \wrt the original topology. Let us now consider the identity
map from the compact space $X$, endowed with its original topology,
to $X$, endowed with the metric topology of $d$. Clearly, this map
is a \cts bijection, thus, since the metric topology of $d$ is
Hausdorff, it is a homeomorphism. \qed

Thus, a topological space is metrisable \iff its \alg of \cts \fns
is separable. This gives us the final generalisation of this
section:
\[
\bg{tabular}{lll} metrisable \cpt space& $\lto$ & separable unital
\calgn.
\end{tabular}
\]

\end{enumerate}

Before we finish it should be noted that the transition from
topology to \alg is not always as smooth as in the examples above.
It quite often happens that there is more than one option for the
generalisation of a topological feature (or differential feature),
and it may not always be obvious which one is the `correct' choice;
we will see an example of this when we come to the generalise the de
Rham calculus in Chapter $2$.

\section{Vector Bundles and Projective Modules}

In this section we shall present the Serre--Swan Theorem following
more or less Swan's original proof in \cite{SWAN} (for a more modern
category style proof see \cite{BROD}). This result will give us a
\nc generalisation of vector bundles, one of the basic objects in
differential geometry. Thus, this section can be seen as our first
venture into \ncg proper.

Much of the material presented in this section will be of use to
us later on when we discuss elementary $K$-theory and when we
present the theory of geometric Dirac operators.

\subsection{Vector Bundles}

\begin{defn} A {\em (complex) vector
bundle} is a triple $(E, \pi_E, X)$, consisting of a \tsp $E$
called the {\em total space}, a \tsp $X$ called the {\em base
space}, a \cts surjective map $\pi_E:E \to X$, called the {\em
projection}, and a complex linear space structure defined on each
{\em fibre} $E_x=\pi_E \inv (\{x\})$, such that the following
conditions hold:
\begin{enumerate}
\item For every point $x \in X$, there is an open neighborhood $U$
of $x$, a natural number $n$, and a homeomorphism
\[
\f_U: \pi_E \inv (U) \to U \by \bC^n;
\]
\stn, for all $x \in X$,
\[
\f_U(E_x)=\{x\} \by \bC^n.
\]
\item The map $\f_U$ restricted to $E_x$ is a linear mapping
between $E_x$ and $\{x\} \by \bC ^ n$; (where $\{x\} \by \bC ^ n$ is
regarded as a linear space in the obvious way).
\end{enumerate}
\end{defn}

We could alternatively define a vector bundle to have a real
linear space structure on each fibre. Such a vector bundle is
called a {\em real \vbdn.} All of the results presented below
would hold equally well in this case. However, since we shall
always work with the complex-valued \fns on a topological space,
it suits us better to work with complex \vbdn s.

The canonical examples of a \vbd are the tangent and cotangent
bundles of a manifold.

Whenever possible, we shall denote a vector bundle $(E,\pi_E,X)$
by $E$ and suppress explicit reference to the projection and the
base space. Also, when no confusion arises, we shall use $\pi$
instead $\pi_E$. An open set of the form $U$ is called a {\em base
\nbdn} and the corresponding homeomorphism $\f_U$ is called the
{\em associated local trivialisation}. Note that the definition
implies that the dimension of the fibres is locally constant.

Let $E$ and $F$ be two \vbds over the same base space $X$. A {\em
bundle map} $f$ from $E$ to $F$ is a \cts mapping ${f:E \to F}$,
\st $\pi \circ f=\pi$, and $f_x$, the restriction of $f$ to $E_x$,
is a linear mapping from $E_x$ to $F_x$. If $f$ is also a
homeomorphism between $E$ and $F$, then it is called a {\em bundle
isomorphism}. A vector bundle $(E,\pi, X )$ is called {\em
trivial} if it is isomorphic to the bundle ${(X \by
\bC^n,\pi,X)}$, for some natural number $n$.

Let $U$ be a base \nbd for $E$ and $F$, and let $n$ and $m$ be the
dimensions of $E$ and $F$ respectively. When restricted to $\pi \inv
(U)$, any bundle map $f:E \to F$ induces a map ${\wt{f}:U \by \bC^n
\to U \by \bC^m}$, defined by $\wt{f}=\ps_U \circ f \circ \f_U
\inv$. This in turn induces a map $\wh f:U \to M_{m \by n}(\bC)$,
which is determined by the formula ${\wt{f}(x,v)=(x,\wh f(x)v)}$,
for $(x,v) \in U \by \bC^n$. As a little thought will verify, if
$\wh f$ is \cts for each such \nbd $U$, then $f$ will be \ctsn. An
important use of this fact arises when the bundle map takes each
$E_x$ to the corresponding fibre $F_x$ by a linear isomorphism. In
this case, $m=n$, and so $f \inv$ will determine a mapping $\wh{f
\inv}$ from $U$ to $M_{n \by n}(\bC)$. Now,~ as a little more
reflection will verify, $\wh{f \inv}(x) = (\wh f(x)) \inv$. Thus, $f
\inv$ is \ctsn. This implies that $f$ is a homeomorphism, which in
turn implies that $f$ is a bundle isomorphism. Thus, we have
established the following lemma.

\begin{lem} \label{1:lem:bundleisom}
Let $E$ and $F$ be two \vbds over the same base space $X$. A \cts
bundle map $f:E \to F$ is an isomorphism \iff it maps $E_x$ to
$F_x$ by a linear isomorphism, for all $x \in X$.
\end{lem}

\bigskip

A {\em subbundle} of $(E,\pi_E,X)$ is a \vbd $(S,\pi_S,X)$, \stn:
\begin{enumerate}
\item $S \sseq E$,
\item $\pi_S$ is the restriction of $\pi_E$ to $S$,
\item the linear structure on each fibre $S_x$ is the linear structure
induced on it by
$E_x$.
\end{enumerate}

\subsubsection{Transition Functions}
Let $E$ be a \vbd over $X$, let $\{U_\a\}$ be an open covering of
$X$ by base \nbdsn, and let $\f_\a$ denote the corresponding
trivialisation maps. Now, consider the following diagram:
\begin{displaymath}
\xymatrix{ & \pi \inv (U_{\a} \cap U_{\b})
\ar[dl]_{\f_{\b}}
\ar[dr]^{\f_{\a}} & \\%
(U_{\a} \cap U_{\b}) \by \bC^n \ar[rr]_{\f_\a \circ \f^{-1}_{\b}}&
& (U_{\a} \cap U_{\b}) \by \bC^n.}
\end{displaymath}
Clearly $\f_\a \circ \f^{-1}_{\b}$ is a bundle map from $(U_{\a}
\cap U_{\b}) \by \bC^n$ to itself. As explained above, this means
that it determines a \cts map from $U_\a \cap U_\b$ to $M_n(\bC)$;
we denote this map by $g_{\a\b}$. In fact, since $\f_\b \circ
\f_\a \inv$ is a bundle map that is inverse to $\f_\a \circ
\f^{-1}_{\b}$, it must hold that $g_{\a\b}(x) \in \GL(n,\bC)$, for
all $x \in U_\a \cap U_\b$. We call the collection
\[
\{g_{ \alpha \beta}: {U_\alpha \cap U_\beta \neq \emptyset\}}
\]
the {\em transition \fns} of the vector bundle for the covering
$\{U_\a\}$.

Three observations about any set of transition \fns can be made
immediately;
\begin{enumerate}
\item $g_{\alpha \alpha}=1$,
\item $g_{\beta \alpha} \circ g_{\alpha \beta}=1$,
\item if $U_\alpha \cap U_\beta \cap U_\gamma \neq
\emptyset$, then $g_{\alpha \beta}\circ g_{\beta \gamma} \circ
g_{\gamma \alpha}=1.$
\end{enumerate}
The last property is known as the {\em cocycle condition}.

The following result is of great importance in the theory of
vector bundles, and shows that the transition \fns for any
covering completely determine the bundle.

\begin{prop} \label{1:prop:transfns}
Given a cover $\{U_\alpha\}$ of $X$, and a \cts map
\[
g_{\alpha \beta}:U_\alpha \cap U_\beta \to \GL(n,\bC),
\]
for every non-empty intersection $U_\a \cap U_\b$, such that the
conditions $1,2,3$ listed above hold, then there exists a vector
bundle (unique up to bundle isomorphism) for which $\{g_{\alpha
\beta}\}$ are the transition functions.
\end{prop}

\subsubsection{Sections}
Let $E$ be a \vbd over a base space $X$. A \cts mapping $s$ from
$X$ to $E$ is called a {\em section} if $\pi (s(x)) = x$, for all
$ x \in X$. The set of sections of $E$ is denoted by $\G(E)$. If
$s,t \in \G(E)$ and $a \in C(X)$, we define
\[
(s+t)(x)=s(x)+t(x), \text{~~~ and ~~~} (sa)(x)=a(x)s(x).
\]
Here it is understood that the addition and scalar multiplication
on the right hand side of each equality takes place in $E_x$. The
mappings $s+t$ and $sa$ are \cts since the composition of either
with a local trivialisation is \cts on the corresponding base
\nbdn. Hence, the mappings $s+t$ and $sa$ are sections. With
respect to these definitions $\G(E)$ becomes a
right-$C(X)$-module. Let $G$ be another \vbd over $X$ and let $f:E
\to G$ be a bundle map. We define $\G(f)$ to be the unique module
mapping from $\G(E)$ to $\G(G)$ \st
\[
[\G(f)(s)](x)=(f \circ s)(x).
\]
Let $U$ be a base \nbd of $x \in X$ and let $\f_U$ be its
associated trivialisation. Consider the set of $n$ \cts mappings
on $U$
\[
\wh{s_i}:U \to U \by \bC^n, ~~~~ x \mto (x,e_i),
\]
where $\{e_i\}_{i=1}^{n}$ is the standard basis of $\bC^n$. If we
define $s_i = \f_U \inv \circ \wh{s_i}$, then the set
$\{s_i(y)\}_{i=1}^{n}$ forms a basis for $E_y$, for every $y \in U$.
If we assume that $X$ is normal, then using the Tietze Extension
Theorem, it can be shown that for each $s_i$, there exists a section
$s'_i \in \G(E)$ \st $s_i$ and $s'_i$ agree on some \nbd of $x$. We
call the set $\{s'_i(x)\}_{i=1}^{n}$ a {\em local base } at $x$.

\subsubsection{Direct sum}

Let $E$ and $F$ be two vector bundles and let $\{U_\a\}$ be the
family of subsets of $X$ that are base \nbds for both bundles. Now,~
if ${U_\a \cap U_\b}$ is a non-empty intersection, then we denote
the corresponding transition \fns for $E$ and $F$ by $g^E_{\a\b}$
and $g^F_{\a\b}$ respectively. Using $g^E_{\a\b}$ and $g^F_{\a\b}$
we can define a matrix-valued \fn on $U_\a \cap U_\b$ by
\[
g_{\a\b}:x \mto \left(%
\begin{array}{cc}
  g^E_{\a\b}(x) & 0 \\
  0          & g^F_{\a\b}(x) \\
\end{array}%
\right).
\]
Clearly, $g_{\a\b}$ satisfies the conditions required by
Proposition (\ref{1:prop:transfns}), for every non-empty
intersection $U_\a \cap U_\b$. Hence, there exists a vector bundle
for which $\{g_{\a\b}\}$ is the set of transition \fnn s. The
nature of the construction of each $g_{\a\b}$ implies that the
fibre of the bundle over any point $x \in X$ will be isomorphic to
$E_x \oplus F_x$. This prompts us to denote the bundle by $E
\oplus F$ and to call it the {\em direct sum} of $E$ and $F$.

Using an analogous construction, we can produce bundles whose
fibres over any $x \in X$ are equal to $E^*_x$, $\text{Hom}(E_x)$,
or $E_x \oby F_x$; we denote these bundles by $E^*$,
$\text{Hom}(E)$, and $E \oby F$ respectively.

\subsubsection{Inner Products and Projections}

An {\em inner product} on a \vbd $E$ is a \cts mapping from
\linebreak ${D(E)=\{(v_1,v_2):v_1,v_2 \in E,\, \pi(v_1)=\pi(v_2)\}}$
to $\bC$ \st its restriction to \linebreak ${E_x \by E_x}$ is an
inner product, for all $x \in X$. Using locally defined inner
products and a partition of unity, an inner product can be defined
on any \vbdn.

Let $S$ be a subbundle of $E$, and let $\< \cdot , \cdot \>$ be an
inner product on $E$. Using $\< \cdot , \cdot \>$ we can define an
orthogonal projection $P_x:E_x \to S_x$, for each $x \in X$. This
defines a map $P:E \to S$, which we shall call the {\em projection}
of $E$ onto $S$. To see that this map is \cts we shall examine the
mapping it induces on $U \by \bC^n$, for some arbitrary base \nbd
$U$. Let $\{t_i\}_{i=1}^n$ be a {\em local basis} over $U$, that is,
let $\{t_i\}_{i=1}^n$ be a set of sections \stn, for each $p \in X$,
$\{t_i(p)\}_{i=1}^n$ is a basis of $E_x$; then consider the map
\[
U \by \bC^n \to D,~~~~(x,v) \mto (\f_U\inv (x,v),t_i(x)).
\]
It is clearly \ctsn, for all $i = 1, \ldots ,n$, and as a result its
composition with the inner product is \ctsn; that is, the map $(x,v)
\mto \<\f_U\inv (x,v),t_i(x)\>$ is \ctsn. Therefore, the induced map
\[
U \by \bC^n \to \bC,~~~~~~(x,v) \mto \sum_{i=1}^n
\<(x,v),t_i(x)\>t_i(x)
\]
is \ctsn. The \cty of the projection easily follows.

\subsection{Standard results}

The following four results are standard facts in vector bundle
theory. As above, their proofs consist of routine arguments
involving local \nbds and local sections. We shall briefly sketch
how the results are established, and  refer the interested reader
to \cite{SWAN}, \cite{ATIY}, or \cite{HUS}.

\begin{lem} \label{1:lem:linindpt} Let $s_1, \ldots ,s_k$ be sections of a
\vbd $E$ such that for some $x \in X$, $s_1(x), \ldots ,s_k(x)$ are
linearly independent in $E_x$. Then there is a \nbd $V$ of $x$ \st
$s_1(y), \ldots ,s_k(y)$ are linearly independent, for all $y \in
V$.
\end{lem}

In fact, this lemma is a simple consequence of the \cty of the
determinant \fnn. It has the following easy corollary.

\begin{cor} \label{1:cor:linindpt}
Let $F$ be a \vbd over $X$ and let $f:E \to F$ be a bundle map. If
${\dim(\im (f_x)) = n}$, then $\dim(\im(f_y)) \geq n$, for all $y$
in some \nbd of $x$.
\end{cor}

Note that in the following theorem $\ker(f)$, the {\em kernel} of a
bundle map $f$, is the topological subspace $\bigcup_{x \in X}
\ker(f_x)$. We define $\im (f)$ similarly.

\begin{thm} \label{thm:subbundle}
Let $f:E \to F$ be a bundle map. Then the following statements are
equivalent:
\begin{enumerate}
\item $\im (f)$ is a subbundle of $F$;%
\item $\ker(f)$ is a subbundle of $E$;%
\item the dimensions of the fibres of $\im (f)$ are locally constant; %
\item the dimensions of the fibres of $\ker(f)$ are locally
constant.
\end{enumerate}
\end{thm}

(Note that the linear structure on each $\im(f_x)$ and $\ker(f_x)$
is understood to be that induced by $F_x$ and $E_x$ respectively.)

Clearly statements $(3)$ and $(4)$ are equivalent and are implied by
either statement $(1)$ or $(2)$. Thus, the theorem would be proved
if $(3)$ could be shown to imply $(1)$ and $(4)$ could be shown to
imply $(2)$. Now,~if one assumes local constancy of the fibres, then
Lemma \ref{1:lem:linindpt} can be used to construct local bases for
$\im (f)$ and $\ker (f)$ at each point of $X$. A little thought will
verify that the existence of a local base implies local triviality,
and the result follows.

\begin{thm} \label{thm:functorf}
If $X$ is a compact Hausdorff space, then for any module
homomorphism $G:\G(E) \to \G(F)$, there is a unique bundle map $g: E
\to F$ such that $G= \G(g)$.
\end{thm}

The first step in establishing this theorem is to show that
$\G(E)/I_x \simeq F_x$, for all $x \in X$, where $I_x=\{f \in
C(X):f(x)=0\}$. Since $F$ clearly induces a map from $\G(E)/I_x$ to
$\G(F)/I_x$, it must now induce a map from $E$ to $F$. This map can
then be shown to satisfy the required properties and the result
follows.

\subsection{Finite Projective Modules}
Recall that a right $A$-module $\E$ is said to be {\em projective}
if it is a direct summand of a free module; that is, if there
exists a free module $\F$ and a module $\E'$, such that
\[
\F = \E \oplus \E '.
\]
Recall also that this is equivalent to the following alternative
definition: A right $A$-module $\E$ is {\em projective} if, given a
surjective homomorphism $\t : \M \lto \N$ of right $A$-modules, and
a homomorphism $\l : \E \lto \N$, there exists a homomorphism
$\wt{\l} : \E \lto \M$ such that $\t \circ \wt{\l} = \l$, or
equivalently, \st the following diagram is commutative
\begin{displaymath}
\xymatrix{  &
\E \ar[dl]_{\wt{\l}} \ar[d]^{\l}\\
\M \ar[r]_{\t} & \N. }
\end{displaymath}

Suppose now that $\E$ is a projective and finitely-generated module
over $A$. Clearly, there exists a surjective homomorphism $\t : A^n
\to \E$, for some natural number $n$. Since $\id_{\E}$ is a
homomorphism from $\E$ to $\E$, the projective properties allow us
to find a map $\wt{\l} : \E \to A^n$ such that $\t \circ \wt{\l} =
\id_{\E}$, or equivalently, \st that the following diagram is
commutative

\begin{displaymath}
\xymatrix{  &
\E \ar[dl]_{\wt{\l}} \ar[d]^{id_{\E}}\\
A^n \ar[r]_{\t} & \E. }
\end{displaymath}
\label{profin}

We then have an idempotent element $p$ of $\End A^n$, given by
\[
p = \wt{\l} \circ \t. \label{projector}
\]
We can see that $p$ is idempotent from
\[
p^2 = \wt{\l} \circ \t \circ \wt{\l} \circ \t = \wt{\l} \circ \t =
p.
\]
This  allows one to decompose the free module $A^n$, in the standard
manner, as a direct sum of submodules,
\[
A^n = \im (p) \oplus \ker(p) = p A^n \oplus (1-p) A^n.
\]
Now,~since $\wt{\l} \circ \t=\id_{pA^n}$ and $\t \circ \wt{\l} =
\id_{\E}$, we have that $\E$ and $pA^n$ are isomorphic as right
$A$-modules. Thus, a module $\E$ over $A$ is finitely-generated and
projective if, and only if, there exists an idempotent $p \in
\End{A^n}$ such that $\E= pA^n$.

\subsection{Serre--Swan Theorem}

As a precursor to the Serre--Swan Theorem, we shall show that if
$X$ is a compact Hausdorff space, and $E$ is a \vbd over $X$, then
$\G(E)$ is a  finitely-generated projective right-$C(X)$-module.

Let $E$ be a vector bundle over $X$, and let $S$ be a subbundle of
$E$. We endow $E$ with an inner product $\<\cdot,\cdot\>$, and we
 denote the projection of $E$ onto $S$ by $P$. If $S_x^\bot$
is the subspace of $E_x$ that is orthogonal to $S_x$, then $S^\bot =
\bigcup_{x \in X}S_x^\bot$ is the kernel of $P$. Since the image of
$P$ is $S$, Theorem \ref{thm:subbundle} implies that $S^\bot$ is a
subbundle of $E$. Now,~the mapping
\[
S \oplus S^\bot \to E,~~~~(v,w) \mto v+w,
\]
is clearly a \cts mapping. Moreover, it is an isomorphism on each
fibre. Thus,  Lemma \ref{1:lem:bundleisom} implies that the two
spaces are isomorphic. We now summarise what we have established
in the following lemma.

\begin{lem} \label{1:lem:orthosubbundle}
Let $E$ be a \vbd equipped with an inner product, and let $S$ be a
subbundle of $E$. If $S_x^\bot$ is the subspace of $E_x$ that is
orthogonal to $S_x$ and $S^\bot=\bigcup_{x \in X} S_x$, then
$S^\bot$ is a subbundle of $E$ and $E \simeq S^\bot \oplus S$.
\end{lem}

\begin{lem} Let $X$ be a compact Hausdorff space and let
$E$ be a vector bundle. Then there exists a surjective bundle map
$f$ from the trivial bundle ${X \by \bC^n}$ to $E$, for some \pve
integer $n$.
\end{lem}
\demo For each $x \in X$, choose a set of sections $s_1^x, \ldots
,s_{k_x}^x \in \G(E)$ that form a local base over $U_x$, a base \nbd
of $x$. A finite number of these \nbds cover $X$. Therefore, there
are a finite number of sections $s_1, \ldots ,s_n \in \G(E)$ \st
$s_1(x), \ldots ,s_n(x)$ span $E_x$, for every $x \in X$.

Now,~$\G(X \by \bC^n)$ is a free module over $C(X)$ generated by
sections  $t_i(x)=(x,e_i)$, where $\{e_i\}_{i=1}^n$ is the standard
basis of $\bC^n$. There is a unique module map from $\G(X \by
\bC^n)$ to $\G(E)$  that maps each $t_i$ onto $s_i$. By Theorem
\ref{thm:functorf} this mapping is induced by a map $f:X \by \bC^n
\to E$. Since,
\[
f(t_i(x))=[\G(f)(t_i)](x)=s_i (x),~~~~\text{for all } x \in X,
\]
it is clear that $f$ is surjective. \qed

\begin{cor} \label{cor:onlyif} If $X$ is a compact Hausdorff space,
then any vector bundle $E$ is a direct summand of a trivial
bundle, and $\G(E)$ is a finitely-generated projective
right-$C(X)$-module.
\end{cor}

\cvtu Let $f:X \by \bC^n \to E$ be the map defined in the previous
lemma. Since $\im (f)=E$, Theorem \ref{thm:subbundle} implies that
$\ker(f)$ is a subbundle of ${X \by \bC^n}$. If we put an inner
product on $X \by \bC^n$, then by Lemma
\ref{1:lem:orthosubbundle},\\
\[
\ker(f) \oplus \ker (f)^\bot \simeq X \by \bC^n.
\]
Restricting $f$ to $\ker(f)^\bot$, we see that it is a linear
isomorphism on each fibre, and so $\ker(f)^\bot \simeq E$.

We can identify $\G(\ker(f) \oplus F)$ and $\G(\ker(f)) \oplus
\G(F)$, using the module isomorphism
\[
\G(\ker(f)) \oplus \G(F) \to \G(\ker(f) \oplus F),~~~~s \oplus t
\mto (s,t);
\]
where $(s,t)(x)=(s(x),t(x))$. Hence, we have that $\G(\ker(f))$ is a
direct summand of the finitely generated free $C(X)$-module $\G(X
\by \bC^n)$. \qed

Finally, we are now in a position to prove the principal result of
this section, the Serre--Swan Theorem. It was first published in
$1962$ \cite{SWAN}, and was inspired by a paper of Serre \cite{SER}
that established an analogous result for algebraic vector bundles
over affine varieties.

\begin{thm}[Serre--Swan] Let $X$ be a compact Hausdorff space. Then
a module $\E$ over $C(X)$ is isomorphic to a module of the form
$\G(E)$ \iff $\E$ is finitely generated and projective.
\end{thm}

\demo If $\E$ is finitely-generated and projective, then, as
explained earlier, there exists an idempotent endomorphism
$p:C(X)^n \to C(X)^n$, with $\E \simeq \im (p)$, for some natural
number $n$. Clearly $C(X)^n$ can be associated with the sections
of the trivial vector bundle $X \by \bC^n$ by mapping $(f_1,
\ldots ,f_n)$ to the section $s$, defined by $s(x)=(x,f_1(x),
\ldots ,f_n(x))$.

By Theorem \ref{thm:functorf}, $p$ is the image under $\G$ of a
bundle map $f: X \by \bC^n \to X \by \bC^n$. Since $p^2=p$, and
since $p(s)= f \circ s$, we have that
\[
f^2 \circ s=p(f \circ s)=p^2(s)=p(s)=f \circ s ,
\]
for all sections $s$. Thus, since $\{s_i(x)\}_{i=1}^n$ spans $(X \by
\bC^n)_x$, for all $x \in X$, it holds that $f^2=f$.

Let us now define the map
\[
(1-f):X \by \bC^n \to X \by \bC^n,~~~~ v \mto v-f(v),
\]
where the addition takes place fibrewise. As usual we denote the
restriction of $f$ to the fibre $(X \by \bC^n)_x$ by $f_x$. Since
$f_x$ is an idempotent linear map, it holds that
$\im(1-f_x)=\ker(f_x)$. Clearly this implies that $\ker (f)= \im
(1-f)$.

Suppose that $\dim (\im (f_x)) = h$ and that $\dim (\ker (f_x))=k$,
then Lemma \ref{1:lem:linindpt} implies that $\dim (\im (f_y)) \geq
h$, and $\dim (\ker (f_y)) \geq k$, for all $y$ in some \nbd of $x$.
However, since $(X \by \bC^n)_x=\im (f_x) \oplus \ker (f_x)$, for
all $x \in X$,
\[
\dim (\im (f_y)) +  \dim (\ker (f_y)) =  \dim (X \by \bC^n)_y = h+k
\]
is a constant. Thus, $\dim (\im (f_y))$ must be locally constant.
This implies that $\im(f)$ is a subbundle of $X \by \bC^n$.

If we make the observation that
\[
\G(\im(f))=\{f \circ s:s \in \G(X \by \bC^n)\},
\]
then we can see that
\[
\G(\im(f))=\im(\G(f)) =\im(p)=\E.
\]
Thus, $\E$ is indeed isomorphic to the sections of a \vbdn.

The proof in the other direction follows from Corollary
\ref{cor:onlyif}. \qed

Thus, the Serre--Swan Theorem shows that the \vbds over \cpt
Hausdorff space $X$ are in one-to-one correspondence with the finite
projective modules over $C(X)$. This motivates us to view
finitely-generated projective modules over \nc \calgs as
non-commutative \vbdsn.

\subsubsection{Smooth Vector Bundles}

Let $E$ be a \vbd over a manifold $X$. It is not too hard to see
that we can use the differential structure of $X$ to canonically
endow $E$ with a differential structure. A routine check will
establish that the local trivialisations of the bundle then become
smooth maps. In general, if $(E,\pi,X)$ is a \vbd \st $E$ and $X$
are manifolds and all the local trivialisations are smooth, then we
call $(E,\pi,X)$ a {\em smooth \vbdn.} {\em Smooth vector bundle
maps} and {\em smooth vector bundle isomorphisms} are defined in the
obvious way. When we speak of the {\em smooth sections of an
ordinary \vbd} we mean the smooth sections of the bundle endowed
with the canonical differential structure discussed above. Clearly,
the tangent and cotangent bundles of a manifold are smooth \vbdn s.
An important point to note is that a direct analogue of Proposition
\ref{1:prop:transfns} holds for smooth \vbdn s.

\bigskip

When a section of $E$ is also a smooth map between $X$ and $E$, then
we call it a {\em smooth section}; we denote the set of smooth
sections by $\G^\infty(E)$. The canonical example of a smooth
section is a smooth vector field over a manifold; it is a smooth
section of the tangent bundle. Now,~just as we gave $\G(E)$ the
structure of a right module over $C(X)$, we can give $\G^\infty(E)$
the structure of a right module over $C^\infty(X)$, where
$C^\infty(X)$ is the \alg of smooth complex-valued \fns on $X$.

It can  be shown, using an argument quite similar to the one
above, that the modules of smooth sections of the smooth vector
bundles over $X$ are in one-to-one correspondence with the
finitely-generated projective modules over $C^{\infty}(X)$; for
details see \cite{VAR}.

\section{Von Neumann Algebras}

In this section we give a brief presentation of \nc measure
theory. We shall not venture too far into the details since, with
the exception of the material presented on \lcsn s, we shall not
return to this area again. It is  introduced here for its
heuristic value only. We refer the interested reader to \cite{MUR}
for the details of von Neumann \algsn, and to \cite{CON} for an
in-depth presentation of \nc measure theory.

\bigskip

The parallels between \nc measure theory and \nc topology are
obvious. In fact, both areas are part of an overall trend in
mathematics towards viewing \fn \algs as special commutative cases
of operator \algn s. This area is loosely known as {\em quantum
mathematics}, an obvious reference to the quantum mechanical origins
of operator theory. An excellent overview of this trend towards
`quantization' can be found in \cite{WEAV}.

\subsubsection{Locally Convex Topological Vector Spaces}

Let $\cP$ be a non-empty family of seminorms on a linear space $X$.
If $x \in V$, $\e \geq 0$, and  $\cP_0$ is a finite subfamily of
$\cP$, then we define
\[
B(x,\cP_0,\e)=\{y \in X: p(x-y) \leq \e,~ p \in P_0\}.
\]
It straightforward to show that
\[
\B=\{B(x,\cP_0,\e):x \in X,~ \cP_0 \text{ a finite subfamily of }
\cP,~ \e \geq 0\}
\]
is a base for a topology on $X$; it is called the {\em topology
generated by} $\cP$. It is clear that a net $(x_\l)_\l$ in $X$
converges to $x$, \wrt this topology, \iff $p(x-x_\l) \to 0$, for
all $p \in \cP$. If $x_\l \to x$ and $y_\l \to y$, then it is easily
seen that $p(x_\l + y_\l - x -y) \to 0$, for all $p \in \cP$. Thus,
addition is \cts \wrt the topology generated by $\cP$. Similarly,
scalar multiplication can be shown to be \ctsn. Hence, when $X$ is
endowed with this topology it is a \tvsn. (Recall that a {\em
topological vector space} is a linear space for which the linear
space operations of addition and scalar multiplication are \ctsn.)
It can be shown that the topology is Hausdorff \iff for each $x \in
X$, there exists a $p \in \cP$ \st $p(x) \neq 0$.

We call a \tvs whose topology is determined by a family of seminorms
a {\em locally convex (topological vector) space} ({\em locally
convex} refers to the fact that $\B$ forms a locally convex base for
the topology).

The simplest example of a \lcs is a normed vector space. The
family of seminorms is just the one element set containing the
norm, and the topology generated is the norm topology.

\subsubsection{Von Neumann Algebras}

If $H$ is a Hilbert space then the {\em strong operator} topology
is the topology generated by the family of seminorms $\{\| \cdot
\|_x: x \in H\}$, where
\[
\|T\|_x=\|Tx\|, ~~~~~T \in B(H).
\]
We denote the strong operator topology by $\t_S$. Thus, $T_{\l}
\to T$ \wrt $\t_S$ \iff ${\|(T_{\l}-T)x\| \to 0}$, for all $x \in
H$.

Recall that if $\t_1, \t_2$ are two topologies on a set $X$ \st \cvg
of a net \wrt $\t_1$ implies \cvg of the net \wrt $\t_2$, then $\t_2
\sseq \t_1$. Suppose now that $A_\l \to A$ in $B(H)$ \wrt the norm
topology. Since  $\|(A_\l-A)x\| \leq \|A_\l - A\|\|x\|$, for all $x
\in X$, we have that $A_\l \to A$ \wrt $\t_S$. Thus, the strong
operator topology is weaker than the norm topology. In fact, it can
be strictly weaker. This happens \iff $H$ is infinite-dimensional.

\begin{defn}
Let $H$ be a Hilbert space. If $\A$ is a $*$-sub\alg of $B(H)$
that is closed \wrt $\t_S$, then we call $\A$ a {\em von Neumann
\algn}.
\end{defn}

Since every subset that is closed \wrt the strong operator
topology is also closed \wrt the norm topology, every von Neumann
\alg is also a \calgn. In fact, it can be shown that every von
Neumann \alg is a unital \calgn.

If $A$ is an \alg and $C$ is a subset of $A$, then we define $C'$,
the {\em commutant} of $C$, to be the set of elements of $A$ that
commute with all the elements of $C$. We define $C''$ to be $(C')'$
and call it the {\em double commutant} of $C$. Clearly $C \sseq
C''$. The following famous theorem shows us that there exists a very
important relationship between the double commutants and von Neumann
\algn s; see chapter $4$ of \cite{MUR} for a proof.

\begin{thm}[von Neumann] \label{1:thm:doublecommutant}
If $\A$ is a $*$-sub\alg of $B(H)$, for some Hilbert space $H$, \st
$\id_{H} \in \A$, then  $\A$ is a von Neumann \alg \iff $\A''=\A$.
\end{thm}

\subsection{Noncommutative Measure Theory} \label{1:sect:NCMT}
\subsubsection{$L^{\infty}$ as a von Neumann Algebra}

Recall that if $(M,\mu)$ is a measure space, then
$L^{\infty}(M,\mu)$ is the \alg of all equivalence classes of
measurable \fns on $M$ that are bounded almost everywhere (two \fns
being equivalent if they are equal almost everywhere). If it is
equipped with the {\em $L^\infty$-norm} defined by
\[
\|f\|_{\infty}=\inf\{C \geq 0:|f|\leq C ~~\ae\},
\]
then it is a unital Banach \algn. We can define an involution on
$L^{\infty}(M,\mu)$ by $f^*=\ol{f}$, and this clearly gives it the
structure of a \calgn.

Recall also that $L^2(M,\mu)$ is the \alg of equivalence classes of
measurable \fns on $M$  \st if $f \in L^2(M,\mu)$, then $|f|^2 $ has
finite integral. If it is equipped with the {\em $L^2$-norm},
defined by setting
\[
\|f\|_{2}=\left(\int_M |f|^2d\mu \right)^{\frac{1}{2}},
\]
then it is a Banach space. In fact, $L^2(M,\mu)$ is a Hilbert
space since the norm is generated by the inner product
$\<f,g\>=\int_Mf\ol g d\mu$.

Let us now assume that $M$ is a \cpt Hausdorff space and that
$\mu$ is a finite positive regular Borel measure. For any $f \in
L^{\infty}(M,\mu)$, consider the mapping
\[
M_f:L^2(M,\mu) \to L^2(M,\mu),~~~~g \mto fg.
\]
Clearly each $M_f$ is a linear mapping, and since
\[
\|M_f g\|_2^2=\int_M |fg|^2 d\mu \leq \|f\|_{\infty}^2 \int_M |g|^2
d\mu = \|f\|_{\infty}^2\|g\|_2^2,
\]
each $M_f$ is bounded by $\|f\|_{\infty}$. In fact, using the
regularity of the measure, we can show that
$\|M_f\|=\|f\|_{\infty}$. The adjoint of $M_f$ is equal to $M_{\ol{
f}}$, and so the mapping
\[
M:L^{\infty}(M,\mu) \to B(L^2(M,\mu)),~~f \mto M_f
\]
is an isometric  $*$-isomorphism between $L^\infty (M ,\mu)$ and
$\L=M(L^\infty (M ,\mu)).$ Moreover, it can be shown that if $T \in
\L'$, then $T = M_f$, for some $f \in L^\infty(M,\mu)$. Thus, by
Theorem~\ref{1:thm:doublecommutant}, $L^\infty (M ,\mu)$ is an
abelian von Neumann \algn.

\subsubsection{General Abelian von Neumann Algebras}

The following theorem shows that all abelian von Neumann are of the
form $L^\infty (M ,\mu)$, for some measure space $(M ,\mu)$. We can
consider this result the analogue of the Gelfand--Naimark Theorem
for von Neumann \algsn; for a proof see \cite{DIX}.

\begin{thm}
Let $\A$ be an abelian von Neumann \alg on a Hilbert space $H$. Then
there exists a locally \cpt Hausdorff space $M$ and a positive Borel
measure $\mu$ on $M$ \st $\A$ is isometrically $*$-isomorphic to
$L^\infty(M,\mu)$.
\end{thm}

It is interesting to note that the space $\M$ is produced as the
spectrum of a \mbox{$C^*$-sub\alg} of $\A$ that is dense in $\A$
\wrt the strong operator topology. Furthermore, if $H$ is separable,
then $M$ can be shown to be \cpt and second-countable.

\subsubsection{Noncommutative Measure Theory}

In the same spirit as \nc topology, we now think of \nc von Neumann
\algs as `\nc measure spaces'. A lot of work has been put into
finding von Neumann algebra structures that correspond to measure
theoretic structures in the commutative case. As would be expected,
this area of mathematics is called {\em \nc measure theory}.

Central to most of this work is the notion of a factor. A {\em
factor} is a von Neumann algebra with a trivial centre, that is, a
von Neumann \alg $\A$ for which ${\A' \cap \A= \bC 1}$. The simplest
example is $B(H)$, for any Hilbert space $H$. Factors have been
classified into three types according to the algebraic properties of
the projections they contain. The important thing about them is that
every von Neumann \alg is isomorphic to a direct integral of
factors. (A {\em direct integral} of linear spaces is a \cts
analogue of the direct sum of  linear spaces.)

Much of Connes' original work was in this area. Building on the
Tomita--Takesaki Theorem, he established a \nc version of the
Radon-Nikodym Theorem. This result furnishes a canonical
homomorphism from the additive group $\bR$ to the group of outer
automorphisms of any \nc von Neumann \algn. It has no parallel in
the commutative case and it inspired Connes to write that `\nc
measure spaces evolve with time'. This work then led on to a
classification of all hyperfinite type III factors (type III being
one of the three types of factors).

Connes has found applications for his results in the study of the
type of singular spaces discussed in the introduction. He has met
with particular success in the study of foliations of manifolds.

\section{Summary} We conclude this chapter by summarising the
algebraic generalisations of the elements of topology,
differential geometry, and measure theory collected above:
\[
\bg{tabular}{lll}
locally compact space              & $\lto$ &  $C^{*}$-algebra, \\
compact space                      & $\lto$ &  unital $C^{*}$-algebra, \\
homeomorphism                      & $\lto$ &  $*$-isomorphism,\\
image of a \fn                     & $\lto$ &  spectrum of a element,\\
\pve \fn                           & $\lto$ &  \pve element,\\
regular Borel complex measure      & $\lto$ &  bounded linear functionals on
a \calgn,\\
one-point compactification of a space        & $\lto$ &  unitisation of a \calgn,\\
closed subset of a \cpt space      & $\lto$ &  quotient of a unital \calgn,\\
open subset of a \cpt space        & $\lto$ &  ideal of a unital \calgn,\\
connected \cpt space               & $\lto$ &  projectionless unital \calgn,\\
metrisable \cpt space              & $\lto$ &  separable unital \calgn,\\
vector bundle                      & $\lto$ &  finite projective module\\
over a locally \cpt space          &        &  over a \calgn,\\
measure space                      & $\lto$ &  von Neumann algebra.\\
\end{tabular}
\]

\begin{rem}
Most of the results of this chapter can be expressed in a very
satisfactory manner using the langauge of category theory. For
example, most of the above correspondences can be viewed as
functors between the respective categories. For a presentation of
this approach see \cite{VAR}.
\end{rem}

\chapter{Differential Calculi}

We are now ready to consider generalised differential structures on
quantum spaces. Following on from the last chapter, we shall begin
with a \cpt manifold $M$, and then attempt to express its structure
in terms of the \alg of its \cts \fnn s.

A natural starting point would be to try and establish an algebraic
relationship between the smooth and the \cts functions of $M$.
Unfortunately, however, there does not appear to be any simple way
of doing this. In fact, there does not seem to be any simple
algebraic properties that characterise $\csm$ at all. The best that
we can do is establish that $\csm$ is dense in $C(M)$ using the
Stone--Weierstrass Theorem. Consequently, we shall use an arbitrary
associative (involutive) algebra $A$ to generalise $\csm$. We could
assume that $A$ is dense in some \mbox{$C^{*}$-algebra}, but there
is no major technical advantage in doing so. Neither is there any
advantage in assuming that $A$ is unital.

\bigskip

We shall begin this chapter by reviewing differential calculus on a
manifold in global algebraic terms. This will naturally lead us to
the definition of a differential calculus: this object is a
generalisation of the notion of the de Rham calculus, and it is of
fundamental importance in noncommutative geometry.

\bigskip

We shall then introduce derivation-based differential calculi, as
formulated by M. Dubois--Violette and J. Madore \cite{DBVMAD1}.
Their work has been strongly influenced by that of J.L. Koszul who
in \cite{KOZ} described a powerful algebraic version of differential
geometry in terms of a general commutative associative algebra. We
present derivation-based calculi here because we shall refer to them
in our discussion of fuzzy physics in Chapter $5$, and because they
provide a pleasingly straightforward example of a noncommutative
differential calculus.

\section{The de Rham Calculus} Let $M$ be an $n$-dimensional
manifold, and let $\csm$ denote the \alg of smooth complex-valued
functions on $M$. Unless otherwise stated, we shall always assume
that the manifolds  we are dealing with are smooth, real, compact,
and without boundary. We define a {\em smooth vector field} $X$ on
$M$ to be a {\em derivation} on $\csm$; that is, a linear mapping on
$\csm$ \stn, for all $f,g \in \csm$,
\[
X(fg)=X(f)g+fX(g).
\]
The set of smooth vector fields on $M$ is denoted by $\X(M)$. We
give it the structure of a (left) $\csm$-module in the obvious way.
We denote the dual module of $\X(M)$ by $\Om^1(M)$, and we call it
the {\em module of differential $1$-forms} over $M$. It is easily
shown that these definitions are equivalent to defining a smooth
vector field to be a smooth section of the tangent bundle of $M$,
and a differential $1$-form to be a smooth section of the cotangent
bundle of $M$.

We denote the $\csm$-module of $\csm$-valued, $\csm$-multilinear
mappings on
\[
\Om^{\by p}(M)\by \X^{\by q}(M)=\underbrace{\Om^1(M) \by \ldots \by
\Om^1(M)}_{p \textrm{ times}}\by \underbrace{\X(M) \by \ldots \by
\X(M)}_{q \textrm{ times}}
\]
by $\fT_q^p(M)$, and we call it the {\em space of rank-$(p,q)$
smooth tensor fields} on $M$. Note that $\fT^0_1(M)=\Om^1(M)$.

An important point about tensor fields is their `locality'. Let
$A_1$ and $A_2$ be two elements of $\Om^{\by p}(M)\by \X^{\by q}(M)$
\st $A_1(p)=A_2(p)$, for some $p \in M$. If $T$ is a rank-$(p,q)$
tensor field, then, using a simple bump-\fn argument, it can be
shown that $T(A_1)(p)=T(A_2)(p)$.

\bigskip

Let $\E$ be a module over $R$, and let $S$ be a multilinear mapping
on $\E^n$ with values in $R$. We say that $S$ is {\em
anti-symmetric} if, for every $e_1, \ldots ,e_n \in \E$,
\[
S(e_1, \ldots ,e_i, \ldots ,e_j, \ldots ,e_n)=-S(e_1, \ldots ,e_j,
\ldots ,e_i, \ldots ,e_n),
\]
for all $1 \leq i <j\leq n$. Consider the submodule of $\fT^0_p(M)$
consisting of the anti-symmetric $p$-linear maps on $\X^1(M)$. We
call it the module of {\em smooth exterior differential $p$-forms}
on $M$, or, more simply, the module of {\em $p$-forms} on $M$; we
denote it by $\Om^p(M)$. Again, it is easily shown that this is
equivalent to defining a $p$-form to be a smooth section of
$p^\mathrm{ th}$-exterior power of the cotangent bundle of $M$. This
equivalent definition implies that $\Om^p(M)=\{0\}$, for $p > n$.

Let us define
\[
\Om(M)=\bigoplus_{p=0}^{\infty} \Om^p(M),
\]
with $\Om^0(M)=\csm$. Then consider the associative bilinear mapping
\[
\LL:\Om(M) \by \Om(M) \to \Om(M), \qquad (\w,\w') \mto \w \wed \w';
\]
where if $\w_p \in \Om^p(M)$, and $\w_q \in \Om^q(M)$, then
\begin{eqnarray}\label{2:eqb:wedgeprod}
\lefteqn{ \w_p \wed \w_q(X_{\pi(1)}, \dots X_{\pi(k)})}
\nonumber\\
&  &  = \sum_{\pi \in \text{Perm}(p+q)}\sgn(\pi)\w_p(X_{\pi(1)},
\dots X_{\pi(k)})\w_q(X_{\pi(p+1)}, \dots X_{\pi(p+q)})
\end{eqnarray}
(Note that $\sum_{\pi \in \text{Perm}(p+q)}$ means summation over
all permutations of the numbers $1,2, \ldots ,p+q$, and $\sgn(\pi)$
is the sign of the permutation $\pi$.) When $\Om(M)$ is endowed with
$\LL$ we call it the {\em \alg of exterior differential forms} on
$M$. Upon examination we see that $\LL$ is {\em graded commutative};
that is, if $\w_p \in \Om^p(M)$ and $\w_q \in \Om^q(M)$, then
\[
\w_q \wed \w_p = (-1)^{pq}\w_p \wed \w_q.
\]

\subsection{The Exterior Derivative}

Consider the canonical mapping from the smooth \fns of $M$ to the
smooth $1$-forms of $M$ defined by
\begin{equation} \label{2:eqn:exterior}
d:\Om^0(M) \to \Om^1(M), \qquad f \mto df,
\end{equation}
where
\[
df(X)=X(f).
\]
This mapping admits a remarkable extension to a linear operator $d$
on $\Om(M)$, \stn,
\[
d(\Om^p(M)) \subset \Om^{p+1}(M), \qquad \text{ for all }  p \geq 0.
\]
It is called the {\em exterior differentiation operator}, and it is
defined on $\Om^p(M)$ by {\setlength\arraycolsep{2pt}
\begin{equation}\label{2:eqn:extderiv}
\begin{array}{rll}
d\w(X_0, \ldots X_{p}) & = & \sum_{i=0}^{p}(-1)^i X_i(\w(X_0,
\ldots,\wh{X_i},\ldots ,X_{p}))\\
 &  & +\sum_{0 \leq i<j\leq
p}(-1)^{i+j}\w([X_i,X_j], X_0, \ldots ,\wh{X_i}, \ldots ,\wh{X_j},
\dots , X_{p}),
\end{array}
\end{equation}
where $\wh{X_i}$ means that $X_i$ is omitted, and
$[X_i,X_j]=X_iX_j-X_jX_i$. (Note that $[X_i,X_j]$ is indeed a
derivation on $\csm$; it is called the {\em Lie bracket } of $X_i$
and $X_j$). It is routine to check that $d\w$ is an anti-symmetric
$(p+1)$-linear mapping on $\X(M)$. Moreover, when $p=0$, the
definition reduces to $df(X)=X(f)$, showing that $d$ does indeed
extend the operator defined in (\ref{2:eqn:exterior}).

The pair $(\Om(M),d)$ is called the {\em de Rham calculus} of $M$.
Using a standard partition of unity argument, it is easy to show
that every element of $\Om^k(M)$ is a sum of the elements of the
form
\[
f_0df_1\wed df_2 \cdots \wed df_k.
\]
Furthermore, it is also easy to show that $d^2=0$, and that
\[
d(\w_p \wed \w_q)=d\w_p \wed \w_q + (-1)^p \w_p \wed d\w_q.
\]

\begin{rem}
All of the above constructions and definitions are equally well
defined if one uses $C^\infty(M;\bR)$, the \alg of smooth
real-valued \fns on $M$, instead of $\csm$. An important point to
note is that if $T_p(M)$ denotes the tangent plane to $M$ at $p$,
and if $T_p(M;\bR)$ denotes the real tangent plane to $M$ at $p$,
then $T_p(M)=T_p(M;\bR) \oby \bC$.
\end{rem}

\section{Differential calculi}

A {\em positively-graded \alg} is an \alg of the form $\Om=
\bigoplus_{n=0}^{\infty}\Om^n$, where each $\Om^n$ is a linear
subspace of $\Om$, and $\Om^p\Om^q \sseq \Om^{p+q}$, for all $p,q
\geq 0$. (Since there will be no risk of confusion, we shall refer
to a positively graded \alg simply as a graded \algn.) If $\w \in
\Om^p$, then we say that $\w$ is of {\em degree $p$}. A {\em
homogenous mapping of degree $k$} from a graded \alg $\Om$ to a
graded \alg $\LL$ is a linear mapping $h:\Om \to \LL$ \st if $\w \in
\Om^p$, then $h(\w) \in \LL^{p+k}$. A {\em graded derivation} on a
graded \alg $\Om$ is a homogenous mapping of degree $1$ \st
\[
d(\w\w')=d(\w)\w'+(-1)^p \w d\w',
\]
for all $\w \in \Om^p$, and $\w' \in \Om$. A pair $(\Om,d)$ is a
{\em differential \alg} if $\Om$ is a graded \alg and $d$ is a
graded derivation on $\Om$ \st  $d^2=0$. The operator $d$ is called
the {\em differential} of the \algn. A {\em differential \alg
homomorphism} $\f$ between two differential \algs $(\Om,d)$ and
$(\LL,\d)$ is an \alg homomorphism between $\Om$ and $\LL$, which is
also a homogenous mapping of degree $0$ that satisfies $\f \circ d =
\d \circ \f$.

\begin{defn}
A {\em differential calculus} over an \alg $A$ is a differential
\alg $(\Om,d)$, \st $\Om^0=A$, and
\begin{equation} \label{2:eqn:partofunity}
\Om^n=d(\Om^{n-1})\oplus A d(\Om^{n-1}),~~~~~~\text{ for all } n
\geq 1.
\end{equation}
\end{defn}

We note that some authors \cite{DGA2} prefer to omit condition
(\ref{2:eqn:partofunity}) from the definition of a differential
calculus.

Clearly $(\Om(M),d)$ is a differential calculus over $\csm$.
However, it should be noted that for a general calculus there is no
analogue of the graded commutativity of classical differential
forms.

The definition has some immediate consequences. Firstly, if a
differential calculus $(\Om,d)$ is unital (as an \algn) then the
unit of $\Om$ must belong to $\Om^0=A$, and so it must be a unit for
$A$. It then follows that
\[
d(1)=d(1.1)=d(1).1+1.d(1)=2.d(1).
\]
Thus, if $1$ is the unit of a differential calculus, then $d(1)=0$.
If $a_0da_1 \in \Om^1$, then
\[
d(a_0 da_1) = da_0 da_1+a_0 d(da_1) =da_0 da_1.
\]
In general, an inductive argument will establish that
\[
d(a_0  da_1  \ldots  da_n)=da_0 da_1  \cdots  da_n.
\]
If $a_0da_1da_2 \in \Om^2$, and $b_0db_1 \in \Om^1$, then
\begin{eqnarray*}
(a_0da_1da_2)( b_0db_1) & = & a_0da_1d(a_2b_0)db_1 - a_0(da_1)a_2db_0db_1\\
                        & = & a_0da_1d(a_2b_0)db_1 - a_0d(a_1a_2)db_0db_1\\
                        &   & +a_0a_1da_2db_0db_1.
\end{eqnarray*}
Using an inductive argument again, it can be established that, in
general,
\begin{eqnarray*}
(a_0da_1  \cdots  da_n)(a_{n+1}da_{n+2} \cdots da_{n+k})
 = &(-1)^n a_0a_1da_2  \cdots  da_{n+k}\\
& + \sum_{r=1}^n (-1)^{n-r} a_0  da_1   \cdots d(a_ra_{r+1}) \cdots
da_{n+k}.
\end{eqnarray*}

A {\em differential ideal} of a differential calculus $(\Om,d)$ is a
two-sided ideal $I \sseq \Om$, \st $I \cap \Om^0 = \{0\}$, and
$d(I)\sseq I$. Let $\pi$ denote the projection from $\Om$ to
$\Om/I$, and let $\wt{d}$ denote the mapping on $\Om/I$ defined by
$\wt{d}(\pi(\w))=\pi(d(\w))$ (it is well defined since $d(I) \sseq
I$). It is easy to see that, with respect to the natural grading
that $\Om/I$ inherits from $\Om$, $(\Om/I,\wt{d})$ is a differential
\algn. Furthermore, since condition (\ref{2:eqn:partofunity}) is
obviously  satisfied, and $\pi(\Om^0) \cong A$, we have that
$(\Om/I, \wt{d})$ is a differential calculus over $A$.

We say that a graded \alg  $\Om=\bigoplus \Om_{n=0}^\infty$ is {\em
$n$-dimensional} if there exists a \pve integer $n$, \st $\Om^n \neq
\{0\}$ and $\Om^m= \{0\}$, for all $m > n$. If $(\Om,d)$ is
$n$-dimensional, for some \pve integer $n$, then we say that it is
{\em finite-dimensional}; otherwise we say that it is {\em
infinite-dimensional}. Given an infinite-dimensional calculus
$(\Om,d)$ there is a standard procedure for abstracting from it a
new calculus of any finite dimension $n$. Define the \alg
$\Om'=\bigoplus_{k=0}^\infty \Om'^k$ by setting $\Om'^k=\Om^k$, if
$k\le n$, and $\Om'^k=\{0\}$, if $k>n$. Then define a multiplication
$\cdot$ on $\Om'$ as follows: for $\w\in \Om^k$, and $\w' \in
\Om^l$, define $\w\cdot \w'=\w\w'$, if $k+l\le n$, and define
$\w\cdot \w'=0$, if $k+l>n$. Define $d'(\w)=d(\w)$, if $k < n$, and
define $d'(\w)=0$, if $k \ge n$. We call $(\Om',d')$ the
$n$-dimensional calculus {\it obtained from $(\Om,d)$ by
truncation}.

\subsubsection{Differential {\large $*$-}Calculi}

Just as for $C(M)$, we can use complex-conjugation to define an \alg
involution on $\csm$. This involution extends to a unique involutive
conjugate-linear map $*$ on $\Om(M)$, \st $d(\w^*)=(d\w)^*$, for all
$\w \in \Om(M)$. However, $*$ is  not anti-multiplicative; if $\w_p
\in \Om^p(M)$, and $\w_q \in \Om^q(M)$, then
\[
(\w_p \wed \w_q)^*=\w_p^* \wed \w^*_q = (-1)^{pq}\w_q^* \wed \w_p^*.
\]

\bigskip

Let $(\Om,d)$ be a differential calculus over a $*$-\alg $A$. Then,
there exists a unique extension of the involution of $A$ to an
involutive conjugate-linear map $*$ on $\Om$, \st $d(\w^*)=(d\w)^*$,
for all $\w \in \Om$. If it holds that
\[
(\w_p\w_q)^*=(-1)^{pq}\w_q^*\w_p^*, \qquad \text{ for all }\w_p \in
\Om^p,~ \w_q \in \Om^q,
\]
then we say that $(\Om,d)$ is a {\em differential $*$-calculus}.

\subsection{The Universal Calculus}

As would be expected from the quite general nature of the
definition, there can exist several distinct differential calculi
over the same \algn. In fact, as we shall see, there even exist
calculi over $\csm$ that are different from $(\Om(M),d)$. Thus, a
differential calculus is not a strict generalisation of the notion
of the  de Rham calculus.

This is an example of a common feature of moving from the
commutative to \nc setting: it often makes more sense to formulate a
\nc version of a classical structure that is more (or sometimes
less) general than the structure one had originally intended to
generalise.

Over any \alg $A$, we can construct a differential calculus
$\Om_u(A)$ called the universal calculus.  Before we present its
construction, it should be said that while the universal calculus
has a central role to play to the theory of differential calculi, it
is in a sense `too large' to be considered as a suitable
generalisation of the de Rham calculus. The most important thing
about $\Om_u(A)$ is its `universal property', which we shall discuss
later.

\subsubsection{Construction of $\Om_u(\A)$}

Let $A$ be an \alg and define $\Om^1_u(A)=\wt{A} \oplus A$; where
$\wt{A}$ is linear space $A \oplus \bC$ endowed with a
multiplication as defined in (\ref{1:eqn:unitmult}). Define the
structure of an $A$-bimodule on $\Om^1_u(A)$ by
\[
x((a+\l1)\oplus b)=(xa+\l x)\oplus b,
\]
and
\[
((a+\l1)\oplus b)y=(a+\l 1)\oplus by - (ab + \l b) \oplus y,
\]
for $a,b,x,y \in A$. The mapping
\[
d_u:A \to \Om_u^1(A),\qquad  a \mto 1 \oplus a
\]
is a derivation since
\[
d_u(ab)= 1 \oplus ab = (1 \oplus a)b + a \oplus b = d_u(a)b +
ad_u(b).
\]
In general, if $E$ is an $A$-bimodule, and $\d$ a derivation from
$A$ to $E$, then we call the pair $(E,\d)$ a {\em first order
differential calculus} over $A$. We can define a mapping $
i_{\d}:\Om^1_u(A) \to E$, by
\[
i_\d((a + \l1)\oplus b) = a\d(b) + \l \d (b).
\]
Clearly, $i_{\d} \circ d_u = \d$. Furthermore, as a straightforward
calculation will verify, $i_{\d}$ is an $A$-bimodule homomorphism.
The existence of this surjective mapping, for any first order
calculus $(E,\d)$, is known as the {\em universal property} of
$(\Om_u^1(A),d_u)$.

Let $\Om_u^n(A)$ denote the $n$-fold tensor product of $\Om_u^1(A)$
over $A$, and let $\Om_u(A)$ denote the tensor \alg of $\Om_u^1(A)$
over $A$. Clearly, $\Om_u(A)$ is canonically a graded \algn, with
$\Om_u^0=A$.

We would like to define an operator $d_u:\Om_u^1(A) \to \Om_u^2(A)$
\stn, for all $a \in A$, ${d_u(d_u a)=0}$. The fact that
\[
(a+\l1) \oplus b=a\oplus b + \l1 \oplus b=ad_ub+\l d_ub,
\]
implies that we must have
\[
d_u[(a+\l1) \oplus b]=(d_ua)(d_ub).
\]

Let us now progress to define the operator
\[
d_u:\Om^2(A) \to \Om^{3}(A),~~~~\w_1 \oby \w_2 \mto (d_u\w_1) \oby
\w_2 -\w_1\oby (d_u\w_2).
\]
It is well defined since (as a routine calculation will verify)
$d_u(\w a) = (d_u\w)a - \w d_ua$, and $d_u(a\w) = (d_u a)\w +
ad_u\w$. In fact, it is the unique operator on $\Om^2_u(A)$ for
which $d_u(d_u\w))=0$, for all $\w \in \Om^1_u(A)$.

It is not very hard to build upon all of this  to define a
square-zero graded derivation $d_u$ on $\Om_u(A)$ that extends each
$d_u$ defined above. Neither is it too hard to see that this
extension is necessarily unique.

A little reflection will verify that $\Om^n_u(A)$ is spanned by the
elements of the form $a_0d_u(a_1) \oby \cdots \oby d_u(a_n)$ and
$d_u(a_1) \oby \cdots \oby d_u(a_n)$, for every positive integer
$n$. Hence, $(\Om_u(A),d_u)$ is a differential calculus over $A$; we
call it the {\em universal differential calculus} over $A$.

From now on, for sake of simplicity, we shall denote $d_u(a_1) \oby
\cdots \oby d_u(a_n)$ by $d_ua_1 \cdots d_ua_n$.

\subsubsection{Some Remarks on $\Om_u(A)$}

Just like $\Om_u^1(A)$, $\Om_u(A)$ also has a `universal property'.
If $(\LL,\d)$ is another differential calculus over $A$, then, by
the universal property of $\Om_u^1(A)$, there exists a bimodule
mapping $i_\d:\Om_u^1(A) \to \LL$ \st $i_\d \circ d_u = \d$. This
extends to a unique surjective differential \alg homomorphism from
$\Om_u(A)$ to $\LL$, whose restriction to $A$ is the identity. The
existence of this mapping for any  differential calculus over $A$ is
known as the {\em universal property} of $(\Om_u(A),d)$. It is
easily seen that the universal property defines the universal
calculus uniquely up to isomorphism.

For any differential calculus $(\LL,\d)$ the kernel of $i_\d:\Om_u
\to \LL$ is a differential ideal. Hence,
$(\Om_u(A)/\ker(i_\d),\wt{d})$ is well defined as a differential
calculus; clearly it is isomorphic to $(\LL,\d)$. This means that
every differential calculus is isomorphic to a quotient of the
universal calculus.

If $A=\csm$, then it is reasonably clear that the universal calculus
of $A$ will not correspond to the de Rham calculus. The most obvious
reasons are that it is neither unital nor finite-dimensional. The
best we can say is that $(\Om(M),d)$ is isomorphic to a quotient of
$(\Om_u(\csm),d)$.

\bigskip

We remark that there is a natural isomorphism of linear spaces
\[
j:\wt{A}\oby A^{\oby n} \to \Om^n_u(A),
\]
defined by setting
\[
j((a_0+\l 1) \oby a_1 \oby \cdots \oby a_n) = a_0da_1 \cdots da_n +
\l da_1 \cdots da_n.
\]
This fact has an important consequence: for any positive integer
$n$, let  $T_1$ be a linear map from $A^{n+1}$ to a linear space
$B$, and let $T_2$ be a linear map from $A^{n}$ to the same space.
Then there exists a linear map $T:\Om^n_u(A) \to B$ \st ~
$T(a_0d_ua_1d_ua_2 \cdots d_ua_n)=T_1(a_0,a_1, \ldots ,a_n)$,~ and ~
$T(d_ua_1d_ua_2 \cdots d_ua_n)=T_2(a_1,a_2, \ldots ,a_n)$, for all
${a_0, a_1, \ldots , a_n \in A}$.

\section{Derivations-Based Differential Calculi}

We shall now present an example of a differential calculus that is
based on derivations. This approach is a direct generalisation of
the construction of the de Rham calculus presented earlier. It was
originally introduced by M. Dubois--Violette and J. Madore and it is
based upon the work of Koszul. Dubois--Violette and  Madore
originally became interested in derivation based calculi because of
their suitability for use in fuzzy physics. In Chapter $5$ we shall
discuss fuzzy physics, and the role derivation based calculi
originally played in it.

\bigskip

Given an \alg $A$ we define the set of {\em vector fields} over $A$
to be $\Der(A)$, the set of derivations on $A$. We give it the
structure of a complex linear space in the obvious manner. To
generalise the fact that $\X(M)$ is a module over $\csm$ is a little
more problematic because of the noncommutativity of $A$; if $X \in
\Der(A)$, and $a \in A$, then, in general, it does not follow that
$aX$ or $Xa$ are in $\Der(A)$. According to Dubois--Violette, the
most natural solution to this problem is to regard $\Der(A)$ as a
(left) module over $Z(A)$, the centre of $A$. In the classical case
$Z(\csm)=\csm$, and so our definition reproduces the original module
structure.

Generalising the classical case directly, we define
$\underline{\Om}^n_{\Der}(A)$ to be the $Z(A)$-module of
anti-symmetric $Z(A)$-multilinear mappings from $\Der(A)^n$ to $A$.
We then define $\underline{\Om}^0_{\Der}(A)=A$, and
$\underline{\Om}_{\Der}(A)=\bigoplus^\infty_{n=0}
\underline{\Om}^n_{\Der}(A)$. We endow $\underline{\Om}_{\Der}(A)$
with a product that is the direct analogue of the classical product
defined in equation (\ref{2:eqb:wedgeprod}). With respect to this
product $\underline{\Om}_{\Der}(A)$ is canonically a graded \algn.
Then we define a differential $d$ on $\underline{\Om}_{\Der}(A)$
that is the direct analogue of the classical exterior
differentiation operator defined in equation (\ref{2:eqn:extderiv}).
Unfortunately, the pair $(\underline{\Om}_{\Der}(A),d)$ is not
necessarily a differential calculus over $A$, since it may happen
that condition (\ref{2:eqn:partofunity}) is not satisfied.

If $M$ is a compact manifold, then it is obvious that
$(\underline{\Om}_{\Der}(\csm),d)$ will coincide with the de Rham
calculus. However, if $M$ is not compact, or more specifically if
$M$ is not para\cptn, then it may happen that the derivation based
construction of the de Rham calculus will no longer coincide with
the standard construction based on \vbdn s. More explicitly, the
module of sections of the $p$-exterior power of the cotangent bundle
of $M$ may be properly contained in the module of anti-symmetric
$\csm$-multilinear maps on $\X(M)^n$. But, it has been observed that
the calculus constructed using \vbds will coincide with the smallest
differential subalgebra  of $(\underline{\Om}_{\Der}(\csm),d)$ that
contains $\csm$ (see \cite{DBVMAD1} and references therein for
details).

This motivates us to consider $(\Om_{\Der}(A),d)$ the smallest
differential subalgebra of $(\underline{\Om}_{\Der}(A),d)$ that
contains $A$. It can be shown that $\Om_{\Der}(A)$ is the canonical
image of $\Om_u(A)$ in $\underline{\Om}_{\Der}(A)$. This means that
it consists of finite sums of elements of the form $a_0da_1\dots
da_n$, and $da_0\dots da_n$. Hence, $(\Om_{\Der}(A),d)$ is a
differential calculus over $A$. Dubois--Violette has proposed
$(\Om_{\Der}(\csm),d)$ as the most natural \nc generalisation of the
de Rham calculus.

Using derivation based calculi we can formulate \nc versions of many
elements of classical differential geometry. We present the
following examples:

\subsubsection{Noncommutative Lie Derivative}

If $X$ is a vector field over an $n$-dimensional manifold $M$, and
$\w \in \Om^n(M)$, then $\L_X(\w)$, the {\em Lie derivative} of $\w$
\wrt $X$, is defined by setting
\begin{equation}
\L_X(\w)=i_X \circ d(\w) + d \circ i_X(\w);
\end{equation}
where $i_X$ is the mapping from $\Om^n(M)$ to $\Om^{n-1}(M)$ defined
by setting
\begin{equation} \label{2:eqn:cartanid}
i_X(\w)(X_1,X_2, \ldots ,X_{n-1})=\w(X,X_1,X_2, \ldots ,X_{n-1});
\end{equation}
we call $i_X(\w)$ the {\em contraction} of $\w$ by $X$.

The Lie derivative is usually presented in terms of the pullback of
the flow of a vector field. It is intuitively thought of as a
`generalised directional derivative' of $\w$ \wrt $X$.

\bigskip

Returning to the \nc world we see that we can effortlessly
generalise $i_X$. This then enables us to define the {\em Lie
derivative} of ${\w \in \underline{\Om}_{\Der}(A)}$ \wrt $X \in
\Der(A)$ by setting
\[
\L_X(\w)  =i_X \circ d(\w) + d \circ i_X(\w).
\]
In addition, it can be shown that $\Om_{\Der}(A)$ is invariant under
contraction by any element of $\Der(A)$. Thus, we can also define a
generalised Lie derivative for $\Om_{\Der}(A)$. Obviously, both
these definitions correspond to the classical Lie derivative when
$A=\csm$.

\subsubsection{Symplectic Structures}

Symplectic manifolds are objects of central importance in
differential geometry and modern physics. The phase space of every
manifold is canonically a symplectic manifold, and modern
Hamiltonian mechanics is formulated in terms of symplectic
manifolds.

If $M$ is a manifold and $d$ is its exterior differentiation
operator, then a form ${\w \in \Om(M)}$ is said to be {\em closed}
if $\w \in \ker(d)$. A $2$-form $\w$ on $M$ is said to be {\em
non-degenerate at} $p \in M$ if, when $\w(X,Y)(p)=0$, for all $Y \in
\X(M)$, it necessarily holds that $X(f)(p)=0$, for all $f \in C(M)$.
A {\em symplectic manifold} is a pair $(M,\omega)$ consisting of a
manifold $M$, and a closed $2$-form $\omega\in\Omega^2(M)$ that is
nondegenerate at each point.

We say that $\w \in \underline{\Om}^2_{\Der}(A)$ is {\em
non-degenerate} if, for all $a \in A$, there is a derivation $\Ham
(a) \in \Der(A)$ \st $\w(X,\Ham(a))=X(a)$, for all $X \in \Der(A)$.
If $\w$ is an ordinary $2$-form on $M$, then it is nondegenerate in
this sense, \iff it is nondegenerate at each point. This motivates
us to define a {\em symplectic structure} on $A$ to be closed
nondegenerate element of $\underline{\Om}^2_{\Der}(A)$.

\subsubsection{Noncommutative connections}

Our next generalisation makes sense for any differential calculus,
not just the derivation based ones. However, we feel that its
inclusion at this point is appropriate.

If $(E,\pi,M)$ is a smooth \vbd over a manifold $M$, then a {\em
connection} for $E$ is a linear mapping
\begin{equation} \label{2:eqn:connection1}
\nabla:\G^\infty(E) \to \G^\infty (E) \oby_{\csm} \Om^1(M),
\end{equation}
that satisfies the {\em Leibniz rule}
\[
\nabla(sf)=(\nabla s)f + s \oby df, \qquad s \in \G^\infty(E), f \in
\csm.
\]

Equivalently, one can define a {\em connection} to be a bilinear
mapping
\begin{equation} \label{1:eqn:connection2}
\nabla:\X(M) \by \G^\infty(E) \to \G^\infty(E),~~~~(X,s) \mto
\nabla_Xs,
\end{equation}
\st $\nabla_X(sf)=sX(f)+ (\nabla_Xs)f,$ and
$\nabla_{(fX)}s=f\nabla_X s$.

A connection $\nabla$ in the sense of (\ref{2:eqn:connection1})
corresponds to the connection in the sense of
(\ref{1:eqn:connection2}) given by
\[
\nabla(X,s)=i_X(\nabla s);
\]
where, in analogy with (\ref{2:eqn:cartanid}), $i_X(s \oby \w)=s
\w(X)$, for $s \in \G^\infty(E)$, $\w \in \Om(M)$.

If the \vbd in question is the tangent bundle, then we see that a
connection is a generalisation of the directional derivative of one
vector field \wrt another. One important application of connections
is that they are used to define curvature tensors for manifolds.

\bigskip

We can use the Serre--Swan Theorem to generalise connections to the
\nc case. If $E$ is a finitely-generated projective right
$A$-module, and $(\Om,d)$ is a differential calculus over $A$, then
a {\em connection} for $E$ is a linear mapping
\[
\nabla: E \to E \oby_A \Om^1,
\]
that satisfies the generalised Leibniz rule
\[
\nabla(sa)=(\nabla s)a + s \oby da, \qquad s \in E, a \in A.
\]

With a \nc generalisation of connections in hand, one can progress
to define a \nc version of curvature. We shall not pursue this path
here, instead, we refer the interested reader to \cite{DGA2}.

\chapter{Cyclic Cohomology and Quantum Groups}

In this chapter we shall motivate and introduce cyclic (co)homology.
This is a \nc generalisation of de Rham (co)homology due to Alain
Connes. In the process of doing so we shall also introduce \nc
generalisations of volume integrals and de Rham currents. Cyclic
(co)homology has been one of Connes' most important achievements. As
just one example of its usefulness we cite the spectacular
applications it has had to one of the central problems in algebraic
topology, the Novikov conjecture; for details see \cite{CON}.

We shall also introduce \cpt quantum groups. These objects, which
generalise \cpt topological groups, were developed by S. L.
Woronowicz independently of Connes. While the relationship between
the two theories is still not very well understood, it seems that
there are deep connections between them. One important link that has
recently emerged is twisted cyclic cohomology.  It is a
generalisation of cyclic cohomology, and it was developed in Cork by
J. Kustermans, G. J. Murphy, and L. Tuset.

\section{Cyclic Cohomology}

\subsection{The Chern Character}

If $X$ is a \cpt Hausdorff space, then a {\em characteristic class}
for $X$ is a mapping from $V(X)$, the family of vector bundles over
$X$, to $H^*(X)$, the singular cohomology of $X$. (We note that this
definition of a characteristic class is quite general. Usually a
mapping from $V(X)$ to $H^*(X)$ is required to satisfy certain
`natural' conditions before it qualifies as a characteristic class.
However, we have no need to concern ourselves with such details
here.) Characteristic classes are of great value in the study of
vector bundles and can be used in their classification. Very loosely
speaking, characteristic classes measure the extent to which a
bundle is `twisted', that is, `how far' it is from being a trivial
bundle. As an example of their usefulness, we cite the very
important role characteristic classes play in the Atiyah--Singer
index theorem. For an accessible introduction to the theory of
characteristic classes see \cite{MADS,HATCH}.

\bigbreak

A prominent characteristic class is the Chern character. In general,
it is defined to be the unique mapping that satisfies a certain
distinguished set of conditions. However, if $M$ is a compact
manifold, then its Chern character admits a useful explicit
description: If $E$ is a \vbdn, then the {\em Chern character}~
${\textrm{ch}:V(M) \to H^*(M)}$ is defined by setting
\begin{eqnarray}\label{3:eqn:chern}
\textrm{ch}(E) = \ol{\tr}(\exp({\nabla}^2/2\pi i)),~~~~~~E \in V(M);
\end{eqnarray}
where $\nabla$ is a connection for $E$. The value of $\text{ch}(E)$
can be shown to be independent of the choice of connection.

To give the reader a little feeling for the Chern character it is
worth our while to take some time to carefully present equation
(\ref{3:eqn:chern}). For some connection $\nabla$ for $E$, the
extension
\[
{\nabla}:\G^\infty(E) \oby_{\csm} \Om(M) \to \G^\infty(E)\oby_{\csm}
\Om(M)
\]
is defined by setting
\[
{\nabla}(s \oby \w)=\nabla(s) \wed \w + s \oby d\w;
\]
and
\[
\exp({\nabla}^2):\G^\infty(E) \to \G^\infty(E) \oby_{\csm} \Om(M)
\]
is defined by setting
\[
\exp({\nabla}^2)=1+{\nabla}^2+{\nabla}^4/2!+ \ldots +
{\nabla}^{2n}/n!.
\]
Clearly, we can consider the restriction of $\exp({\nabla}^2)$ to
$\G^\infty (E)$ as an element of $\End(\G^\infty(E)) \oby_{\csm}
\Om(M)$. Thus, if we define $\ol{\tr}$ to be the unique mapping
\[
\ol{\tr}:\End(\G^\infty(E)) \oby_{\csm} \Om^{2n}(M) \to \Om^{2n}(M),
\]
for which $\ol{\tr}(A \oby \w)=\tr(A)\w$, then equation
(\ref{3:eqn:chern}) is well defined. It can be shown that
$\textrm{ch}(E)$ is always a closed form. Therefore, because of the
de Rham theorem, we can consider $\ch(E)$ as an element of $H^*(M)$.

\bigskip

Let $X$ be a \cpt Hausdorff space and let $V(X)$ denote the set of
isomorphism classes of vector bundles over $X$. We note that \wrt
the direct sum operation, $V(X)$ is an abelian semigroup. Now, for
any abelian semigroup $S$, there is a standard construction for
`generating' a group from $S$. Let $F(S)$ be the free group
generated by the elements of $S$, and let $E(S)$ be the subgroup of
$F(S)$ generated by the elements of the form $s+s'-(s\oplus s)$;
where $\oplus$ denotes addition in $S$, and $+$ denotes the addition
of $F(S)$. The quotient group $F(S)/E(S)$ is called the {\em
Grothendieck group} of $S$, and it is denoted by $G(S)$. (It is
instructive to note that $G(\bN)=\bZ$.) For any \cpt Hausdorff space
$X$, we define $K_0(X)=G(V(X))$. The study of this group is known as
{\em topological $K$-theory}. It is a very important tool in
topology, and was first used by Grothendieck and Atiyah. For
examples of its uses see \cite{ATIY}. If $X$ is a \cpt smooth
manifold, then an important point for us to note is that the above
construction works equally well for smooth vector bundles.
Pleasingly, the group produced turns out to be isomorphic to
$K_0(X)$.

It is not too hard to show that the Chern character can be extended
to a homomorphism from $K_0(M)$ to $H^*_{dR}(M)$ (see Section
$3.1.3$ for a definition of $H^*_{dR}(M)$) . This extension plays a
central role in topological $K$-theory.

\bigskip

There also exists an algebraic version of $K$-theory. Let $A$ be an
arbitrary associative \alg and define $P[A]$ to be the semigroup of
finitely generated projective left $A$-modules, where the semigroup
addition is the module direct sum. Then define $K_0(A)$ to be the
Grothendieck group of $P[A]$. (It more usual to give a definition of
$K_0$ in terms of projections in $M_n(A)$. However, the above
definition is equivalent, and it is more suited to our needs.)

Recalling the Serre--Swan Theorem, we see that if $M$ is a \cpt
manifold, then $K_0(C(M))=K_0(\csm)=K_0(M)$. The important thing
about the algebraic formulation, however, is that it is well defined
for any associative \alg $A$. The study of the
${\text{$K_0$-groups}}$ of \algs is called {\em algebraic
$K$-theory}. It has become a very important tool in the study of \nc
\calgn s and it has been used to establish a number of important
results; see \cite{MUR} for details.

There also exists a theory that is, in a certain sense, dual to
$K$-theory; it is called $K$-homology. While $K$-theory classifies
the vector bundles over a space, \mbox{$K$-homology} classifies the
Fredholm modules over a space. (A Fredholm module is an operator
theoretic structure, introduced by Atiyah, that axiomizes the
important properties of a special type of differential operator on a
manifold called an elliptic operator; for details see \cite{ATIY}.)
Just like $K$-theory, $K$-homology has a straightforward \nc
generalisation. It has also become a very important tool in the
study of \nc \mbox{\calgn s}. (In recent years, G. Kasparov unified
$K$-theory and $K$-homology into a single theory called
\mbox{$KK$-theory}.)

\bigskip

The existence of a \nc version of $K$-theory prompts us to ask the
following question: Could a \nc generalisation of de Rham cohomology
be defined; and if so, could equation (\ref{3:eqn:chern}) be
generalised to the \nc setting? Connes showed that the answer to
both these questions is yes. In the late $1970$'s he was
investigating foliated manifolds. (Informally speaking, a foliation
is a kind of `clothing' worn on a manifold, cut from a `stripy
fabric'. On each sufficiently small piece of the manifold, these
stripes give the manifold a local product structure.) Canonically
associated to a foliated manifold is a \nc \alg called the foliation
algebra of the manifold. While studying the $K$-homology of this
\algn, Connes happened upon a new cohomology theory, and a means of
associating one of its cocycles to each Fredholm module. After
investigating this cohomology theory, which he named cyclic
cohomology, Connes concluded that it was a \nc generalisation of de
Rham homology. He also concluded that its associated homology theory
was a \nc generalisation of de Rham cohomology. As we shall discuss
later, Connes then proceeded to define \nc generalisations of the
Chern mapping.

In this section we shall show how cyclic (co)homology is
constructed, we shall examine its relationship with classical de
Rham (co)homology, and we shall briefly discuss Connes' \nc Chern
characters.

\subsection{Traces and Cycles} \label{3:sec:tracescycles}

While investigating the $K$-homology of the foliation \alg of
foliated manifolds, Connes developed a method for constructing \nc
differential calculi from Fredholm modules. (Connes' method for
constructing calculi from Fredholm modules is much the same as his
method for constructing calculi from spectral triples; see Chapter
$4$). Using operator traces he then constructed a canonical linear
functional on the $n$-forms of this calculus. Connes' functional had
properties similar to those of the volume integral of a manifold. By
axiomizing these properties he came up with the notion of a graded
trace. This new structure is considered to be a \nc generalisation
of the volume integral. We shall begin this section by introducing
graded traces. Then we shall show how they easily propose cyclic
cohomology theory.

\subsubsection{Graded Traces}

If $M$ is a manifold and $\w$ is a form on $M$, then we say that
$\w$ is {\em non-vanishing} at $p \in M$, if there exists an $f \in
\csm$ \st $\w(f)(p) \neq 0$. If $M$ is $n$-dimensional, then we say
that it is {\em orientable} if there exists an $n$-form $\w \in
\Om^n(M)$ that is  non-vanishing at each point. When $M$ is
orientable it is well known that one can define a complex-valued
linear mapping $\int$ on $\Om^n(M)$ that, in a certain sense,
generalises the ordinary $n$-dimensional volume integral; we call
$\int$ the {\em volume integral} of $M$. An important result about
volume integrals is Stokes' Theorem; it says that if $M$ is a
manifold, then $\int$ vanishes on $d(\Om^{n-1}(M))$; that is, if $\w
\in \Om^{n-1}(M)$, then  $\int d\w = 0$.

\bigskip

We shall now generalise volume integrals to the \nc setting. Let $A$
be an \alg and let $(\Om,d)$ be an $n$-dimensional differential
calculus over $A$. If $n$ is a \pve integer and $\int$ is a linear
functional on $\Om^n$, then we say that $\int$ is {\em closed} if
$\int d\w = 0$, for all $\w \in \Om^{n-1}$. Since $ da_1da_2 \cdots
da_n=d(a_1da_2 \cdots da_n)$, it must hold that $\int da_1da_2
\cdots da_n = 0$, for all closed functionals $\int$. Moreover, if
$\w_p \in \Om^p$, $\w_q \in \Om^q$, and $p+q=n-1$, then, since
$d(\w_p\w_q)=d(\w_p)\w_q+(-1)^p \w_p d\w_q$, it must hold that
\begin{equation} \label{2:eqn:didentity}
\int d(\w_p)\w_q=(-1)^{p+1}\int\w_p d(\w_q).
\end{equation}
An {\em $n$-dimensional graded trace} $\int$ on $\Om$ is a linear
functional on $\Om^n$ \stn, whenever $p+q=n$, we have that
\[
\int \w_p\w_q =(-1)^{pq} \int \w_q\w_p,
\]
for all $\w_p \in \Om^p(A)$, $\w_q \in \Om^q(A)$. Clearly, closed
graded traces generalise volume integrals. A point worth noting is
that the generalisation is not strict, that is, for a manifold $M$,
there can exist closed graded traces on $\csm$ that are not equal to
the volume integral. Another point worth noting is that the
definition of a closed graded trace reintroduces a form of graded
commutativity.

If $\Om$ is an $n$-dimensional calculus over $A$ and $\int$ is an
closed $n$-dimensional graded trace on $\Om$, then the triple
$(A,\Om,\int)$ is called an {\em $n$-dimensional cycle} over $A$.
Cycles can be thought of as \nc generalisations of orientable
manifolds.

\subsubsection{Cyclic Cocycles}

Let $A$ and $(\Om,d)$ be as above, and let $\int$ be a closed,
$n$-dimensional graded trace on $\Om$. We can define a
complex-valued linear \fnn al on $A^{n+1}$ by
\begin{equation}\label{3:eqn:cyclicint}
\f(a_0, a_1, \ldots ,a_n)=\int a_0da_1 \cdots da_n.
\end{equation}
From the graded commutativity of $\int$ we see that
\begin{equation*}
\int a_0da_1da_2 \cdots da_n = (-1)^{n-1} \int da_na_0da_1 \cdots
da_{n-1}.
\end{equation*}
Equation (\ref{2:eqn:didentity})then implies that
\[
\int a_0da_1da_2 \cdots da_n = (-1)^n \int a_{n}da_0da_1 \cdots
da_{n-1}.
\]
This means that
\begin{eqnarray} \label{eqn:2:cycliccocycle}
\f(a_0,a_1, \ldots ,a_n) = (-1)^n \f(a_n, a_0, \ldots ,a_{n-1}).
\end{eqnarray}

\bigbreak

For every positive integer $n$, let $C^n(A)$ denote the space of
complex-valued multilinear mappings on $A^{n+1}$. As we shall see
later, it is common to consider the sequence of mappings
\begin{equation} \label{2:hochdiff}
b_n:C^n(A) \to C^{n+1}(A),~~~~\ps \mto b(\ps);
\end{equation}
where
\begin{eqnarray}\label{3:defn:Hochcobdd}
b_n(\ps)(a_0, \ldots ,a_{n+1}) & = &\sum_{i=0}^n (-1)^i\ps(a_0,\dots,a_{i-1} ,a_ia_{i+1},a_{i+2},\dots,a_{n+1})\nonumber\\
                                               &   & ~~~+(-1)^{n+1}\ps(a_{n+1}a_0,a_1,\dots,a_n).
\end{eqnarray}
We call this sequence, the sequence of {\em Hochschild coboundary
operators}.

If $\f$ is defined as in equation (\ref{3:eqn:cyclicint}), then it
turns out that $b_n(\f)=0$. To show this we shall need to use the
following identity:
\[
\sum_{i=1}^n (-1)^i da_1\cdots d(a_i a_{i+1})\cdots da_{n+1}=-
a_1da_2\cdots da_{n+1} + (-1)^n da_1\cdots da_n(a_{n+1}).
\]
(It is not too hard to convince oneself of the validity of this
identity. For example, if we take the instructive case of $n=3$,
then we see that
\begin{eqnarray*}
-d(a_1a_2)da_3+da_1d(a_2a_3) & = & -d(a_1)a_2da_3-a_1da_2da_3+da_1d(a_2)a_3+d(a_1)a_2da_3\\
                             & = & -a_1da_2da_3+d(a_1)a_2da_3.)
\end{eqnarray*}
With this identity in hand we see that
\begin{eqnarray*}
b_n(\f)(a_0,\dots,a_{n+1}) & = & \int a_0a_1da_2 \cdots da_{n+1}+\sum_{i=1}^n(-1)^i\int a_0da_1 \cdots d(a_ia_{i+1}) \cdots da_{n+1} \\
                       &   & ~~~+(-1)^{n+1}\int a_{n+1}a_0da_1 \cdots da_{n}\\
                       & = & \int a_0a_1da_2\cdots da_{n+1}-\int a_0a_1da_2 \cdots da_{n+1}\\
                       &   & ~+(-1)^n \int a_0da_1 \cdots da_{n}(a_{n+1}) +(-1)^{n+1}\int a_{n+1}a_0da_1 \cdots da_n.\\
                       & = & 0.
\end{eqnarray*}

\bigskip

In general, if $A$ is an \algn, and $\ps$ is an $(n+1)$-multilinear
mapping on $A$ \st $\ps$ satisfies equation
(\ref{eqn:2:cycliccocycle}), and $b_n(\ps)=0$, then we call it an
{\em $n$-dimensional cyclic cocycle}. Interestingly, it turns out
that every cyclic cocycle arises from a closed graded trace; that
is, if $\ps$ is an $n$-dimensional cyclic cocycle, then there exists
an $n$-dimensional cycle $(A,\Om,\int)$  \st
\[
\ps(a_0,a_1,\ldots,a_n)=\int a_0 da_1 \ldots da_n,
\]
for all $a_0,a_1, \ldots a_n \in \Om$. To see this take $\Om_u$(A),
the universal differential calculus over $A$, and consider the
linear functional
\[
\int:\Om_u^n(A) \to \bC,~~~~~ \w \mto \int \w;
\]
where, if $\w = da_1 \cdots da_n$, then $\int \w = 0$, and if $\w =
a_0da_1 \cdots da_n$, then \\ $\int a_0da_1 \cdots da_n = \ps(a_0,
\cdots a_n)$. By definition, $\int$ is closed. Using the fact that
$b_n(\ps)=0$, and that $\ps$ satisfies equation
(\ref{eqn:2:cycliccocycle}), it is not too hard to show that $\int$
is also a graded trace. Hence, if we denote by $\Om'_u(A)$ the
truncation of $\Om_u(A)$ to an $n$-dimensional calculus, then
$(A,\Om'_u(A),\int)$ is an $n$-dimensional cycle from which $\ps$
arises.

\bigskip

It appears that cyclic cocycles are of great importance in \ncgn.
Moreover, all the pieces are now in place to define a cohomology
theory that is based on them. In the following sections we shall
carefully present this theory. We shall begin with an exposition of
Hochschild (co)homology, and then progress to cyclic (co)homology.
We shall also show why cyclic (co)homology is considered a \nc
generalisation of de Rham homology.

\subsection{(Co)Chain Complex (Co)Homology}

Before we begin our presentation of Hochschild and cyclic
(co)homology it seems wise to recall some basic facts about
(co)chain complex (co)homology in general.

A {\em chain complex} $(C_\ast,d)$ is a pair consisting of a
sequence of modules $C_\ast=\{C_n\}_{n=0}^\infty$ (all over the same
ring), and a sequence of module homomorphisms $d=\{d_n:C_n \to
C_{n-1}\}_{n=1}^\infty$, satisfying ${d_n d_{n+1}=0}$. A chain
complex is usually represented by a diagram of the form
\[
\cdots \overset{d_4}{\lto} C_3 \overset{d_3}{\lto} C_2
\overset{d_2}{\lto} C_1 \overset{d_1}{\lto} C_0.
\]
Each $d_n$ is called a {\em differential operator}; for ease of
notation we shall usually omit the subscript and write $d$ for all
differentials. The elements of each $C_n$ are called \mbox{{\em
$n$-chains}}. Those chains that are elements of
$Z_n(C_\ast)=\ker(d_n)$ are called {\em $n$-cycles}, and those that
are elements of $B_n(C_\ast)=\im(d_{n+1})$ are called {\em
$n$-boundaries}. Since the composition of two successive
differentials is $0$, it holds that $B_n(C_\ast) \sseq Z_n(C_\ast)$.
For any positive integer $n$, we define the {\em $n^{\rm
th}$-homology group} of a complex $(C_\ast,d)$ to be the quotient
\[
H_n(C_\ast)=Z_n(C_\ast)/B_n(C_\ast),
\]
We define the $0^{\rm th}$-homology group of $C_\ast$ to be
$C_0/B_0(C_\ast)$. (Note that the homology groups are actually
modules, the use of the term group is traditional.) If
${Z_n(C_\ast)=B_n(C_\ast)}$, for all $n$, then it is clear that the
homology groups of the complex are trivial; we say that such a
complex is {\em exact}.

\bigskip

A cochain complex is a structure that is, in a certain sense, dual
to the structure of a chain complex. A {\em cochain complex}
$(C^\ast,d)$ is a pair consisting of a sequence of modules
$C^\ast=\{C^n\}_{n=0}^\infty$, and a sequence of module
homomorphisms ${d=\{d^n:C^n \to C^{n+1}\}_{n=1}^\infty}$, satisfying
${d^{n+1}d^n =0}$. A chain complex is usually represented by a
diagram of the form
\[
\cdots \overset{d^3}{\lfrom} C^3 \overset{d^2}{\lfrom} C^2
\overset{d^1}{\lfrom} C^1 \overset{d^0}{\lfrom} C^0.
\]
Each $d^n$ is called a {\em boundary operator}; as with
differentials, we shall usually omit the superscript. The
definitions of {\em $n$-cochains, $n$-cocycles}, and {\em
co-boundaries} are analogous to the chain complex definitions. For
$n >0$, the {\em $n^{\rm th}$-cohomology groups} are defined in
parallel with the homological case; the $0^{\rm th}$-group is
defined to be $Z^0(C^*)$.

Given a chain complex $(C_\ast,d)$, we can canonically associate  a
cochain complex to it: For all non-negative integers $n$, let $C^n$
be the linear space of linear functionals on $C_n$, and let $d^n:C^n
\to C^{n+1}$ be the mapping defined by $d^n(\f)(a)=\f(d_n(a))$,  $\f
\in C^n,\, a \in C_{n+1}$. Clearly, the sequence of modules and the
sequence of homomorphisms form a cochain complex. We denote this
complex by $(C^n,d)$, and we say that it is the cochain complex {\em
dual} to $(C_n,d)$. This approach is easily amended to produce a
chain complex dual to a cochain complex.

If $(C_\ast,d)$ and $(D_\ast,d')$ are two chain complexes, then a
{\em chain map} \linebreak ${f:(C_\ast,d) \to (D_\ast,d')}$ is a
sequence of maps $f_n:C_n \to D_n$, such that, for all positive
integers $n$, the following diagram commutes
\begin{displaymath}
\xymatrix{
C_n \ar[r]^{d} \ar[d]_{f_n} & C_{n-1} \ar[d]^{f_{n-1}}\\
D_n \ar[r]^{d'} & D_{n-1}.}
\end{displaymath}
It is clear that a chain map brings cycles to cycles and boundaries
to boundaries. This implies that it induces a map from $H_n(C_\ast)$
to $H_n(D_\ast)$; we denote this map by $H_n(f)$. A {\em chain
homotopy} between two chain maps  ${f,g:C \to D}$ is a sequence of
homomorphisms $h_n:C_n \to D_{n+1}$ satisfying
\[
d'_{n+1}h_n+h_{n-1}d_n=f_n-g_n;
\]
(this is usually written more compactly as $d'h+hd=f-g$). If there
exists a chain homotopy between two chain maps, then we say that the
two maps are {\em homotopic}. Note that, for any $a \in Z_n(C)$,
$d'_{n+1}(h_n(a))+h_{n-1}(d_n(a)) \in B_n(D)$. Thus, for any two
homotopic maps $f$and $g$, $H(f)$ and $H(g)$ are equal. If the
identity map on a chain complex $(C_\ast,d)$ is homotopic to the
zero map, then the chain is called {\em contractible}; the map $h$
is called a {\em contracting homotopy}. A complex that admits a
contracting homotopy is obviously exact.

Clearly, directly analogous results hold for cochain complexes.

\subsubsection{Examples}

The standard example of a cochain complex is the {\em de Rham
complex} of a manifold $M$:
\[
\cdots \overset{d}{\lfrom} \{0\} \overset{d}{\lfrom} \Om^n(M)
\overset{d}{\lfrom} \cdots \overset{d}{\lfrom}  \Om^1(M)
\overset{d}{\lfrom} \Om^0(M).
\]
The cohomology  groups of this complex are called the {\em de Rham
cohomology groups} of $M$; we denote them by $H_{dR}^n(M)$. They are
very important because, as we shall  now see, they provide an
intimate link between the differential structure of $M$ and its
underlying topology.

The typical example of a chain complex is the {\em simplicial
complex} of a manifold $M$: an {\em $n$-simplex} in $M$ is a
diffeomorphic image of the standard $n$-simplex in $\bR^n$; and the
$n$-cochains of the simplicial complex are the formal sums of
$n$-simplices with complex coefficients. It is easy to use the
canonical boundary operator on the standard $n$-simplices to define
a operator from the $n$-simplices to the ${(n-1)\text{-simplices}}$.
This operator can then be extended by linearity to a boundary
operator on the \mbox{$n$-cochains}. The homology groups of the
simplicial complex are called the {\em singular homology groups} and
are denoted by $H_n(M)$. The {\em singular cohomology groups}
$H^n(M)$ are constructed by duality, as explained above.

In one of the fundamental theorems of differential geometry, de Rham
showed that the singular cohomology groups and the de Rham
cohomology groups coincide, that is, $H_{dR}^n(M) \cong H^n(M)$, for
all $n \geq 0$. This is a remarkable result since singular
cohomology is a purely topological object while de Rham cohomology
is derived from a differential structure. For example, $H^0_{dR}(M)=
\bigoplus_{i=0}^k\bC$, where $k$ is the number of connected
components of $M$; we note that the value of $n$ is a purely
topological invariant.

Another fundamental result about de Rham cohomology is Poincar\'e
duality: It states that if $M$ is an $n$-dimensional oriented
manifold, then the $k \th$-cohomology group of $M$ is isomorphic to
the $(n-k)\th$-homology group of $M$, for all \linebreak  ${k=0,1,
\ldots ,n}$.

\bigbreak

Another example of a cochain complex is
\[
\cdots \overset{d_u}{\lfrom} \Om_u^2(A) \overset{d_u}{\lfrom}
\Om_u^1(A) \overset{d_u}{\lfrom} A.
\]
This complex is exact. To see this consider the unique sequence of
maps \linebreak $h=\{h_n:\Om_u^n \to \Om_u^{n+1}\}$ for which
\[
h(a_0da_1 \cdots da_n)=0, \qquad \text{ and } \qquad h(da_0da_1
\cdots da_n)=a_0da_1 \cdots da_n.
\]
It is easily seen that $dh+hd=1$. Hence, $h$ is a contracting
homotopy for the complex.

Let $[\Om_u(A),\Om_u(A)]$ denote the smallest subspace of $\Om_u(A)$
containing the elements of the form
\[
[\w_p,\w_q]=\w_p\w_q-(-1)^{pq}\w_q\w_p,
\]
where $\w_p \in \Om_u^p(A), \w_q \in \Om_u^q(A)$. It was shown in
\cite{KAR2} that $[\Om_u(A),\Om_u(A)]$ is invariant under the action
of $d_u$. Thus, $d_u$ reduces to an operator on the linear quotient
$\Om_u(A)/[\Om_u(A),\Om_u(A)]$. We shall denote the image of
$\Om_u^n(A)$, under the projection onto this quotient, by
$\Om^n_{dR}(A)$. In general, the cohomology of the complex
\[
\cdots \overset{d_u}{\lfrom} \Om_u^2(A) \overset{d_u}{\lfrom}
\Om_u^1(A) \overset{d_u}{\lfrom} A.
\]
is not trivial, and its cohomology groups have some very useful
properties.

\subsection{Hochschild (Co)Homology}

The Hochschild cohomology of associative \algs was introduced by G.
Hochs-child in \cite{HOCH}. One of his original motivations was to
formulate a cohomological criteria for the separability of \algn s.

\bigskip

Let $A$ be an \algn. For every non-negative integer $n$, define
$C_n(A)$, the space of {\em Hochschild $n$-chains} of $A$, to be the
$(n+1)$-fold tensor product of $A$ with itself; that is, define
\begin{equation} \label{3:defn:tensorprod}
C_n(A) =A^{\oby(n+1)}.
\end{equation}
We denote the sequence $\{C_n(A)\}_{n=0}^\infty$ by $C_\ast(A)$. For
all positive integers $n$, let $b'$ be the unique sequence of linear
maps $b'=\{b':C_n(A) \to C_{n-1}(A)\}_{n=0}^\infty$ for which
\begin{equation} \label{2:eqn:barcomplex}
b'(a_0 \oby \ldots \oby a_n)  =  \sum_{i=0}^{n-1}(-1)^i(a_0 \oby
\ldots  \oby a_i a_{i+1} \oby \dots \oby a_n).
\end{equation}

Due to a simple cancellation of terms it holds that $b'^2=b' \circ
b'=0$. Hence, $(C_\ast(A),b')$ is well defined as a complex; it is
known as the {\em bar complex} of $A$. If $A$ is unital, then the
bar complex is exact. To see this consider the unique sequence of
linear maps $s=\{s:C_n(A) \to C_{n+1}(A)\}_{n=0}^\infty$, for which
\[
s(a_0 \oby \ldots \oby a_n)= 1 \oby a_0 \oby \ldots \oby a_n.
\]
It is not too hard to show that $b's+sb'=1$. Hence, $s$ is a
contracting homotopy for the complex. It has been noted by Wodzicki
\cite{WOD} that non-unital algebras whose bar complex is exact have
useful properties. Algebras with this property are now called {\em
H-unital} (homologically unital).

Define $b$ to be the unique sequence of linear maps $b=\{b:C_n(A)
\to C_{n-1}(A)\}_{n=0}^\infty$ for which
\begin{equation} \label{3:defn:hochbound}
b(a_0 \oby \ldots \oby a _{n-1} \oby a_n)= b'(a_0 \oby \ldots \oby a
_{n-1} \oby a_n)+(-1)^n(a_na_0 \oby \ldots \oby a_{n-1}).
\end{equation}
We call $b$ the {\em Hochschild boundary operator}. It is easy to
conclude from the fact that $b'^2=0$, that $b^2=0$. Hence,
$(C_\ast(A),b)$ is a complex; we call it the {\em Hochschild chain
complex} of $A$. The $n^{\rm th}$-homology group of the Hochschild
chain complex is called the {\em $n^{\rm th}$-Hochschild homology
group} of $A$; it is denoted by $HH_n(A)$. (There also exists a more
general formulation of the Hochschild chain complex in terms  of a
general $A$-bimodule $M$; when $M=A^*$ the two formulations
coincide.)

\bigskip

Note that the image of  $b$ in $A$ is $[A,A]$, the {\em commutator
sub\alg} of $A$; that is, the smallest sub\alg containing all the
elements of the form $[a_0,a_1]=a_0a_1-a_1a_0$. Thus,
$HH_0(A)=A/[A,A]$. If $A$ is commutative, then $HH_0(A)=A$.

To gain some familiarity with the definition, we shall calculate the
Hochschild cohomology groups of $\bC$. Firstly, we note that
$\bC^{\oby n+1} \simeq \bC$. Thus, $b$ reduces to $0$ if $n$ is odd,
and to the identity if $n$ is even. This means that $HH_0(\bC)=\bC$,
(confirming that $HH_0(A)=A$ when $A$ is commutative), and
$HH_n(\bC)=0$, if $n \geq 1$.

\bigskip

As would be expected, the cochain complex dual to $(C_\ast(A),b)$ is
called the {\em Hochschild cochain complex}; it is denoted by
$(C^\ast(A),b)$. The $n^{\rm th}$-cohomology group of $(C^*(A),b)$
is called the $n^{\rm th}$-Hochschild cohomology group; it is
denoted by $HH^n(A)$. Note that $b$, the sequence of Hochschild
coboundary operators, is the same as the sequence defined by
equation (\ref{3:defn:Hochcobdd}).

Just like $HH_0(A)$, $HH^0(A)$ has a simple formulation. Recall
that, by definition, $HH^0(A)$ is the set of $0$-cocycles of the
complex. Now,~ $\t \in Z^0(C^*(A))$ \iff $\t(a_0a_1-a_1a_0)=0$, for
all $a_0,a_1 \in A$. Thus, $HH^0(A)$ is the space of traces on $A$.

Just as in the homological case, the $0^{\rm th}$-cohomology group
of the Hochschild complex of $\bC$ is $\bC$, while all the others
groups are trivial.

\subsubsection{Projective Resolutions}
In general, direct calculations of Hochschild homology groups can be
quite complicated. Fortunately, however, there exists an easier
approach involving projective resolutions. If $\E$ is an $R$-module,
then a {\em projective resolution} of $\E$ is an exact complex
\[
\cdots \overset{b_{3}}{\lto} P_{2} \overset{b_{2}}{\lto} P_1
\overset{b_{1}}{\lto} P_0 \overset{b_{0}}{\lto} \E,
\]
where each $P_i$ is a projective $R$-module. Let $A$ be a unital
\algn. Define $A^{\op}$, the {\em opposite \alg} of $A$, to be the
\alg that is isomorphic to $A$ as a linear space by some isomorphism
$a \mto a^o$, and whose multiplication is given by $a^o b^o=(ba)^o$.
Let us denote $B=A \oby A^{\op}$, and let us give $A$ the structure
of a $B$-module by setting $(a\oby b^o)c=acb$. If $(P_\ast,b)$ is a
projective resolution of $A$ over $B$, then
\[
\cdots \overset{1 \oby b}{\lto} A \oby_B P_{3} \overset{1 \oby
b}{\lto} A \oby_B P_2 \overset{1 \oby b}{\lto} A \oby_B P_1
\overset{1 \oby b}{\lto} A \oby P_0
\]
is clearly a complex. Using chain homotopies it can be shown that
the homology of this complex is equal to the Hochschild homology of
$A$. (Those readers familiar with the theory of derived functors
will see that this result can be more succinctly expressed as
$HH_n(A)=\Tor^B_n(A,A)$.) Thus, by making a judicious choice of
projective resolution of $A$, the job of calculating its Hochschild
homology groups can be considerably simplified. A similar result
holds for Hochschild cohomology.

\subsubsection{Hochschild--Konstant--Rosenberg Theorem}

The following theorem is a very important result about Hochschild
homology that comes from algebraic geometry. It inspired Connes to
prove Theorem \ref{2:thm:ctsHoch}, which can be regarded as the
analogue of this result for manifolds. (We inform the reader
unacquainted with algebraic geometry that any introductory book on
the subject will explain any of the terms below that are
unfamiliar.)

\begin{thm} [Hochschild--Konstant--Rosenberg]
Let $Y$ be a smooth complex algebraic variety, let $\O [Y]$ be the
ring of regular functions on $Y$, and let $\Om^q(Y)$ be the space of
algebraic $q$-forms on $Y$. Then there exists a map $\chi$, known as
the {\em Hochschild--Konstant--Rosenberg map}, that induces an
isomorphism
\[
\chi:HH_q(\O[Y]) \to \Om^q(Y), ~~~~~\text{ for all } q \geq 0.
\]
\end{thm}

\subsubsection{Continuous Hochschild (Co)Homology}

The Hochschild cohomology of an \alg is a purely algebraic object,
that is, it does not depend on any topology that could possibly be
defined on the \algn. However, there exists a version of Hochschild
cohomology, called \cts Hochschild cohomology, that does take
topology into account. For applications to noncommutative geometry
it is often crucial that we consider topological \algn s. For
example, while the Hochschild cohomology groups of the \alg of
smooth functions on a manifold are unknown, their continuous version
have been calculated by Connes, as we shall see below.

The continuous Hochschild homology is defined for complete locally
convex \algs and not for general topological algebras: A {\em
locally convex \alg} $A$ is an \alg endowed with a locally convex
Hausdorff topology for which the multiplication $A \by A \to A$ is
\cts. (The use here of the word {\em complete} requires some
clarification. A net $\{x_\l\}_{\l \in D}$ in a locally convex
vector space with a Hausdorff topology is called a {\em Cauchy net}
if, for every open set $U$ containing the origin, there exists a
$\LL \in D$, \stn, for $\k,\l \geq \LL$, $x_\k-x_\l \in U$. The
space is {\em complete} if every Cauchy net is convergent.)

To define the \cts Hochschild homology we shall need to make use of
 Grothendieck's projective tensor product of two locally
convex vector spaces \cite{GROTH} which is denoted by $\oby_\pi$. We
shall not worry too much about the details of the product, but we
shall make two important points about it: Firstly, the projective
tensor product of two locally convex vector spaces (with Hausdorff
topologies) is a complete locally convex vector space (with a
Hausdorff topology); and secondly, the projective tensor product
satisfies a universal property, for jointly \cts bilinear maps, that
is analogous to the universal property that the algebraic tensor
product satisfies for bilinear maps.

If $A$ is a locally convex \alg then we define $C^{\rm cts}_n(A)$,
the space of {\em \cts Hochschild $n$-chains}, to be the completion
of the $(n+1)$-fold projective tensor product of $A$ with itself.
Because of the universal property of $\oby_\pi$, and the joint \cty
of multiplication, each Hochschild boundary operator $b$ has a
unique extension to a \cts linear map on the \cts Hochschild
$n$-chains. The $n^{\rm th}$-homology group of the complex $(HH^{\rm
cts}_n,b)$ is called the {\em \cts Hochschild $n^{\rm th}$-homology
group} of $A$; it is denoted by $HH^{\rm cts}_n(A)$.

We define $C^n_{\rm cts}(A)$, the space of {\em \cts Hochschild
$n$-cochains} of $A$, to be the continuous dual of $C_n^{\rm
cts}(A)$. Note that the Hochschild differentials canonically induce
a unique sequence of linear operators on the \cts cochains. The
$n^{\rm th}$-cohomology group of the resulting complex is called the
{\em \cts Hochschild $n^{\rm th}$-cohomology group} of $A$; it is
denoted by $HH_{\rm cts}^n(A)$.

We saw that the Hochschild (co)homology groups of an \alg could be
computed using projective resolutions. A similar result holds true
for \cts Hochschild (co)homology. However, in the \cts version one
must use topological projective resolutions, for details see
\cite{CON} and references therein.

\subsubsection{Continuous Hochschild Cohomology of $\csm$}

The \alg of smooth \fns on a manifold $M$ can be canonically endowed
with a complete, locally convex topological \alg structure (the
seminorms are defined using local partial derivatives). With a view
to constructing its \cts Hochschild cohomology groups, Connes
constructed a topological projective resolution of $\csm$ over $\csm
\oby_\pi \csm$ \cite{CONHOCH}. Each of the modules of the resolution
arose as the module of sections of a vector bundle over $M \by M$.
(In fact, each of the modules arose as the module of sections of the
vector bundle pullback of $\bigwedge^n(M)$, for some positive
integer $n$, through the map $M \by M \to M,~ (p,q) \mto q$.) Connes
then used the aforementioned result relating Hochschild cohomology
and projective resolutions (or more correctly its \cts version) to
conclude the following theorem. It identifies the \cts Hochschild
$n^{\textrm{th}}$-homology group of $\csm$ with $\D_n$, the \cts
dual of $H_{dR}^n(M)$. The elements of $\D_n$ are called the {\em
$n$-dimensional de Rham currents}.
\begin{thm} \label{2:thm:ctsHoch}
Let $M$ be a smooth compact manifold. Then the map
\[
HH_{\textrm{cts}}^n(\csm) \to \D_n(M),~~~~ \f \mto C_\f,
\]
where
\[
C_\f(f_0df_1   \cdots df_n)=\frac{1}{k!} \sum_{\pi \in
\Perm(p+q)}(-1)^{\sgn{(\pi)}} \f(f_0,f_1, \ldots ,f_n),
\]
is an isomorphism.
\end{thm}

As a corollary to this result it can be shown that
\[
HH^{\textrm{cts}}_n(\csm) \simeq  \Om^n(M).
\]
Thus, Hochschild homology groups generalise differential forms.

\subsection{Cyclic (Co)Homology}

Let $A$ be an \algn, and let $C^n(A)$ and $b$ be as above. Define
\linebreak ${\l=\{\l:C^n(A) \to C^n(A)\}_{n=0}^\infty}$ to be the
unique sequence of linear operators for which
\[
\l(\f)(a_0 \oby \cdots  \oby a_{n-1} \oby a_n)=(-1)^n\f (a_n \oby
a_0 \oby \cdots \oby a_{n-1}),~~~~~ \f \in C^n(A).
\]
If $\f \in C^n(A)$, and $\l(\f)=\f$, then we call $\f$ a {\em cyclic
$n$-cochain}. We denote the subspace of cyclic $n$-cochains by
$C^n_\l(A)$.

Let $b'$ be the operator defined in equation
(\ref{2:eqn:barcomplex}). A direct calculation will show that $(1-\l
)b=b'(1-\l )$. Thus,~ since $C^n_\l(A)=\ker(1-\l)$, we must have
that ${b(C^n_\l(A)) \sseq C^{n+1}_\l(A)}$. Thus, if we denote
$C^*_\l(A)=\{C_\l^n(A)\}_{n=0}^{\infty}$, then the pair
$(C^*_\l(A),b)$ is a subcomplex of the Hochschild cochain complex;
we call it the {\em cyclic complex} of $A$. The
\mbox{$n\th$-cohomology} group of $(C^*_\l(A),b)$ is called the {\em
\mbox{$n\th$-cyclic} cohomology group} of $A$; we denote it by
$HC^n(A)$ ($HC$ stands for homologie cyclique). Note that each
cocycle of the cyclic complex is a cyclic cocycle, as defined in
Section \ref{3:sec:tracescycles}.

To get a little feeling for these definitions, we shall calculate
the cyclic cohomology groups of $\bC$: If $f$ is a non-zero
Hochschild $n$-cochain of $\bC$, then it is clear that $\l(f)=f$
\iff $n$ is even. Thus, $C^{2n+1}_{\l}(\bC) =\{0\}$, and
$C^{2n}_{\l}(\bC) = \bC$. This means that the cyclic complex of
$\bC$ is
\[
\cdots \overset{0}{\lfrom} \bC \overset{\id}{\lfrom} 0
\overset{0}{\lfrom} \bC \overset{\id}{\lfrom} 0 ,
\]
and
\[
HC^{2n}(\bC)=\bC, \text{~~and~~} HC^{2n+1}(\bC)=0.
\]

\subsubsection{Cyclic Homology}
As we stated earlier, there is also a cyclic homology theory. Define
\linebreak $\l=\{\l:C_n(A) \to C_n(A)\}_{n=0}^\infty$ to be the
unique sequence of linear maps for which
\[
\l(a_0 \oby \cdots \oby a_n)=(-1)^n(a_n \oby a_0 \oby \cdots \oby
a_{n-1});
\]
and define the {\em space of cyclic $n$-chains} to be
\[
C^{\l}_n(A)=C_n(A)/\im(1-\l).
\]
Using an argument analogous to the cohomological case, it can be
shown that $b(C_n^\l(A)) \sseq C^\l_{n-1}(A)$. Thus, if we denote
$C^{\l}_*(A)=\{C^\l_n(A)\}_{n=0}^\infty$, then the pair
$(C^{\l}(A),b)$ is a subcomplex of the Hochschild chain complex; we
call it the {\em cyclic complex} of $A$. The $n \th$-homology group
of this complex is called the {\em $n \th$-cyclic homology group} of
$A$; it is denoted by $HC_n(A)$.

The theory of cyclic homology was developed after the cohomological
version. In \cite{LODAYQ}, Loday and Quillen related it to the Lie
\alg homology of matrices over a ring.

\subsubsection{Connes' Long Exact Sequence}

If $I:C_\l^n(A) \to C^{n}(A)$ is the inclusion map of the space of
cyclic $n$-cochains into the space of Hochschild $n$-cochains, and
$\pi:C^n(A) \to C^n(A)/C_\l^n(A)$ denotes the projection map, then
the sequence
\[
0 \lto C^n_\l(A) \overset{I}{\lto} C^n(A) \overset{\pi}{\lto}
C^n(A)/C^n_\l(A) \lto 0
\]
is exact. By a standard result of homological \algn, this sequence
induces a long exact sequence of cohomology groups
\[
\cdots HC^n(A) \overset{I}{\lto} HH^n(A) \overset{B^n}{\lto}
HC^{n-1}(A) \overset{S^n}{\lto} HC^{n+1}(A) \overset{I}{\lto}
\cdots.
\]
The operators $S_n$ are called the {\em periodicity maps}, the
operators $B_n$ are called the {\em connecting homomorphisms}, and
the sequence is known as {\em Connes' long exact sequence}, or the
{\em $SBI$-sequence}. The $SBI$ sequence is of great importance in
the calculation of cyclic cohomology groups because it ties cyclic
cohomology up with the tools of homological \alg available to
calculate Hochschild cohomology from projective resolutions.

There also exists a more efficient approach to calculating cyclic
cohomology groups that is based on the cyclic category (this is an
abelian category introduced by Connes). Connes used this new
structure to show that cyclic cohomology can be realised as a
derived functor. For details see \cite{CON}.

\bigskip

The periodicity maps $S_n:HC^{2n}(A) \to HC^{n+2}(A)$ define two
directed systems of abelian groups. Their inductive limits
\[
HP^0(A)=\lim_\to HC^{2n}(A), \qquad \qquad HP^1(A)=\lim_\to
HC^{2k+1}(A),
\]
are called the {\em periodic cyclic cohomology groups} of $A$. (See
\cite{MUR} for details on directed systems and inductive limits.)

\bigskip

For cyclic homology there also exists a version of Connes' long
exact sequence,
\[
\cdots HC_n(A) \overset{I}{\lfrom} HH_n(A) \overset{B_n}{\lfrom}
HC_{n-1}(A) \overset{S_n}{\lfrom} HC_{n+1}(A) \overset{I}{\lfrom}
\cdots.
\]
Moreover, there exists two {\em periodic cyclic homology} groups
$HP_0(A)$ and $HP_{1}(A)$, whose definitions are, in a certain
sense, analogous to the cohomological case. (We pass over the exact
nature of these definitions since their presentation would require
an excessive digression. Details can be found in \cite{CON}.)

\subsubsection{Continuous Cyclic (Co)Homology}

Beginning with the definition of \cts Hochschild cohomology, it is
relatively straightforward to formulate \cts versions of cyclic, and
periodic cyclic (co)homology for locally convex \algn s. Building on
the proof of \mbox{Theorem \ref{2:thm:ctsHoch}}, Connes established
the following result \cite{CON}.

\begin{thm}
If $M$ is an $n$-dimensional manifold, then
\[
HP_0^{\rm cts}(\csm) \oplus  HP_1^{\rm cts}(\csm) \simeq
\bigoplus_{i=0}^n H_{dR}^{i},
\]
and
\[
HP^0_{\rm cts}(\csm) \oplus  HP^1_{\rm cts}(\csm) \simeq
\bigoplus_{i=0}^n H_{i},
\]
\end{thm}

This result justifies thinking of periodic cyclic homology as a \nc
generalisation of de Rham cohomology, and thinking of cyclic
cohomology as a \nc generalisation of singular homology.

\subsection{The Chern--Connes Character Maps}

As we discussed earlier, Connes discovered a method for constructing
differential calculi from Fredholm modules while he was studying
foliated manifolds. Moreover, he also discovered a way to associate
a cyclic cocycle to each Fredholm module. In fact, this cocycle was
none other than the cocycle corresponding to the module's calculus
and Connes' trace functional. For any \alg $A$, this association
induces a mapping
\[
Ch^0:K^0(A) \to HP^0(A);
\]
it is called the {\em Chern--Connes character mapping} for $K^0$. We
should consider it as a `dual' \nc Chern character.

Let us now return to the question we posed much earlier: Does there
exists a \nc Chern character that maps $K_0(A)$ to $HP_0(A)$? By
generalising \mbox{equation (\ref{3:eqn:chern})} directly to the \nc
case Connes produced just such a mapping
\[
Ch_0:K_0(A) \to HP_0(A);
\]
it is called the {\em Chern--Connes character map} for $K_0$. The
existence of $Ch_0$ is the second principle reason why periodic
cyclic homology is considered to be a \nc generalisation of de Rham
cohomology. (We inform readers familiar with the $K_1$ and $K^1$
groups that, for any \alg $A$, there also exist Chern--Connes maps
$Ch_1:K_1(A) \to HP_1(A)$ and $Ch^1:K^1(A) \to HP^1(A)$.)

Connes used these powerful new tools to great effect in study of
foliated manifolds. He then went on to apply them to a number of
other prominent problems in mathematics. The Chern maps are tools of
fundamental importance in Connes' work. Some mathematicians have
even described them as the ``backbone of \ncgn''. For an in-depth
discussion of the Connes--Chern maps see \cite{CON} and references
therein; for a more accessible presentation see \cite{BROD}.

\section{Compact Quantum Groups}

A {\em character} on a locally compact group $G$ is a continuous
group homomorphism from $G$ to the circle group $T$. If $G$ is
abelian, then set of all characters on $G$ is itself canonically a
locally compact abelian group; we call it the {\em dual group} of
$G$ and denote it by $\wh{G}$. The following result was established
in $1934$ by Pontryagin; for a proof see \cite{RUD}.

\begin{thm}[Pontryagin duality] \label{3:thm:pontryagin}
If $G$ is a locally \cpt abelian group, then the dual of $\wh{G}$ is
isomorphic to $G$.
\end{thm}

Unfortunately, the dual of a non-abelian group is not itself a
group. Thus, in order to extend Pontryagin's result one would have
to work in a larger category that included both groups and their
duals. The theory of quantum groups has its origin in this idea.
(The term {\em quantum group} is only loosely defined. It is usually
taken to mean some type of \nc generalisation of a locally \cpt
topological group.) G. I. Kac, M. Takesaki, M. Enock and J.-M.
Schwartz did pioneering work in this direction, see \cite{ENOCK} and
references therein. One of the more important structures to emerge
during this period was Kac algebras, a theory based on von Neumann
algebras. However, it became apparent in the $1980$s that the then
current theories did not encompass new and important examples such
as S. L. Woronowicz's $SU_q(2)$ \cite{WORSU2} and others found in V.
G. Drinfeld's work \cite{DRIN}. In response a number of new
formulations emerged. Guided by the example of $SU_q(2)$, Woronowicz
developed a new theory based on $C^*$-algebras that he called \cqg
theory. It is Woronowicz's approach that we shall follow here.
Woronowicz later went on to define a more general structure called
an algebraic quantum group. In this setting an extended version of
Pontryagin duality was established, for further details see
\cite{VANDEALE2}. Recently, another approach called locally \cpt
quantum group theory has also emerged. It generalises  a large
number of previous formulations including Kac \algn s, \cqgn s, and
algebraic quantum groups, for details see \cite{KUSVAES}.

\subsection{Compact Quantum Groups}

Recall that a {\em topological group} $G$ is a group endowed with a
Hausdorff topology \wrt which the group multiplication $(x,y) \mto
xy$ and the inverse map $x\mapsto x^{-1}$ are continuous. If $G$ is
\cpt as a topological space, then we say that $G$ is a {\em \cpt
group}. {\em Topological,} and {\em \cpt semigroups} are defined
similarly. Theorem \ref{1:thm:cptsurjective} tells us that we can
recover $G$, as a topological space, from  $C(G)$. However, it is
clear that the semigroup structure of $G$ cannot be recovered; two
\cpt semigroups can be homeomorphic as topological spaces without
being homomorphic as semigroups. With a view to constructing a \nc
generalisation of \cpt groups, we shall attempt below to express the
group structure of $G$ in terms of $C(G)$. However, before we do
this we shall need to present some facts about \calg tensor
products.

\bigbreak

If $H$ and $K$ are two Hilbert spaces, then it is well known that
there is a unique inner product $\<\cdot,\cdot\>$ on their algebraic
tensor product $H \oby K$ \st
\begin{equation*}
\< x \oby y,x' \oby y'\> = \<x,x'\>\<y,y'\>,~~~~~x,x' \in H,~ y,y'
\in K.
\end{equation*}
The Hilbert space completion of $H \oby K$ \wrt the induced norm is
called the {\em Hilbert space tensor product} of $H$ and $K$; we
shall denote it by $H \obyss K$. If $S \in B(H)$ and $T \in B(H)$,
then it can be shown that there exists a unique operator $S \obyss T
\in B(H \obyss K)$ \st
\[
(S \obyss T)(x \oby y)=S(x) \oby T(y), ~~~~~ x \in H, y \in K.
\]
Moreover, it can also be shown that
\begin{eqnarray}\label{3:eqn:optensornorm}
\|S \obyss T\|=\|S\|\|T\|.
\end{eqnarray}

If $A$ and $B$ are two $*$-\algs and $A \oby B$ is their \alg tensor
product, then it is well known that one can define an \alg
involution $*$ on $A \oby B$, \stn, for $a \oby b \in A \oby B$,\,
$(a \oby b)^*=a^* \oby b^*$. This involution makes $A \oby B$ into a
\mbox{$*$-\alg} called the {\em $*$-\alg tensor product} of $A$ and
$B$. If $\A$ and $\B$ are two \calgs and $(U,H)$ and $(V,K)$ are
their respective universal representations, then it can be shown
that there exists a unique injective $*$-\alg homomorphism ${W:\A
\oby \B \to B(H \obyss K)}$, \st ${W(a \oby b)=U(a) \obyss \,
V(b)}$, for $a \in \A$, $b \in \B$. We can use $W$ to define a norm
$\|\cdot\|_*$ on $\A \oby \B$ by setting
\[
\|c\|_* =  \|W(c)\|, ~~~~~~c \in \A \oby \B.
\]
We call it the {\em spatial $C^*$-norm}. Clearly,
$\|cc^*\|_*=\|c\|_*^2$. Hence, the $C^*$-\alg completion of $\A \oby
\B$ \wrt $\|\cdot\|_*$ exists. We call this \calg the {\em spatial
tensor product} of $\A$ and $\B$, and we denote it by $\A \obys \B$.
If we recall that $U$ and $V$ are isometries, then it is easy to see
that equation (\ref{3:eqn:optensornorm}) implies that
\begin{eqnarray}\label{3:eqn:sptensornorm}
\|a \oby b\|_*=\|a\|\|b\|, ~~~~~ a \in \A,~b \in \B.
\end{eqnarray}

Now,~ if $\A,\B,\C,\D$ are \calgn s, and $\f:\A \to \C$ and $\ps:\B
\to \D$ are two $*$-\alg homomorphisms, then a careful reading of
Chapter $6$ of \cite{MUR} will verify that ${\f \oby \ps:\A \oby \B
\to \C \oby \D}$ is \cts \wrt $\|\cdot \|_*$. Hence, it has a unique
extension to a \cts linear mapping from $\A \obys \B$ to $\C \obys
\D$; we denote this mapping by $\f \obys \ps$.

\bigbreak

Let $G$ be a \cpt semigroup, and consider the injective $*$-\alg
homomorphism \linebreak ${\pi:C(G) \oby C(G) \to C(G \by G)}$
determined by ${\pi(f \oby g)(s,t)=f(s)g(t)}$, \linebreak for ${f,g
\in C(G)},~{s,t \in G}$. Since equation (\ref{3:eqn:sptensornorm})
implies that $\|f \oby g\|_*=$ \linebreak $\|\pi(f \oby
g)\|_\infty$, it is easy to use the Stone--Weierstrass theorem to
show that $\pi$ has a unique  extension to a \mbox{$*$-isomorphism}
from $C(G)\obys C(G) \to C(G \by G)$. We shall, from here on,
tacitly identify these two \calgn s.

Consider the $*$-\alg homomorphism
\[
\DEL: C(G) \to C(G \by G), ~~~~ f \mto \DEL(f),
\]
where
\[
\DEL(f)(r,s)=f(rs),~~~~ r,s \in G;
\]
we call it the {\em composition with the multiplication} of $G$. 
As a straightforward examination will verify,
\[
[(\id \obys \DEL)\DEL (f)](r,s,t)=f(r(st)),~~~ \text{and}~~~ [(\DEL
\obys \id)\DEL (f)](r,s,t)=f((rs)t),
\]
for all $r,s,t \in G$. Thus, since $C(G)$ separates the points of
$G$, the associativity of its multiplication is equivalent to the
equation
\begin{equation}
\label{2:eqn:coassociativity} (\id \obys \DEL)\DEL =  (\DEL \obys
\id)\DEL.
\end{equation}

This motivates the following definition: A pair $(\A,\DEL)$ is a
{\em \cpt quantum semigroup} if $\A$ is a unital \calg and $\DEL$ is
a {\em comultiplication} on $\A$, that is, if $\DEL$ is a unital
$*$-\alg homomorphism from $\A$ to $\A \obys \A$ that satisfies
equation (\ref{2:eqn:coassociativity}). 
We say that $(C(G),\DEL)$ is the {\em classical} \cqsg associated to
$G$.

Let $(\A,\DEL)$ be an {\em abelian} \cpt quantum semigroup, that is,
a \cqsg whose \calg is abelian. By the Gelfand--Naimark Theorem
${\A=C(G)}$, where $G$ is the \cpt Hausdorff space $\Om(\A)$. Thus,
we can regard $\DEL$ as a mapping from $C(G)$ to $C(G \by G)$. By
modifying the argument of the second generalisation in Section
$1.2.1$, one can show that there exists a unique \cts mapping $m:G
\by G \to G$, \st $\DEL(f)=f \circ m$, for $f \in C(G)$. The fact
the $\DEL$ is a comultiplication implies that $m$ is associative.
Thus, $G$ is a semigroup and $(\A,\DEL)$ is the \cpt quantum
semigroup associated to it. All this means that the abelian \cqsgn s
are in one-to-one correspondence with the \cpt semigroups.

\bigskip

It is natural to ask whether or not we can identify a `natural'
subfamily of the family of \cqsgn s whose abelian members are in
one-to-one correspondence with the \cpt groups. Pleasingly, it turns
out that we can. Let $\A$ be a \calg and consider the two subsets of
$\A \wh{\oby} \A$ given by
\begin{equation}\label{2:sudset:dense1}
\DEL(\A)(1 \obys \A) = \{\DEL(a)(1 \obys b) : a,b \in \A\},
\end{equation}
and
\begin{equation}\label{2:sudset:dense2}
\DEL(\A)(\A \obys 1) = \{\DEL(a)(b \obys 1) : a,b \in \A\}.
\end{equation}
When $\A=C(G)$, for some \cpt semigroup $G$, then the linear spans
of (\ref{2:sudset:dense1}) and (\ref{2:sudset:dense2}) are both
dense in $C(G)\obys C(G)$ \iff $G$ is a group. With a view to
showing this consider $T$ the unique automorphism of $C(G) \obys
C(G)$ for which ${T(f)(s,t)=f(st,t)}$, for ${f \in C(G \by G),}~s,t
\in G$. If we assume that $G$ is a group, then this map is
invertible; its inverse $T \inv$ is the unique endomorphism for
which \linebreak ${T \inv (f)(s,t)=f(st \inv , t)}$. Since $T$ is
\ctsn, ${T(C(G) \oby C(G))}$ is dense in $C(G) \obys C(G)$. If we
now note that
\begin{equation*}
T(g \oby h)(s,t) = g \oby h(st,t)= g(st)h(t)=[\DEL(g)(1\oby
h)](s,t),
\end{equation*}
for all $g,h \in C(G)$, then we see that the linear span of
$\DEL(C(G))(1 \obys C(G))$ is indeed dense in $C(G) \obys C(G)$. A
similar argument will establish that $\DEL(C(G))(C(G) \obys 1)$ is
dense in $C(G) \obys C(G)$. Conversely, suppose that the linear
spans of (\ref{2:sudset:dense1}) and (\ref{2:sudset:dense2}) are
each dense in ${C(G) \obys C(G)}$, and let $s_1,s_2,t \in G$ \st
$s_1t = s_2t$. As a little thought will verify, the points $(s_1,t)$
and $(s_2,t)$ are not separated by the elements of the linear span
of $\DEL(C(G))(1 \obys C(G))$. Thus, the density of the latter in
$C(G) \obys C(G)$ implies that $s_1 = s_2$. Similarly, the density
of $\DEL(C(G))(C(G) \obys 1)$ in $C(G) \obys C(G)$ implies that if
$st_1=st_2$, then $t_1=t_2$. Now, it is well known that a \cpt
semigroup that satisfies the left- and right-cancellation laws is a
\cpt group. Consequently, $G$ is a \cpt group. This motivated
Woronowicz \cite{WOR0} to make the following definition.

\begin{defn}
A compact quantum semigroup $(\A, \DEL)$ is said to be a {\em
compact quantum group} if the linear spans of $(1 \obys \A)\DEL(\A)$
and $(\A \obys 1)\DEL(\A)$ are each dense in $\A \obys \A$.
\end{defn}

Motivated by the classical case, the density conditions on the
linear spans of $(1 \obys \A)\DEL(\A)$ and $(\A \obys 1)\DEL(\A)$
are sometimes called {\em Woronowicz's left- and right-cancellation
laws} respectively.

\subsubsection{The Haar State}

Let $G$ be a locally compact topological group, and let $\mu$ be a
non-zero regular Borel measure on $G$. We call $\mu$ a {\em left
Haar measure} if it is {\em invariant under left translation}, that
is, if $\mu(gB) = \mu(B)$, for all $g \in G$, and for all Borel
subsets $B$. A {\em right Haar measure} is defined similarly. It is
well known that every locally compact topological group $G$ admits a
left and a right Haar measure that are unique up to positive scalar
multiples. If the left and right Haar measures coincide, then $G$ is
said to be {\em unimodular}. It is a standard result that every \cpt
group is unimodular. Thus, since we shall  only consider \cpt groups
here, we shall speak of {\em the} Haar measure. Furthermore, the
regularity of $\mu$ ensures that $\mu(G)$ is finite. This allows us
to work with the {\em normalised} Haar measure, that is, the measure
$\mu$ for which $\mu(G)=1$. We define the {\em Haar integral} to be
the integral over $G$ \wrt $\mu$. It is easily seen that the left-
and right-invariance of the measure imply {\em left-} and {\em
right-invariance} of the integral, that is,
\begin{equation}
\label{chp5:Haarintinv} \int_Gf(hg)dg = \int_Gf(gh)dg = \int_G f(g)
dg,
\end{equation}
for all $h \in G$, and for all integrable \fns $f$. (Note that since
there is no risk of confusion, we have suppressed explicit reference
to $\mu$.)

Considered as a linear mapping on $C(G)$, $\int$ is easily seen to
be a \pve linear mapping of norm one. In general, if $\f$ is a \pve
linear mapping of norm one on a \calgn, then we call it a {\em
state}. Now, as direct calculation will verify,
\[
(\id \; \obys \int_G) \DEL(f)(h)=\int f(hg)dg,
\]
and
\[
(\int_G \obys \; \id) \DEL(f)(h)=\int f(gh) dg.
\]
Thus, the left- and right-invariance of the integral is equivalent
to the equation
\[
(\id \; \obys \int_G) \DEL(f) = (\int_G \obys \; \id) \DEL(f)=\int_G
f d\mu.
\]

This motivates the following definition: Let $(\A,\DEL)$ be a \cqg
and let $h$ be a state on $\A$. We call $h$ a {\em Haar state} on
$\A$ if
\[
(\id \obys h) \DEL(a)=(h \obys \id) \DEL (a) = h(a),
\]
for all $a \in A$.

The following result is of central importance in the theory of \cqgn
s. It was first established by Woronowicz under the assumption that
$\A$ was separable \cite{WOR1}, and it was later proved in the
general case by Van Daele \cite{VANDEALE}.

\begin{thm}
If $(\A,\DEL)$ is a \cqgn, then there exists a unique Haar state on
$\A$.
\end{thm}

A state $\f$ on a $C^*$-algebra $\A$ is called {\em faithful} if
$\ker(\f)=\{0\}$. Obviously, the Haar state of any classical \cqg is
faithful. However, this does not carry over to the \nc setting. It
is a highly desirable property that $h$ be faithful. Necessary and
sufficient conditions for this to happen are given in \cite{MURTUS}.

\subsubsection{Hopf Algebras}

In the definition of a \cqg we find no mention of a generalised
identity nor any generalisation of the inverse of an element. This
might cause one to suspect that identities and inverses have no
important quantum analogues. However, this is certainly not the
case.

Consider the {\em classical co-unit} $\e:C(G) \to \bC$, defined by
$\e(f)=f(e)$, where $e$ is the identity of $G$. Direct calculation
will verify that
\begin{equation} \label{2:eqn:co-unit}
(\id \wh{\oby} \e)\DEL (f)=(\e \wh{\oby} \id)\DEL(f) = f.
\end{equation}

Consider also the {\em classical anti-pode} $S:C(G) \to C(G)$
defined by setting \linebreak  ${(Sf)(t)=f(t \inv)}$, for $t \in G$.
If $m$ is the linearisation of the multiplication of $C(G)$ on
${C(G) \oby C(G)}$, then
\begin{equation} \label{2:eqn:anti-pode}
m(S \wh{\oby} \id) \DEL (f)= m(\id \wh{\oby} S) \DEL (f) = \e (f) 1.
\end{equation}
In this regard $C(G)$ bears comparison to a standard structure in
mathematics.
\begin{defn}
Let $A$ be a unital $*$-\algn, let $A \oby A$ be the $*$-\alg tensor
product of $A$ with itself, and let $\DEL:A \to A \oby A$ be a
unital $*$-\alg homomorphism \st $(\DEL \oby \id)\DEL=(\id \oby
\DEL)\DEL$. The pair $(A,\DEL)$ is called a {\em Hopf $*$-\alg} if
there exist linear maps $\e:A \to \bC$ and $S:A \to A$ \st
\begin{eqnarray}
(\e \oby \id)\DEL(a)=(\id \oby \e) \DEL (a)= a,
\end{eqnarray}
\[
m(S \oby \id)\DEL(f)=m(\id \oby S)\DEL(f)=\e(f) 1,
\]
for all $a \in A$, where $m$ is the linearisation of the
multiplication of $A$ on $A \oby A$.
\end{defn}
When no confusion arises, we shall, for sake of simplicity, usually
denote a Hopf \alg $(A,\DEL)$ by $A$. We call $\e$ and $S$ the {\em
co-unit} and {\em anti-pode} of the Hopf \algn. It can be shown that
the co-unit and anti-pode of any Hopf \alg are unique. It can also
be shown that the co-unit is multiplicative and that the anti-pode
is anti-multiplicative. A standard reference for Hopf \algs is
\cite{ABE}.

\bigbreak

If $f \in C(G)$, then in general $\DEL(f)$ will not be contained in
$C(G)\oby C(G)$. Thus, $(C(G),\DEL)$ is not a Hopf $*$-\alg.
However, there exists a distinguished unital $*$-sub\alg of $C(G)$
for which this problem does not arise. Let
\[
U: G \to M_n(\bC),~~~~~ g \mto U_g,
\]
be a unitary representation of $G$, and let $U_{ij}(g)$ denote the
$ij \th$-matrix element of $U_g$. Clearly, $u_{ij}:g \to U_{ij}(g)$
is a \cts \fn on $G$. We call any such \fn arising from a
finite-dimensional unitary representation a {\em polynomial \fn} on
$G$. We denote by $\Pol(G)$ the smallest $*$-sub\alg of $C(G)$ that
contains all the polynomial \fnn s. Since $U_{gh}=U_gU_h$,
\[
\DEL(u_{ij})(g,h)=\sum_{k=1}^n u_{ik}(g)u_{kj}(h)=\sum_{k=1}^nu_{ik}
\oby u_{kj} (g,h).
\]
Thus, $\DEL(u_{ij}) \in \Pol(G) \oby \Pol(G)$. Moreover, since $U$
is unitary, ${S(u_{ij})=\ol{u_{ji}} \in \Pol(G)}$. Hence, the pair
$(\Pol(G),\DEL)$ is a Hopf $*$-\alg \wrt the restrictions of $S$ and
$\e$ to $\Pol(G)$. Using the Stone--Weierstrass Theorem it can be
shown that $\Pol(G)$ is dense in $C(G)$.

\bigbreak

Using quantum group corepresentations, the natural quantum analogue
of group representations, Woronowicz established the following very
important generalisation of $(\Pol(G),\DEL)$ to the quantum setting.

\begin{thm} \label{2:thm:WoroHopf}
Let $(\A,\DEL)$ be a \cqgn. Then there exists a unique Hopf $*$-\alg
$(A,\Phi)$ \st $A$ is a dense unital $*$-sub\alg of $\A$ and  $\Phi$
is the restriction of $\DEL$ to $A$.
\end{thm}

We call $(A,\Phi)$ the Hopf $*$-\alg {\em underlying} $G$ (or simply
the Hopf \alg underlying $G$ when we have no need to consider its
$*$-structure.)

The fact that the anti-pode and co-unit are only defined on a dense
subset of the \calg of a \cqg might seem a little unnatural. In
fact, earlier formulations of the compact quantum group definition
included generalisations of the classical co-unit and anti-pode that
were defined on all of $\A$. However, examples were later to emerge
that would not fit into this restrictive framework, and the greater
generality of the present definition was required to include them.
Probably the most important such example is $SU_q(2)$.

\subsubsection{Quantum $SU(2)$}

Recall that the {\em special unitary group} of order $2$ is the
group
\[
SU(2)=\{A \in M_2(\bC):A^*=A \inv,~\det(A)=1\}.
\]
Recall also that
\[
SU(2) = \left\{\left(
  \begin{array}{cc}
    z & -\ol{w} \\
    w & \ol{z} \\
  \end{array}
\right):z,w \in \bC, |z|^2+|w|^2=1\right\}.
\]
Consider $\a'$ and $\g'$ two \cts \fns on $SU(2)$  defined by
\[
\a'\left(
  \begin{array}{cc}
    z & -\ol{w} \\
    w & \ol{z} \\
  \end{array}
\right) =z, ~~~~~~~~~ \text{and}~~~~~~~~~ \g'\left(
  \begin{array}{cc}
    z & -\ol{w} \\
    w & \ol{z} \\
  \end{array}
\right)=w.
\]
It can be shown that the smallest $*$-sub\alg of $C(SU(2))$
containing $\a'$ and $\g'$ is $\Pol (SU(2))$. Clearly, $\Pol
(SU(2))$ is a commutative \algn.

\bigskip

We shall now construct a family of not necessarily commutative \algs
$\{A_q\}_{q \in I}$, $I=[-1,1]\bs \{0\}$, \st when $q=1$, the
corresponding \alg is $\Pol (SU(2))$. Define $A_q$ to be the
universal unital $*$-\alg generated by two elements $\a$ and $\g$
satisfying the relations
\[
\a^* \a + \g^*\g=1, ~~~~~~~~\a \a^* + q^2 \g \g^*=1,
\]
\[
\g \g^*=\g^* \g, ~~~~~~~~q\g \a = \a \g, ~~~~~~~q\g^* \a = \a \g^*.
\]

It is straightforward to show that $A_1$ is commutative, and that
the relations satisfied by its generators are also satisfied by
$\a',\g' \in \Pol(SU(2))$. This means that there exists a unique
surjective $*$-\alg homomorphism $\ta$ from $A_1$ to $\Pol(SU(2))$
\st $\ta(\a)=\a'$ and $\ta(\g)=\g'$. Now, if $\l$ is a character on
$A_1$, then
\[
[\l]= \left(
           \begin{array}{cc}
             \l(\a) & -\l(\g^*) \\
             \l(\g) & ~\l(\a^*) \\
           \end{array}
         \right) \in SU(2).
\]

Clearly $\ta(x)[\l]=\l(x)$ if $x=\a$ or $\g$. This immediately
implies that ${\ta(x)[\l]=\l(x)}$, for all $x \in A_1$. Thus, if the
characters separate the points of $A_1$, then $\ta$ is injective. In
\cite{WORSU2} Woronowicz showed that $A_1$ can be embedded into a
commutative \mbox{\calgn}. If we recall that we saw in Chapter $1$
that the characters of a commutative \calg always separate its
elements, then we can see that Woronowicz's result implies that the
characters of $A_1$ do indeed separate the points of $A_1$. Hence,
$\Pol(SU(2))$ is isomorphic to $A_1$.

Now, we would like to give each $A_q$ the structure of a Hopf \algn.
As some routine calculations will verify, there exist unital
$*$-\alg homomorphisms \linebreak ${\DEL:A_q \to A_q \oby A_q}$,
$\e:A_q \to \bC$, and a unital \alg anti-homomorphism $S:A_q \to
A_q$ \stn:
\[
\DEL(\a)=\a \oby \a - q \g^* \oby \g, ~~~~~~~~\DEL(\g)= \g \oby \a +
\a^* \oby \g,
\]
\[
S(\a) = \a^*, ~~~~~  S(\a^*)=\a, ~~~~~ S(\g) = -q \g, ~~~~~
S(\g^*)=-q \inv \g^*,
\]
\[
\e(\a)=1, ~~~~~~~~~~~~~~~~~~\e(\g)=0.
\]
As some more routine calculations will verify, the pair $(A_q,\DEL)$
is a Hopf $*$-\alg \linebreak with $\e$ and $S$ as counit and
antipode respectively. If we denote the composition with the
multiplication of $SU(2)$ by $\Phi$, then it is not too hard to show
that when $q=1$, $\DEL$ and $\Phi$ coincide on $\Pol(SU(2))$.

Let us now define a norm on each $A_q$ by
\[
\|a\|_u = \sup_{(U,H)} \| U(a)\|,~~~~ a \in A_q;
\]
where the supremum is taken over all pairs $(U,H)$ with $H$  a
Hilbert space and $U$ a unital $*$-representation of $A_q$ on $H$.
Let us denote the completion of $A_q$ \wrt $\| \cdot \|_u$ by
$\A_q$. It can be shown that $\DEL$ has a unique extension to a \cts
mapping $\DEL_u : \A_q \to \A_q \obys \A_q$ \st the pair
$(\A_q,\DEL_u)$ is a \cqgn. It is called {\em quantum $SU(2)$} and
it is denoted by $SU_q(2)$.

\bigbreak

It is easy to show that the mapping
\[
\Om(\A_1) \to SU(2),~~~~~~\t \mto \left(
                                   \begin{array}{cc}
                                     \t(\a) & -\t(\g) \\
                                     \t(\g) & \t(\a^*) \\
                                   \end{array}
                                 \right)
\]
is a homeomorphism. Thus, since $\A_1$ is clearly commutative, the
Gelfand--Naimark theorem implies that ${\A_1 = C(SU(2))}$. If we
recall that $\DEL$ and $\Phi$ coincide on the dense subset
$\Pol(SU(2)) \sseq C(SU(2))$, then we can see that $(\A_1,\DEL)$ is
the classical \cqg associated to $SU(2)$.

When $q \neq 1$ what we have is a purely quantum object. Each such
structure is a prototypical example of a \cqg and is of central
importance in the theory of \cqgn s.

\subsection{Differential Calculi over Quantum Groups}

Earlier, we stated that Woronowicz introduced the concept of a \cqg
in the context of a general mathematical movement to extend
Pontrgagin duality. This is not the whole truth: Woronowicz was also
heavily influenced by physical considerations. In the early $1980$s
examples of quantum groups arose in the work of certain Leningrad
based physicists studying the inverse scattering problem
\cite{RESH}. Woronowicz (as well as Drinfeld and many others) took
considerable inspiration from these structures. Furthermore, the
main reason Woronowicz was interested in quantum groups in the first
place was because he felt that they might have applications to
theoretical physics.

All the topological groups generalised by these physicists were Lie
groups. (Recall that a {\em Lie group} $G$ is a topological group
endowed with a differential structure, \wrt which the group
multiplication $(x,y) \mto xy$ and the inverse map $x \mapsto
x^{-1}$ are smooth.) However, none of the quantum groups that
appeared in the physics literature had a generalised differential
structure associated to them. Woronowicz felt that the introduction
of differential calculi into the study of quantum groups would be of
significant benefit to the theory. Consequently, he formulated the
theory of covariant differential calculi: Let $(A,\DEL)$ be the Hopf
algebra underlying a \cqg $G$, let $\e$ be the counit of $(A,\DEL)$,
and let $(\Om,d)$ be a differential calculus over $A$. If there
exists an \alg homomorphism ${\DEL_\Om : \Om \to A \otimes \Om}$ \st
\[
{(\id_A \otimes d)\DEL_\Om}=\DEL_\Om d;~
\DEL_{\Om}(a)=\DEL(a),\text{ for all }a\in A; ~\text{ and }d1=0,
\]
then we call the triple ${(\Om,d,\DEL_\Om)}$ a {\it left-covariant
differential calculus} over $G$. The fact that $\Om$ is generated as
an algebra by the elements $a$ and $da$, for $a\in A$, means that
there can only exist one such left action $\DEL_\Om$ making
${(\Om,d,\DEL_\Om)}$ a left-covariant calculus. (For this reason we
shall usually write $(\Om,d)$ for $(\Om,d,\DEL_\Om)$.) Differential
calculi that are left-covariant are of prime importance in theory of
differential calculi over quantum groups. (There also exists an
analogous definition of a right-covariant calculus. Calculi that are
both left- and right-covariant are called {\em bicovariant}.)

Let $G$ be a \cpt Lie group and let $(\Om_{\Pol}(G),d)$ denote the
de Rham calculus over the \alg of polynomial functions of $G$; the
construction of this calculus is the same as the construction of the
ordinary de Rham calculus except that one uses $\Pol(G)$ instead of
$\csm$. It is not too difficult to show that $(\Om_{\Pol}(G),d)$ can
be endowed with the structure of a left-covariant calculus over $G$.
Woronowicz gave the first example of a left-covariant differential
calculus over a non-classical \cqg in his seminal paper
\cite{WORSU2}. It is a $3$-dimensional calculus over $SU_q(2)$, and
it is the prototypical example for the theory.

\bigskip

Let $(\Om,d)$ be  a left-covariant differential calculus  over the
Hopf \alg $A$ underlying a \cqgn; and let $\int$ be a linear
functional on $\Om^n$. If
\[
(\id \oby \int)\DEL_\Om(\w)=(\int \w)1,~~~~~\text{for all } \w \in
\Om^n,
\]
then we say that $\int$ is {\em left-invariant}. If whenever $\w \in
\Om$ is such that $\int \w\w'=0$, for all $\w' \in \Om$, we
necessarily have $\w=0$, then we say that $\int$ is {\em
left-faithful}. A little thought will verify that both these
definitions are classically motivated. The natural linear
functionals to study on differential calculi over \cqgn s are the
left-invariant, left-faithful, linear functionals.

\bigskip

Recall that earlier in this chapter we studied graded traces on
differential calculi as generalisations of volume integrals. It
would be quite pleasing if `most' of the natural examples of closed
linear functionals (that is, all the closed, left-invariant,
left-faithful, linear functionals) on left-covariant differential
calculi were of this form. However, this is not the case. When
Woronowicz constructed his differential calculus over $SU_q(2)$, he
also constructed a canonical \mbox{$3$-dimensional}, closed,
left-invariant, left-faithful, linear functional on it. This linear
functional was not a graded trace but a twisted graded trace.
\begin{defn} Let $(\Om,d)$ be a differential calculus over an \alg $A$ and
let $\int$ be a linear functional on $\Om^n$. We say that $\int$ is
an {\em $n$-dimensional twisted graded trace} if there exists a
differential \alg automorphism $\s:\Om \to \Om$ of degree $0$ \stn,
whenever $p+q=n$,
\[
\int \w_p \w_q = (-1)^{pq} \int \s(\w_q) \w_p,
\]
for all $\w_p \in \Om^p,~\w_q \in \Om^q$. We say that $\s$ is a {\em
twist automorphism} associated  to $\int$.
\end{defn}

Clearly, if $\int$ is a twisted graded trace with $\id_A$ as an
associated twist automorphism, then $\int$ is a graded trace. There
may exist more than one twist automorphism for a twisted graded
trace. However, if the twisted graded trace is left-faithful, then
its twist automorphism is unique (this a sufficient but not a
necessary condition).

The definition of a twisted graded trace is analogous to the
definition of a KMS state on a \calgn: A {\em KMS state} $h$ on a
\calg $\A$ is a state for which there exists an \alg automorphism
$\s$, defined on a dense $*$-subalgebra of $\A$, \st $h(ab) =
h(\s(b)a)$, for all elements $a$ and $b$ in the $*$-subalgebra.

It turns out that Woronowicz's linear functional is not an isolated
case, that is, there exist many other examples of closed,
left-invariant, left-faithful, linear functionals that are twisted
graded traces, but not graded traces. The following very pleasing
result, due to Kustermans, Murphy, and Tuset \cite{KUSMURTUS}, shows
why twisted graded traces are so important in the theory of
differential calculi over \cqgn s.
\begin{thm}\label{3:thm:tgtrace}
Let $(\Om,d)$ be an $n$-dimensional left-covariant differential
calculus over a \cqgn; and let $\int:\Om^n \to \bC$ be an
\mbox{$n$-dimensional} linear functional. If $\int$ is closed,
left-invariant, and left-faithful, then it is necessarily a twisted
graded trace.
\end{thm}

We  note that the notion of a left-covariant differential calculus
is well defined over any Hopf \algn, not just those associated to
\cqgn s. However, for Theorem \ref{3:thm:tgtrace} to hold we must
assume the existence of a Haar integral on the \algn. A unital
linear functional $h$ on a Hopf \alg is said to be a {\em Haar
integral} if
\[
(\id_A \oby h)\DEL(a)=(h \oby \id_A) \DEL(a)=h(a)1,~~~~\textrm{for
all } a \in A.
\]

\bigskip

Before we leave this section it is interesting to note that
Kustermans, Murphy, and Tuset's work on twisted graded traces led
them to a new method for constructing left-covariant differential
calculi over Hopf \algn s. In their approach one essentially starts
with a twisted graded trace and then constructs a differential
calculus. (Woronowicz's construction ran in the other direction.)
Their method seems to be a more natural approach than others, and in
\cite{KUSMURTUS} they used it to reconstruct Woronowicz's
$3$-dimensional calculus in an entirely different manner.

\subsection{Twisted Cyclic Cohomology}

The fact that the natural linear functionals to study on
differential calculi over quantum groups are twisted graded traces,
and not graded traces, poses a natural question: Can one construct a
 cohomology theory from twisted graded traces in the same way that
we constructed cyclic cohomology from graded traces? It turns out
that one can.

Let $A$ be a unital \alg and let $\s$ be an \alg automorphism of
$A$. As before, for any positive integer $n$, we let $C^n(A)$ denote
the linear space of complex-valued multilinear maps on $A^{n+1}$,
and we define $C^*(A)=\{C^n(A)\}_{n=0}^\infty $. Let us introduce
the unique sequence of maps $b_\s=\{b_{\s}:C^n(A) \to
C^{n+1}(A)\}_{n=0}^{\infty}$, called the {\em twisted Hochschild
coboundary operators}, \st for $\f \in C^n(A)$ and $a_0, \ldots ,a_n
\in A$,
\begin{eqnarray*}
(b_\s \f)(a_0, \ldots ,a_{n+1})= & \sum_{i=0}^n (-1)^i \f(a_0,
\ldots ,a_{i-1}, a_i a_{i+1}, a_{i+2}, \ldots ,a_{n+1})\\
                                  & + (-1)^{n+1}\f(\s(a_{n+1})a_0,a_1,
\ldots
                                  ,a_n).
\end{eqnarray*}

Just as for the ordinary Hochschild coboundary operators, a
straightforward calculation will show that $b_{\s}^2=0$. Let us also
introduce the unique sequence of maps $\l_\s=\{\l_\s:C^n(A) \to
C^n(A)\}_{n=0}^{\infty}$, called the {\em twisted permutation
operators}, \st for  $\f \in C^n(A)$ and $a_0, \ldots ,a_n \in A$,
\[
\l_\s(\f)(a_0,a_1,\dots,a_n)=
(-1)^n\f(\s(a_n),a_0,a_1,\dots,a_{n-1}).
\]
Define $C^n(A,\s)={\{ \f \in C^n(A): \l_\s^{n+1}(\f)=\f\}}$ (it is
instructive to note that \linebreak $\l_\s^{n+1}(\f)(a_0, \ldots
,a_n)$ $=\f(\s(a_0), \ldots , \s(a_n))$). It can be shown that
 \linebreak $b_\s(C^n(A,\s)) \sseq C^{n+1}(A,\s)$. Thus, if we denote
$C^*(A,\s)=\{C^n(A,\s)\}_{n=0}^\infty$, then the pair
$(C^*(A,\s),b_\s)$ is a cochain complex. We denote its
$n^\text{th}$-cohomology group by $HH^n(A,\s)$, and we call it the
{\em \mbox{$n^\text{th}$-twisted} Hochschild cohomology group} of
$(A,\s)$. Clearly, the twisted Hochschild cochain complex of
$(A,\id_A)$ is equal to the Hochschild cochain complex of $A$.

Let us now define $C^n_\l(A,\s)= \{\f:\f \in C^n(A,\s),
\l_\s(\f)=\f\}$, and $C^*_\l(A)=\{C^n_\l(A,\s)\}_{n=0}^\infty$. It
can be shown that $b_\s(C^n_\l(A,\s)) \sseq C^{n+1}_\l(A,\s)$, and
so, the pair $(C^*_\l(A),\s)$ is a subcomplex of the twisted
Hochschild cochain complex. We call the \mbox{$n\th$-cohomology}
group of this complex the {\em \mbox{$n\th$-twisted} cyclic
cohomology group} of $(A,\s)$, and we denote it by $HC^n(A,\s)$.
Furthermore, we denote the set of $n$-cocycles of the complex by
$Z^n_\l(A,\s)$; we call its elements {\em twisted cyclic
\mbox{$n$-cocycles}}. Clearly, when $\s = \id_A$, the twisted cyclic
cochain complex and the cyclic cochain complex coincide.

\bigskip

Recall that if $\int$ is an $n$-dimensional closed graded trace on
an $n$-dimensional differential calculus over $A$, then the
mulitilinear map $\f:A^{n+1} \to \bC$ defined by setting
\[
\f(a_0,a_1, \ldots ,a_n)=\int a_0da_1 \cdots da_n,
\]
is a cyclic $n$-cocycle. In the twisted setting we have the
following result.

\begin{thm} \label{3:thm:twistedcycliccocycle}
Let $(\Om,d)$ be an $n$-dimensional differential calculus over a
unital algebra $A$ and suppose that $\int$ is an $n$-dimensional,
closed, twisted graded trace on $\Om$. Define the function,
${\f\colon A^{n+1}\to \bC}$, by setting
\[
{\f(a_0,\dots,a_n)}={\int a_0da_1\cdots da_n}.
\]
Let $\s$ be an automorphism of $A$ for which $\int \s(a)\w=\int \w
a$, for all $a \in A$, $\w \in \Om^n$. Then it holds that $\f \in
Z^n_{\l}(A,\s)$. We call $\f$ the {\em twisted cyclic $n$-cocycle}
associated to $(\Om,d)$ and $\int$.
\end{thm}

Recall also that if $\f$ is a cyclic $n$-cocycle of a unital \alg
$A$, then there exists an $n$-dimensional cycle $(\Om,d,\int)$ \st
\[
{\f(a_0,\dots,a_n)}={\int a_0da_1\cdots da_n},
\]
for all $a_0, a_1, \ldots ,a_n$. This result generalises to the
twisted case.

\begin{thm} \label{3:thm:everytwistedcocucle}
Let $\s$ be an automorphism of a unital algebra~$A$ and let
\linebreak ${\f \in Z^n_{\l}(A,\s)}$, for some integer $n \ge 0$.
Then there exists an $n$-dimensional differential calculus $(\Om,d)$
over~$A$ and an $n$-dimensional, closed, twisted graded trace $\int$
on $\Om$ such that $\f$ is the twisted cyclic $n$-cocycle {\em
associated} to $(\Om,d)$ and $\int$.
\end{thm}

The proofs of Theorem \ref{3:thm:twistedcycliccocycle} and Theorem
\ref{3:thm:everytwistedcocucle} amount to suitably modified versions
of the proofs in the non-twisted case, for details see
\cite{KUSMURTUS}.

It is very interesting to note that Connes' ${\text{$S$-$B$-$I$}}$
sequence generalises directly to the twisted setting. It relates the
twisted Hochschild, and twisted cyclic cohomologies; for details see
\cite{KUSMURTUS}. Moreover, it is possible to take the cyclic
category technique for calculating cyclic cohomology, and adapt it
for use in the twisted setting; for details see \cite{HADKRA}.

All this is a very pleasing generalisation of Connes' work. The
straightforwardness with which everything carries over to the
twisted setting suggests that a twisted version of the Connes--Chern
maps could be also be constructed. This area is the subject of
active research.

\bigskip

In \cite{CONMOS} Connes and Moscovici introduced a version of cyclic
cohomology theory for Hopf algebras. For a discussion of the
relationship between Hopf cyclic cohomology and twisted  cyclic
cohomology see \cite{HAJ}.

\subsection{Twisted Hochschild Homology and Dimension Drop}

Recall that in our construction of $SU_q(2)$ we produced a family of
Hopf \algs dependant upon a parameter $q$; when $q=1$, the
corresponding Hopf \alg was equal to $\Pol(SU(2))$. There exists
many other examples of $q$-parameterised families of Hopf \algs that
give the polynomial \alg of a Lie group when $q=1$. Such a
parameterised family is called a {\em quantisation} of the Lie
group. The Hopf \alg corresponding to a particular value of $q$ is
known as its {\em $q$-deformed polynomial \algn}, or simply a {\em
quantum group}. A wealth of examples can be found in \cite{SCHMU};
most of these appeared for the first time in the physics literature.

\bigskip

Motivated by the Hochschild--Kostant--Rosenberg Theorem, and Theorem
\ref{2:thm:ctsHoch}, we make the following definition: For any \alg
$A$, we define its {\em Hochschild dimension} to be
\[
\sup\{n : HH_n(A) \neq 0\}.
\]
For a quantisation of a Lie group it can happen that the Hochschild
\linebreak  dimension of the polynomial \alg of the group is greater
than the \linebreak Hochschild dimension of the deformed \algn s.
This occurrence is known as {\em Hochschild} {\em dimension drop}.
For example, take $SL_q(2)$ the standard $q$-deformed \alg of
$SL(2)$ $=\{A \in M_n(2): \det(A)=1\}$; it is a Hopf algebra
$(\Pol(SL_q(2)),\DEL)$, where $\Pol(SL_q(2))$ is the \alg generated
by the symbols $a$, $b$, $c$, $d$ with relations
\[
ab=qba, \quad ac=qca, \quad bd=qdb, \quad cd=qdc, \quad bc=cb,
\]
\begin{equation}
\label{ad=} ad-qbc=1, \quad da - q^{-1} bc=1,
\end{equation}
for $q$ not a root of unity. It can easily be shown that it admits a
Haar integral. The Hochschild homology of $SL_q(2)$ was calculated
in \cite{MAS}: Its Hochschild dimension was found to be $1$, whereas
the Hochschild dimension of $\Pol(SL(2))$ is $3$. Other examples of
Hochschild dimension drop can be found in the work of Feng and
Tsygan \cite{FENTSY}. This loss of homological information has led
many to believe that the Hochschild and cyclic theories are
ill-suited to the study of quantum groups, and that a generalisation
of Hochschild homology should be introduced to overcome it.

\bigskip

Now, there also exists a twisted version of Hochschild homology. Let
$A$ be an \alg and let $\s:A \to A$ be an \alg automorphism. Denote
by $C^n(A,\s)$ the quotient $C^n(A)/\im(1-\s)$; where $\im(1-\s)$
denotes the image of the map $(1-\s):C_n(A) \to C_n(A)$ defined by
setting
\[
(1-\s)(a_0 \oby a_1 \oby \cdots \oby a_n) = a_0 \oby a_1 \oby \cdots
\oby a_n - \s(a_0) \oby \s(a_1) \oby \cdots \oby \s(a_n).
\]
A routine calculation will show that the pair $(C^*(A,\s),b)$ is a
subcomplex of $(C^*(A),b)$. We call the \mbox{$n^{th}$-homological}
group of this subcomplex the {\em twisted Hochschild
$n^{th}$-homological group} of $(A,\s)$, and we denote it by
$HH_n(A,\s)$. (It is interesting to note that it is not very hard to
build upon these definitions and define a twisted version of cyclic
homology.)

\bigskip

In \cite{HADKRA} Tom Hadfield, a former postdoctoral assistant to
Gerard Murphy in Cork,  and Ulrich Kr\"ahmer, calculated the twisted
Hochschild and cyclic homology of $SL_q(2)$. Following Feng and
Tsygan, they carried out these calculations using \nc Kozul
resolutions. They showed that for certain automorphisms, the
corresponding twisted Hochschild dimension was $3$, equal to the
Hochschild dimension of $\Pol(SL(2))$. In fact, the simplest of
these automorphisms arose in a natural way from the Haar integral on
$SL_q(2)$. The two authors would later generalise this result to one
that holds for all $SL_q(N)$ \cite{HADKRA2}. Hadfield \cite{HAD}
went on to verify that a similar situation holds for all Podl\'e{s}
quantum spheres, and Andrzej Sitarz established analogous findings
for quantum hyperplanes \cite{SIT}. Moreover, quite recently, Brown
and Zhang \cite{BROW} have produced some very interesting new
results in this area.

\chapter{Dirac Operators}

The Schr\"odinger equation for a free particle,
\[
i\hslash \frac{\del}{\del
t}\psi=-\frac{\hslash^2}{2m}\nabla^{2}\psi,
\]
is based on $E = \frac{p^2}{2m}$, the non-relativistic relation
between momentum and kinetic energy, and not on the relativistic
one,
\begin{eqnarray} \label{5:eqn:relenergy}
E=c(m^2c^2+p^2)^{\frac{1}{2}}.
\end{eqnarray}
Furthermore, the fact that the time derivative is of first order and
the spatial derivatives are of second order implies that the
equation is not Lorentz invariant. One of the first attempts to
construct a quantum mechanical wave equation that was in accord with
special relativity was the Klein--Gordan equation,
\[
-\hslash^2 \frac{\del^2}{\del t^2} \ps
=c^2(m^2c^2-\hslash^2\nabla^2) \ps.
\]
It is obtained by canonically quantizing the square of both sides of
\mbox{equation (\ref{5:eqn:relenergy})}. While the equation is
Lorentz invariant, it does have some problems: it allows solutions
with negative energy; and that which one would wish to interpret as
a probability distribution turns out not to be positive definite.
Dirac hoped to overcome these shortcomings by reformulating the
equation so that the time derivative would be of first order. He did
this by rewriting it as
\[
\frac{i \hslash}{c} \frac{\del}{\del t}\ps=\left(m^2c^2 - \hslash^2
\nabla^2 \right)^\frac{1}{2}\ps.
\]
The obvious problem with above expression is that the right hand
side is ill-defined. Dirac assumed that it corresponded to a first
order linear operator of the form
\begin{equation} \label{4:eqn:Diraceleceqn}
(mc A_0 + \hslash \sum_{i=1}^3 A_i  \frac{\del}{\del x_i});
\end{equation}
where each  $A_{i}$ is a matrix, \st
\begin{equation}\label{eqn:anticommutate1}
A_0^2=1,~~~~~~~~~~A_iA_0=-A_0A_i, ~~~~ i=1,2,3,
\end{equation}
and
\begin{equation} \label{eqn:anticommutate2}
A_{i}A_{j}+A_{j}A_{i}=-2\d_{ij}, ~~~~i=1,2,3.
\end{equation}
He then found examples of such matrices in $M_4(\bC)$, namely;
\[
\begin{tabular}{ccc}
$A_{0}=\left(
\begin{array}{cc}
1_2 & 0 \\
0 & -1_2
\end{array}
\right) ;$ & $A _{i}=\left(
\begin{array}{cc}
0 & \s_i \\
-\s_i & 0
\end{array}
\right), ~~~~i=1,2,3;$
\end{tabular}
\]
where $1_2$ is the identity of $M_2(\bC)$, and $\s_i$ are the {\em
Pauli spin matrices}:
\[
\begin{tabular}{ccc}
$\s_{1}=\left(
\begin{array}{cc}
0 & 1 \\
1 & 0
\end{array}
\right) ,$ & $\s _{2}=\left(
\begin{array}{cc}
0 & -i \\
i & 0
\end{array}
\right) ,$ & $\protect\s _{3}=\left(
\begin{array}{cc}
1 &  0 \\
0 & -1
\end{array}
\right). $
\end{tabular}
\]
(It can be shown that $n=4$ is the lowest value for which solutions
to equations (\ref{eqn:anticommutate1}) and
(\ref{eqn:anticommutate2}) can be found in $M_n(\bC)$.) For all this
to make sense it must be assumed that $\ps$ in equation
(\ref{4:eqn:Diraceleceqn}) takes values in $\bC^4$.

While Dirac's reformulation of the Klein--Gordan equation is only
suitable for describing electrons, or more correctly
spin-$\frac{1}{2}$ particles, it was still a great success. Firstly,
the allowed solutions no longer have badly behaved probability
densities. Also, a natural implication of the equation is the
existence of electron spin. (This quantity had previously required a
separate postulate.) But, despite Dirac's efforts, the equation does
allow apparently `unphysical' negative energy solutions. Dirac
initially considered this a `great blemish' on his theory. However,
after closer examination of these solutions he proposed that they
might actually correspond to previously unobserved `antielectron
particles'. According to Dirac these particles would have positive
electrical charge, mass equal to that of the electron, and when an
electron and an antielectron came into contact they would annihilate
each other with the emission of energy according to Einstein's
equation $E=mc^2$. Experimental evidence would later verify the
existence of such particles; we now call them positrons. Thus was
introduced the notion of antimatter. We shall return to this topic
in our discussion of quantum field theory in Chapter $5$.

\bigskip

An important point to note is that the operator
\begin{equation}\label{4:eqn:Diracssquare}
\sum_{i=1}^3 A_i \frac{\del}{\del x_i}
\end{equation}
is in fact a square root of the Laplacian
\[
-\nabla^2:C^\infty(\bR^4,\bC^4) \to C^\infty(\bR^4,\bC^4).
\]
Moreover, it is easy to produce many similar examples of Laplacian
square roots. As we shall see below, one can use Clifford \algs to
unite all these examples into a formal method for constructing
square roots for any Laplacian; we shall call these `square root'
operators, Dirac operators.

The definition and construction of a geometric, or non-Euclidean,
version of Dirac operators is more involved. In fact, it was not
until the $1960$s that examples began to appear. The first was the
K\"ahler--Dirac operator in $1961$, and the second was the
Atiyah--Singer--Dirac operator in $1962$. We shall present
 both of these operators in this chapter.

\bigbreak

At this stage one may well ask what connection there is between
Dirac operators and \ncgn: Dirac operators are important in \ncg
because of Connes' discovery that all the structure of a \cpt
Riemannian (spin) manifold $M$ can be re-expressed in terms of an
algebraic structure based on a Dirac operator. This structure, which
is known as a spectral triple, admits a straightforward \nc
generalisation that can then be considered as a `\nc Riemannian
manifold'.

In the first section of this chapter we shall present the basic
theory of Dirac operators, and in the second section we shall
provide an overview of Connes' theory of spectral triples.

\bigbreak

Recall that in the previous chapter we saw that the borderline
between \cqgn s and cyclic (co)homology is a very active area of
research. The same is true of the borderline between \cqgn s and
spectral triples. In the final section of this chapter we shall
present some of the work that was done in Cork to construct
generalised Dirac operators on \qgn s, and we shall discuss how it
relates to Connes' theory.

\section{Euclidean and Geometric Dirac Operators}

As we stated above, the formal method for constructing Euclidean
Dirac operators involves a special type of \alg called a Clifford
\algn. These \algs generalise the properties of Dirac's matrices
given in equation (\ref{eqn:anticommutate2}).

\subsection{Clifford Algebras}

Let $V$ be a linear space, over $\bK$ ($\bK =\bR$ or $\bC$), endowed
with a symmetric bilinear form $B$. We shall denote by $J$ the
smallest two-sided ideal of $\T(V)$, the tensor algebra of $V$, that
contains all elements of the form $v \otimes v + B(v,v)1$, for $v
\in V$. The algebra $\Cl(V)=\T(V)/J(V)$ is called the {\em Clifford
algebra} of $V$. For sake of convenience we shall denote the coset
$(v_{1} \otimes \cdots \otimes v_{n}+J)$ by $v_{1} \cdots v_{n}$.

If we denote the canonical injection of $V$ into $\Cl(V)$ by $j$,
then it must hold that $j(v)^2=-B(v,v)1$; this is a very important
property of $\Cl(V)$. If $A$ is another unital \alg over $\bK$ for
which there exists a linear map $j':V \to A$ \st
\begin{equation}\label{3:eqn:univprop}
j'(v)^2=-B(v,v)1,
\end{equation}
then it is not very difficult to show that there exists a unique
linear homomorphism $h:\Cl(V) \to A$ \st $j'=h \circ j$; or
equivalently, \st the following diagram commutes:
\begin{displaymath}
\xymatrix{ V \ar[r]^j \ar[dr]_{j'} & \Cl(V) \ar[d]^h\\
 & A. }
\end{displaymath}
This fact is known as the {\em universal property} of $\Cl(V)$, and
it defines $\Cl(V)$ uniquely.

From now on, for sake of convenience, we shall suppress any
reference to $j$ and not distinguish notationally between $V$ and
its image in $\Cl(V)$.

\bigskip

Let us assume that $V$ is finite-dimensional and that
$\{e_i\}_{i=1}^{n}$ is an orthonormal basis  for the space. It
easily follows from the definition of $\Cl(V)$ that
\begin{equation}\label{4:eqn:Cliffrelations}
e_i^2=-1, \text{~~~~~~~~ and ~~~~~~~~} e_ie_j=-e_je_i,
\end{equation}
for all $ i,j = 1, \ldots ,n$, $i \neq j$. Thus, Clifford \algs
generalise the properties of Dirac's matrices given in equation
(\ref{eqn:anticommutate2}). The relations in
(\ref{4:eqn:Cliffrelations}) imply that the set
\[
S=\{1,e_{i_1} \ldots e_{i_k}: 1 \leq i_1 < i_2 < \ldots < i_k \leq
n,~1 \leq k \leq n\}
\]
spans $\Cl(V)$. Moreover, its elements can be shown to be linearly
independent. Hence, it forms a basis for $\Cl(V)$ of dimension
$2^n$. An important consequence of this fact is that $\Cl(V)$ is
linearly isomorphic to the exterior \alg of $V$. The unique mapping
that sends
\begin{equation} \label{4:map:cliffext}
v_1 \wed \cdots \wed v_k  \mto  \sum_{\pi \in \Perm(k)}{\sgn(\pi)}
v_{\pi(1)} \cdots v_{\pi(k)},
\end{equation}
for $k=1,2, \ldots ,k$, is a canonical isomorphism between the two
spaces.

\bigskip

The linear map
\[
V \to \Cl(V), ~~~~~~~ v \mto -v,
\]
satisfies equation (\ref{3:eqn:univprop}). Hence, by the universal
property of Clifford algebras, it extends to an algebra automorphism
$\chi : \Cl(V) \to \Cl(V)$. Clearly, $\chi$ is an involution
operator, that is, $\chi^2=1$. This means that its eigenvalues are
$\pm 1$, and that one can decompose $\Cl(V)$ into positive and
negative eigenspaces $\Cl^+(V)$ and $\Cl^-(V)$. As a moment's
thought will verify, $Cl^+(V)$ is spanned by products of even
numbers of elements of $V$, and $Cl^-(M)$ is spanned by products of
odd numbers of elements of $V$.

\subsubsection{Examples and Representations}

If we take $V=\bR$, with multiplication as the bilinear form, then
$\Cl(\bR)$ has $\{1,e_1\}$ as a basis, where $1$ denotes the
identity of the Clifford \algn, and $e_1$ denotes the identity of
$\bR$. Now $e_1^2=-1$, therefore $\Cl(\bR)$ is isomorphic to $\bC$.

Take $V=\bR^2$ with the Euclidean inner product as the bilinear
form. If $\{e_1,e_2\}$ is the standard basis of $\bR^2$, then the
set $\{1,e_1,e_2,e_1e_2\}$ forms a basis for $\Cl(\bR^2)$. If we
denote
\[
\bg{tabular}{lll} $i=e_1$, & $j=e_2$, & $k=e_1e_2$, \ed{tabular}
\]
then the following relations are satisfied:
\[
\begin{tabular}{llll}
$ij=k,$ & $jk=i,$ & $ki=j,$ & $i^2=j^2=k^2=-1$.
\end{tabular}
\]
Thus, we see that $\Cl(\bR ^2) \simeq \bH$, the \alg of
quaternions.

Now let $\{e_i\}_{i=1}^n$ be the standard basis of $\bC^n$, and let
$B$ be the unique symmetric bilinear form on $\bC^n$ for which
$B(e_i,e_j)=\d_{ij}$. It can be shown that for any natural number
$k$, then
\begin{equation} \label{4:eqn:Cliffrep1}
\Cl(\bC^{2k}) \simeq M_{2^k}(\bC);
\end{equation}
and
\begin{equation} \label{4:eqn:Cliffrep2}
\Cl(\bC^{2k+1}) \simeq {M_{2^k}(\bC) \oplus M_{2^k}(\bC)}.
\end{equation}
This \alg representation of $\Cl(\bC^n)$ is called the {\em spin
representation}. Elementary representation theory now implies that
all irreducible representations of $\Cl(\bC^{2k})$, and
$\Cl(\bC^{2k+1})$, are of dimension $2^k$.

Unfortunately, the representation theory of $\Cl(\bR^n)$, the
Clifford \alg of $\bR^n$ endowed with the Euclidean inner product,
is not as straightforward. However, it does follow similar lines.

\subsection{Euclidean Dirac Operators}

We are now ready to define a generalised version of Dirac's
operator. (Note that in this definition we consider $\bR^n$ as
equipped with its usual inner product.)

\begin{defn} \label{4:defn:flatdirac}
Let $c:\Cl(\bR^n) \to \End(\bC^m)$ be an \alg representation of
$\Cl(\bR^n)$ on $\bC^m$, and let $\{e_i\}$ be the standard basis of
$\bR^n$. The operator
\[
D:C^{\infty}(\bR^n,\bC^m) \to
C^{\infty}(\bR^n,\bC^m), ~~~~~ f \mto \sum_{i=1}^n c(e_i)\frac{\del f}{\del x_i},
\]
is called the {\em Dirac operator associated to $c$}.
\end{defn}
The relations in (\ref{4:eqn:Cliffrelations}) easily imply that
$D^2= -\nabla^2$.

As an example, let us construct a Dirac operator for the Laplacian
on $C^\infty(\bR^2, \bC^2)$. We saw above that $\Cl(\bC^2) \simeq
\bH$. Now a routine calculation will show that there exists a unique
homomorphism of real \algs that maps $i \mto i\s_1$, $j \mto i\s_2$,
and $k \mto -i\s_3$. The Dirac operator associated to this
representation is
\[
D=\left(%
\begin{array}{cc}
  0                & i\del_1+\del_2  \\
  i\del_1-\del_2   & 0 \\
\end{array}%
\right),
\]
where $\del_i=\frac{\del}{\del x_i}$. As a straightforward
calculation will verify, the square of $D$ is indeed the Laplacian.
Alternatively, by mapping $i \mto i\s_1$, $j \mto i\s_3$, and $k
\mto i\s_2$, we get the Dirac operator
\[
D=i\left(%
\begin{array}{cc}
  \del_2  &  \del_1 \\
  \del_1  & -\del_2 \\
\end{array}%
\right).
\]
Again, a straightforward calculation will verify that $D$ squares to
give the Laplacian.

Although the details are a little more involved, it is not too hard
to show that the operator of equation (\ref{4:eqn:Diracssquare})
also fits into our framework.

\subsection{Geometric Dirac Operators}\label{sec:4:geomdirac}

The Dirac operators constructed above can be viewed as operating on
the smooth sections of the trivial bundle $\bR^n \by \bC^m$. This
makes it natural to consider the idea of constructing generalised
Dirac operators that would operator on the sections of \vbds over
manifolds. We call such operators {\em geometric} Dirac operators.

\subsubsection{Riemannian Manifolds}

The construction of a geometric Dirac operator for a manifold $M$
requires a choice of Riemannian metric tensor for $M$. Thus, we
shall need to recall some details about Riemannian manifolds.

If $g$ is a real, or complex, rank-$(0,2)$ tensor field on an
$n$-dimensional manifold $M$, then by our comments in Chapter $2$ on
the locality of tensor fields, $g$ will induce a bilinear form on
each real, or complex, tangent plane of $M$ respectively. A {\em
Riemannian metric tensor} is a real rank-$(0,2)$ tensor field  ${g
\in \fT^0_2(M;\bR)}$ such that, for each $p\in M$, the induced
bilinear form ${g_p:T_p(M;\bR) \times T_p(M;\bR) \to \bR}$ is an
inner product. Clearly, $g_p$ has a unique extension to a
complex-valued symmetric bilinear form on the complex tangent plane
$T_p(M)$, which we shall also denote by $g_p$. It is routine to show
that any inner product on $T(M;\bR)$, in the sense of Section $1.3$,
induces a Riemannian metric tensor on $M$. This means that one can
find a Riemannian metric tensor for any manifold $M$. A pair $(M,g)$
consisting of a manifold $M$ and a Riemannian metric $g$ is called a
{\em Riemannian manifold}.

Using $g$ we can define two mutually inverse $\csm$-module
isomorphisms between $\Om^1(M)$ and $\X(M)$. The maps, known as the
{\em flat}, and {\em sharp, musical isomorphisms} respectively, are
\[
\flat:\X(M) \to \Om^1(M),~~~~X \mto X^\flat;
\]
where $X^\flat(Y)=g(X,Y)$, for all $Y \in \X(M)$; and
\[
\sharp:\Om^1(M) \to \X(M),~~~~\w \mto \w^\sharp;
\]
where $g(\w^\sharp,Y)=\w(Y)$, for all $Y \in \X(M)$. (Note that a
simple local argument will show that $\sharp$ is well defined.) We
can use the sharp musical isomorphism to define a rank-$(2,0)$
tensor field $g\inv$ on $M$ by
\[
g\inv:\Om^1(M) \by \Om^1(M),~~~~(\w_1,\w_2) \mto g(\w_1^\sharp,
\w_2^\sharp).
\]
This can then be extended to a unique symmetric bilinear mapping
\[
{g:\Om^p(M) \by \Om^p(M) \to \csm},
\]
by setting
\[
g(\w_1 \wed  \cdots \wed \w_p,\w'_1 \wed \cdots \wed \w'_p)=\det[g
\inv(\w_i,\w_j')]_{ij}.
\]

\subsubsection{Clifford Bundles and Clifford Modules}

The analogue of the Clifford \alg of $\bR^n$ is a smooth \alg bundle
over $M$ called the Clifford bundle of $M$; the definition of an
{\em \alg bundle} is essentially the same as that of a vector bundle
except that the fibres are no longer linear spaces but \algn s, and
all linear mappings are replaced by \alg mappings. Let us denote the
set $\bigcup_{p \in M} \Cl(T_p(M))$ by $\Cl(M)$ (where $\Cl(T_p(M))$
is the Clifford \alg of the tangent plane of $M$ at $p$, with $g_p$
as the bilinear form) and define a projection $\pi:\Cl(M) \to M$ in
the obvious way. Recalling the isomorphism induced by the mappings
in (\ref{4:map:cliffext}), we see that we can endow $\Cl(M)$ with a
unique topology that makes it a smooth \alg bundle that is
isomorphic, as a smooth vector bundle, to the exterior bundle of
$M$. We call $\Cl(M)$ the {\em Clifford bundle} of $M$.

Equivalently, one can define the Clifford bundle of a manifold using
a transition \fn argument. Let $\{U_\a\}$ be an open covering of $M$
by base \nbds for $T(M)$, and let $\{g_{\a\b}:U_\a \cap U_\b \to
\GL(2^n,\bC)\}_{\a\b}$ be the corresponding set of transition
functions. For each $p \in U_\a \cap U_\b$, we can use the
Riemannian metric of $M$ to endow the domain and codomain of each
$g_{\a\b}(p)$ with canonical symmetric bilinear forms. Obviously,
each $g_{\a\b}(p)$ is isometric \wrt these bilinear forms. Using the
universal property of $\Cl(V)$, it is not too difficult to show that
there exists a unique \alg homomorphism
${\cl(g_{\a\b}(p)):\Cl(\bR^n) \to \Cl(\bR^n)}$ \st the following
diagram is commutative:
\begin{displaymath}
\xymatrix{\bR^n \ar[d]_{j} \ar[r]^{g_{\a\b}(p)}&
\bR^n  \ar[d]^{j'}\\
\Cl(\bR^n) \ar[r]_{\cl(g_{\a\b}(p))} & \Cl(\bR^n). }
\end{displaymath}
We can use this fact to define a set of functions
\[
\wh{g}_{\a\b}:U_\a \cap U_\b \to \GL(2^n,\bC),~~~~ p \to \cl
(g_{\a\b}(p)),~~~~U_\a \cap U_\b \neq \emptyset.
\]
This set can easily be shown to satisfy the conditions of the smooth
analogue of Proposition \ref{1:prop:transfns}. The associated smooth
vector bundle, endowed with the obvious smooth \alg bundle
structure, is isomorphic to $\Cl(M)$.

There also exists a more formal construction of the Clifford bundle
in terms of principle bundles. For details on this approach see
\cite{LAW}.

\bigskip

Finally, let us introduce the analogue of a Clifford \alg
representation: a {\em Clifford module} for $M$ is a pair $(E,c)$
where $E$ is a smooth vector bundle over $M$, called the {\em spinor
bundle}, and $c$ is a module homomorphism from $\G(\Cl(M))$ to
$\End(\G(E))$. Smooth sections of a spinor bundle are called {\em
 spinors}. Using Theorem \ref{thm:functorf}, it is not
hard to show that if $b \in \G^\infty(\Cl(M))$, then, for every
spinor $s$, $c(b)s$ is also a spinor.

We note that since $\X(M)$ is canonically a subset of $\G(\Cl(M))$,
any Clifford module $(E,c)$ induces a module homomorphism $c:\X(M)
\to \End(\G(E))$ by restriction. In turn, this homomorphism induces
a homomorphism $c:\Om^1(M) \to \End(\G(E))$ defined by setting
$c(\w)=c(\w ^\sharp)$, for $\w \in \Om^1(M)$.

\subsubsection{Geometric Dirac Operators}

The pieces are now in place to define a generalised Dirac operator.

\begin{defn} \label{4:defn:geomdirac}
Let $(E,c)$ be a Clifford module for a manifold $M$, and let
$\nabla$ be a connection for $E$. The associated {\em Dirac
operator} $D$ is
\[
D=\wh{c} \circ \nabla:\G^\infty(E) \to \G^{\infty}(E);
\]
where $\wh{c}:\G^{\infty}(E) \oby \Om^1(M) \to \G^\infty(E)$ is the
unique module homomorphism for which
\[
\wh{c}(s \oby \w)=c(\w)s.
\]
\end{defn}

\bigbreak

A  local orthonormal basis $\{E_i\}_{i=1}^n$ of $T(M)$ over a \nbd
$U$ is a local basis over $U$, \st $g(E_i,E_j)|_U=\d_{ij}1|_U$, for
$i=1,2, \ldots ,n$. It is clear that a local orthonormal basis can
be constructed from any local basis. Now if $\{E_i\}_{i=1}^n$ is a
local orthonormal basis of $T(M)$ over a \nbd $U$, then for any $X
\in \X(M)$, $X|_U=\sum_{k=1}^m g( E_i,X) E_i(p)|_U$. Therefore, if
$(E,c)$ is a Clifford module and $s \in \G^\infty(E)$, then
\[
\nabla_X s|_U = \sum_{i=1}^n g(E_i, X) \nabla_{E_i} s\big |_U.
\]
Consequently,
\[
\nabla s|_U=\sum_{i=1}^n \nabla_{E_i} s\otimes E_i^{\flat}\big |_U .
\]
This means that if $D$ is the  Dirac operator associated to $(E,c)$
and $\nabla$, then {\setlength\arraycolsep{2pt}
\begin{eqnarray}
D (s)|_U  & = & \wh{c}\left(\sum_{i=1}^n \nabla_{E_i} s \otimes
E_i^{\flat} s\right)\big |_U=\sum_{i=1}^n c((E_i^{\flat})^\sharp)
\nabla_{E_i} s \big |_U \\
          & = & \sum_{i=1}^n c(E_i) \nabla_{E_i} s\big |_U.
\end{eqnarray}

An immediate consequence of this is that if $M=\bR^n$, $E=\bR^n \by
\bC^m$, and $\nabla$ is the canonical connection for $E$, then the
Dirac operators of Definition \ref{4:defn:flatdirac} and Definition
\ref{4:defn:geomdirac} coincide.

\subsubsection{Hodge Theory and the K\"ahler--Dirac Operator}

An interesting example of a geometric Dirac operator comes from
noting that $\Cl(M)$ is itself canonically a spinor bundle.
Moreover, since the exterior bundle of $M$ is isomorphic to $\Cl
(M)$ as a vector bundle, $\LL(M)$ is also canonically a spinor
bundle.

Now for every Riemannian manifold $M$, there is a unique connection
$\nabla^g$ for the tangent bundle that is {\em compatible with the
metric}, that is,
\[
X(g(Y,Z))= g(\nabla^g_X Y,Z) + g(Y,\nabla^g_XZ), ~~~~X,Y,Z \in
\X(M);
\]
and {\em torsion-free}, that is,
\[
\nabla^g_XY-\nabla^g_YX=[X,Y].
\]
We call $\nabla^g$ the {\em Levi-Civita connection}. It induces a
connection for $\Om^1(M)$, also denoted by $\nabla^g$, that is
defined by
\[
[\nabla^g_X \w](Y) = X\big(\w(Y)\big) - \w(\nabla^g_X Y),~~~~\w \in
\Om(M).
\]
(A simple local argument will show that $\nabla^g$ is well defined).
This connection can then be extended to a unique connection
$\ol{\nabla}$ for $\Om(M)$  \st
\[
\ol{\nabla}(\w_1 \wed \cdots \wed \w_n) = \sum_{i=1}^n \w_1 \wed
\cdots \wed \ol{\nabla} \w_i \wed \cdots \wed \w_n,
\]
for $\w_i \in \Om^1(M)$. We call $\ol{\nabla}$ the {\em
Levi--Civita} connection for $\Om(M)$.

The Dirac operator associated to $\Om(M)$ and $\ol{\nabla}$ is
called the {\em K\"{a}hler--Dirac operator}. It has a pleasing
representation in terms of the Hodge codifferential, which we shall
now introduce.

\bigbreak

As is well known, there exists a unique invertible linear mapping
$\ast:\Om(M) \to \Om(M)$ called the {\em Hodge operator} such that
if $\w_1,\w_2 \in \Om^p(M)$ then
\[
\w_1 \wed \ast \w_2 = g(\w_1,\w_2)d\mu;
\]
where $d\mu$ is the Riemannian volume form, (see Section
(\ref{sec:4:ncint}) for details). Clearly,  $\ast$ must map
$p$-forms to $(n-p)$-forms.

Using the Hodge operator, we can endow $\Om(M)$ with an inner
product by defining
\begin{equation}\label{4:eqn:hodgeinnerproduct1}
\<\w,\w'\>=\int \ol{w'} \wed \ast \w,
\end{equation}
if $\w$ and $\w'$ have the same degree, and
\begin{equation}\label{4:eqn:hodgeinnerproduct2}
\<\w,\w'\>=0
\end{equation}
otherwise. If $\w \in \Om^{p-1}(M)$ and $\w' \in \Om^{p}(M)$, then
\[
d(\ol{\w} \wed \ast \w')=(d\ol{\w})\wed \ast \w'+(-1)^{p-1}\ol{\w}
\wed \ast (\ast \inv d \ast) \w'.
\]
Stokes' Theorem now implies that
\[
\int d\ol{\w} \wed \ast \w' = (-1)^p \int \ol{\w} \wed (\ast \inv d
\ast) \w'.
\]
If we denote
\begin{equation}\label{4:eqn:Hodgecodifferential}
d^*= (-1)^p(\ast \inv d \ast),
\end{equation}
then $\<\w',d\w\>=\<d^*\w',\w\>$. Hence, $d^*$ is the adjoint of
$d$. We call $d^*$ the {\em Hodge codifferential}.

\bigskip

It can be shown that the K\"{a}hler--Dirac operator is equal to $d +
d^*$, see \cite{LAW,VAR} for details. Since $d^2=(d^*)^2=0$, it
holds that $(d+d^*)^2=dd^*+d^*d$. We call $\nabla=dd^*+d^*d$ the
{\em Hodge--Laplacian} operator. If $\nabla(\w)=0$, then we call
$\w$ a {\em harmonic form}. The space of harmonic forms is denoted
by $\Om_\nabla(M)$. An important result involving the differential,
codifferential, and Laplacian is the {\em Hodge Decomposition
Theorem}:
\[
\Om(M)= \Om_\nabla(M) \oplus d(\Om(M)) \oplus d^*(\Om(M)).
\]

\subsection{Spin Manifolds and Dirac Operators}

Let $M$ be an $n$-dimensional Riemannian manifold and let $(S,c)$ be
a spinor module over $M$. Theorem \ref{thm:functorf} implies that,
for any $p \in M$, $c$ induces an action of $\Cl(M)_p$ on $S_p$. A
spinor module for which this action is irreducible, for all $p$, is
called an {\em irreducible} spinor bundle. If the dimension of $M$
is $2k$, or $2k+1$, then equations (\ref{4:eqn:Cliffrep1}) and
(\ref{4:eqn:Cliffrep2}) imply that any irreducible spinor bundle
will have rank $2^k$.

Irreducible spinor bundles do not exist for every Riemannian
manifold. However, there does exist a distinguished type of
Riemannian manifold for which one always does: the Riemannian spin
manifolds. Usually, a {\em spin manifold} is defined to be an
$n$-dimensional orientable Riemannian manifold whose
$SO(n)$-principle bundle of oriented orthonormal frames can be
`lifted' to a $\spin(n)$-principle bundle. (The group $\spin(n)$ is
the universal covering group of $SO(n)$ and it arises as a subspace
of $\Cl(\bC^n)$.) The restriction of the spin representation, given
in equations (\ref{4:eqn:Cliffrep1}) and (\ref{4:eqn:Cliffrep2}), to
$\spin(n)$ gives a representation of $\spin(n)$. The vector bundle
associated to this representation is an irreducible spinor bundle.

The ability to `lift' the bundle of oriented orthonormal frames to a
$\spin(n)$-bundle can be shown to be equivalent to the vanishing of
the second Stiefel--Whitney cohomological class of $T(M)$. This
means that a spin manifold can be alternatively defined as
Riemannian manifold whose second Stiefel--Whitney class vanishes.

A very thorough introduction to the principle bundle approach to
spin manifolds can be found in \cite{LAW}.

\bigskip

When $M$ is \cptn, there also exists an operator theoretic
formulation of the definition. While the idea for this approach
originally came from Connes, it was Roger Plymen \cite{PLYM} who
first published a written account of it. As one would expect, it is
the operator theoretic approach that we shall follow here. We shall
make the assumption that $M$ is of even dimension. The odd
dimensional case follows along very similar lines but it has some
added technical difficulties, for details see \cite{PLYM} or
\cite{VAR}.

\subsubsection{Morita Equivalence}

Let $\B$ be a \calgn, and let $\E$ be a right $B$-module. We call
$\E$ a {\em pre-Hilbert $\B$-module} if there exists a map $(\cdot
\,,\cdot):\E \by \E \to \B$, called the {\em Hilbert module inner
product}, \stn, for all $x,y,z \in \E$, $b \in \B$, $\l \in \bC$,
\begin{enumerate}
\item $(x,y+z)=(x,y)+(x,z)$;

\item $(x,yb)+(x,y)b$;

\item $(x,y)=(y,x)^*$;

\item $(x,x) \geq 0; \text{ and if } (x,x)=0, \text{ then } x=0$.
\end{enumerate}
We can define a norm on $\E$ by setting $\|x\|=\sqrt{\|(x,x)\|}$,
for $x \in \E$. If $\E$ is complete \wrt this norm, then we say that
$\E$ is a {\em Hilbert $\B$-module}. Clearly, every Hilbert space
(with a right linear inner product) is a Hilbert $\bC$-module. Every
\calg $\A$ can be given a Hilbert $\A$-module structure by defining
$(a,b)=a^*b$, $a,b \in \A$. (It is not too hard to see that Hilbert
modules generalise the notion of a Hilbert bundle; for details see
\cite{LAND}.)

A Hilbert $\B$-module is called {\em full} if the closure of the
linear span of ${\{(x,y):x,y \in \E\}}$ is $\B$.

Let $\E$ be a Hilbert module and let $T$ be a module mapping from
$\E$ to $\E$. We say that $T$ is {\em adjointable} if there exists a
module mapping $T^*:\E \to \E$ satisfying $(Tx,y)=(x,T^*y)$, for all
$x,y \in \E$. We call $T^*$ the {\em adjoint} of $T$. It is routine
to verify that if an operator is adjointable, then its adjoint is
necessarily unique, and that $(T^*)^*=T$. Furthermore, a
straightforward application of the closed graph theorem will show
that any adjointable operator is bounded. However, unlike the
special case of Hilbert space operators, not all bounded module maps
are adjointable. We denote the space of adjointable module maps on
$\E$ by $L(\E)$.

For $x,y \in \E$, consider the mapping
\[
\ta_{x,y}:\E \to \E,~~~~~z \mto x(y,z).
\]
It is easy to see that $\ta_{x,y}$ is adjointable, for all ${x,y \in
\E}$, with $\ta^*_{x,y}=\ta_{y,x}$. Let $K(\E)$ denote the closure
in $L(\E)$ of the linear span of ${\{\ta_{x,y}:x,y \in \E\}}$. We
call an element of $K(\E)$ a {\em \cpt operator}.

Let $\A,\B$ be two \calgn s. If there exists a full Hilbert
$\B$-module \st \linebreak ${\A \simeq K(\E)}$, or a full Hilbert
$\A$-module \st $\B \simeq K(\E)$, then we say that $\A$ and $\B$
are {\em Morita equivalent}. We call $\E$ an {\em
$\A$-$\B$-equivalence bimodule}, or a {\em $\B$-$\A$-equivalence
bimodule}, depending on which case holds. (An interesting fact is
that is that if two \algs are Morita equivalent, then their
Hochschild and cyclic (co)homology groups are the same.)

\subsubsection{The Morita Equivalence of $C(M)$ and $\G(\Cl(M))$}

We can endow $\G(\Cl(M))$ with a norm by defining
${\|b\|_\infty=\sup\{b:\|b(x)\|, x \in M\}}$, for ${b \in
\G(\Cl(M))}$; where  $\|\cdot \|$ is the unique $C^*$-norm on
$\Cl(M)_p \simeq M_{2^n}(\bC)$. It is straightforward to show that
$\G(\Cl(M))$ is a \calg \wrt this norm. A {\em spin$^c$ structure}
for $M$ is a pair $(\w,\S)$ consisting of an orientation $\w$, and a
$\G(\Cl(M))$-$C(M)$-equivalence bimodule $\S$. A manifold endowed
with a spin$^c$ structure is called a {\em spin$^c$ manifold}.

It is by no means guaranteed that an arbitrary manifold can be
endowed with a spin$^c$ structure. Those manifolds that can be so
endowed can be categorised in cohomological terms. Specifically, an
oriented Riemannian manifold can be equipped with a spin$^c$
structure \iff the Dixmier--Douady class of its Clifford bundle is
zero. (The Dixmier--Douady class of $\Cl(M)$ is a characteristic
class of $\Cl(M)$ that takes values in $H^3(M,\bZ)$; it is equal to
the third integral Stiefel--Whitney class of $T(M)$; for details see
\cite{PLYM}. This subtle interaction of spin$^c$ structures with the
underlying topology of the manifold is one of the reasons why they
are so interesting.)

Now if $M$ is a spin$^c$ manifold and if $(\w,\S)$ is its spin$^c$
structure, then it can be shown that $\S$ is projective and
finitely-generated. Thus, by the Serre--Swan Theorem, there exists a
\vbd $S$ over $M$ \st $\S=\G(S)$. (It is not too hard to see that
the Hilbert module inner product of $\S$ can be induced by a \vbd
inner product on $S$.) Using Theorem \ref{thm:functorf}, it is easy
to show that $S$ is a spinor bundle. Moreover, using Dixmier--Douady
theory, $S$ can be shown to be an irreducible spinor bundle
\cite{VAR,PLYM}.

\subsubsection{Spin Manifolds and the Atiyah--Singer--Dirac
Operator} 

While it is certainly possible to construct Dirac operators that act
on the smooth sections of the spinor bundle associated to a spin$^c$
manifold, our interest lies in a more specific structure. Before we
introduce this structure, however, we shall need to define a new
operator: Let ${\chi:\G(\Cl(M)) \to \G(\Cl(M))}$ be the unique
linear mapping \stn, for $b \in \G(\Cl(M))$,
\[
\chi(b):p \mto \chi(b(p))
\]
(where on the right hand side $\chi$ is the operator defined in
\mbox{Section $4.1.1$}).

\begin{defn}
A {\em spin structure} on an orientable Riemannian manifold $M$ is a
triple $(\w,\G(S),C)$, where $(\w,\G(S))$ is a spin$^c$ structure on
$M$, and $C$ is a bijective module endomorphism of $\G(S)$ \st
\begin{enumerate}
\item $C(bsf)=\chi(\ol{b})(Cs)\ol{f},~~~~f \in C(M),\, b \in
\G(\Cl(M)),\, s \in \G^\infty(S)$;
\end{enumerate}
and, if $(\cdot\,,\cdot)$ is the Hilbert module inner product of
$S$, then
\begin{enumerate}
\item[2] $(Cs,Cs')=(s',s),~~~~s,s' \in \G(S)$.
\end{enumerate}
A {\em spin manifold} is an orientable \cpt Riemannian manifold
endowed with a spin structure.
\end{defn}

For sake of clarity, it is worthwhile to show exactly what is meant
by $\ol{b}$, the `complex conjugate' of $b \in \G(\Cl(M))$. Since
$T_p(M) \simeq T_p(M;\bR) \oplus iT_p(M;\bR)$, there exists a
canonical complex conjugation operator on $T_p(M)$, defined by
$\ol{(v_p,\l v_p)}=(v_p,\ol{\l}v_p)$, for $v_p \in \Cl(M;\bR),\, \l
\in \bC$. This has a unique extension to an antilinear mapping on
$\T(T_p(M))$, the tensor \alg of $T_p(M)$, which in turn descends to
an antilinear mapping on $\Cl(T_p(M)$. This `complex conjugation' on
$\Cl(T_p(M)$ then induces a complex conjugation on $\G(\Cl(M))$ in
an obvious manner; we denote the image of ${b \in \G(\Cl(M))}$ under
this mapping by $\ol{b}$.

\bigskip

One may well ask why we are interested in the existence, or not, of
a spin structure for a manifold. Unfortunately, it is a little
difficult to give an intuitive way of looking at the operator $C$
without engaging in an excessive digression. We shall, instead,
attempt to justify its introduction as follows. Classically, there
are two principal motivations: firstly, we have the important
properties of the (Atiyah--Singer--)Dirac operator that is
canonically associated to each spin structure. In the general
spin$^c$ case such a well-behaved Dirac operator is not guaranteed
to exist. Secondly, we have the formulation of spin structures in
terms of principle bundles. In this setting the condition
corresponding to the existence of the operator $C$ is much more
natural; for details see \cite{VAR}. Finally, from a noncommutative
point of view, we are interested in spin structures because of
Rennie's Spin Manifold Theorem, as discussed below. A version of
Rennie's Theorem is not known to hold in the general spin$^c$ case.

\bigskip

Let $(\w,\G(S),C)$ be a spin structure for a Riemannian manifold
$M$. It can be shown that there exists a unique connection
$\nabla^S$ for $\G^\infty(S)$ called the {\em spin connection}, \st

\begin{enumerate}
\item $\nabla^S$ commutes with $C$;
\item $\nabla_X^S(c(\w)s)=c(\nabla_X^g \w)s+c(\w) \nabla_X^Ss,~~~~\w \in \Om^1(M),\, s \in \G^\infty(S)$;
\end{enumerate}
and, if $(\cdot\,,\cdot)$ is the Hilbert module inner product of
$S$, then
\begin{enumerate}
\item[3] $(\nabla^S_X s_1,s_2)+(s_1,\nabla_X^S s_2)=X(s_1,s_2),~~~~s_1,s_2 \in
\G^\infty(S)$.
\end{enumerate}

\bigskip

The Dirac operator associated to $S$ and  $\nabla^S$ is called the
{\em Atiyah--Singer--Dirac operator} and it is denoted by $\Dirac$.
It was introduced by Atiyah and Singer while they were working on
their famous index theorem.

The Atiyah--Singer--Dirac operator has many important applications
in modern mathematics and physics. In recent years, for example, it
has been used in the study of $4$-dimensional manifolds through the
Seiberg--Witten invariants. Details on the mathematical applications
of Dirac operators can be found in \cite{LAW,FRIED}; and details on
some of its physical applications can be found in \cite{ESPO,NAKA}.

\subsection{Properties of the Atiyah--Singer--Dirac Operator}

Let $S$ be the irreducible spinor bundle associated to a spin
manifold $M$ and let $(\cdot,\cdot)$ be the Hilbert module inner
product of $S$. If $d\mu$ is the Riemannian measure on $M$, then we
can define an inner product on $\G^\infty(S)$ by
\begin{eqnarray*}\label{4:defn:spininnerprod}
\<s_1,s_2\>=\int (s_2,s_1)d\mu,~~~~s_1,s_2 \in \G^\infty(S).
\end{eqnarray*}
We denote by $L^2(S)$ the Hilbert space completion of $\G^\infty(S)$
\wrt the norm induced by this inner product. We call $L^2(S)$ the
{\em Hilbert space of square-integral spinors} on $M$. In this
context $\Dirac$ becomes a linear operator defined on a dense
subspace of $\H$. It can be shown that $\Dirac$ is always unbounded.

\bigskip

If $A$ is an operator on a Hilbert space $H$, then the {\em graph}
of $A$ is the set
\[
G(A) = \{(x, A x): x \in \dom(A)\} \subseteq H \oplus H.
\]
If $\<\cdot,\cdot\>$ is the inner product of $H$, then one can
define an inner product on $H \oplus H$ by setting
\[
\<(x,y),(u,v)\>=\<x,u\>+\<y,v\>.
\]
If the closure of $G(A)$ \wrt the norm induced by $\<\cdot,\,
\cdot\>$ is also the graph of an operator $B$, then we say that $A$
is closable, and we call $B$ the {\em closure} of $A$. Obviously,
the closure of $A$ will extend $A$.

It can be shown that $\Dirac$ is a closable operator \cite{VAR}.
From now on, we shall always use $\Dirac$ to denote the closure of
the Dirac operator, and we shall refer to the closure of the Dirac
operator simply as the Dirac operator.

The theory of unbounded operators is notoriously problematic, to the
extent that substantial results about general unbounded operators
are rare. It is only when an unbounded operator is closed (that is,
equal to its closure) that it becomes somewhat `manageable'.

Given a densely defined linear operator $A$ on a Hilbert space $H$,
its {\em adjoint} $A^*$ is defined as follows: the domain of $A^*$
consists of all vectors $x \in H$ such that the linear map
\[
\dom(A) \to \bC, ~~~~~ y \mapsto \langle x , A y \rangle,
\]
is a bounded linear functional. Since the $\dom(A)$ is dense in $H$,
each such functional will extend to a unique bounded linear
functional defined on all $H$. Now if $x$ is in the domain of $A^*$,
then the Riesz representation theorem implies that there is a unique
vector $z \in H$ such that
\[
\langle x , A y \rangle  = \langle z   , y \rangle, \quad \text{ for
all } y \in \dom(A).
\]
It is routine to show that the dependence of $z$ on $x$ is linear.
We define $A^*$ to be the unique linear operator for which $A^*x=z$.
An operator $A$ is said to be {\em self-adjoint} if $G(A)=G(A^*).$
Self-adjoint operators have many very useful properties that are not
necessarily possessed by non-self-adjoint operators. For example, it
can easily be shown that every self-adjoint operator is closed.

In $1973$ Wolf \cite{WOLF} showed that $\Dirac$ is a self-adjoint
operator.

\bigskip

We can regard any $f \in \csm$ as a linear operator on $\G^\infty(S)
\sseq L^2(S)$ that acts by multiplication. This representation of
$\csm$ is obviously faithful. Furthermore, the definition of the
inner product on $\G^\infty(S)$ implies that $f$ is bounded with
norm $\|f\|_\infty$. Hence, it extends to a unique bounded linear
operator on all of $L^2(S)$ with norm $\|f\|_\infty$. We shall not
distinguish notationally between $f$ and this operator.

Now
\[
[\Dirac,f] s = \hat{c} (\nabla ^S (s f)) - f\hat{c}(\nabla^S s) =
\hat{c} ( \nabla ^S (s f) - (\nabla ^Ss)f ) = \hat{c} (s \oby df) =
c(df)s,
\]
for all $f \in C^{\infty}(M)$, $s \in \G^{\infty}(S)$. Thus, we have
that
\begin{equation}\label{eqn:4:Diraccomm}
[\Dirac,f]= c(df),
\end{equation}
for all $f \in \csm$. As we shall see, this formula allows one to
recover the differential structure of $M$ from $\Dirac$. Moreover,
it can easily be used to show that, for all $f \in \csm$, the
operator norm of the densely defined operator $[\Dirac,f]$ is equal
to $\|(df)^\sharp\|_{\infty}$; where $\|\X\|_{\infty}=\sup_{p \in
M}g_p(X(p),X(p))$, for $X \in \X(M)$. Consequently, each such
operator is bounded and has a unique extension to a bounded linear
operator defined on all of $L^2(S)$. We shall not distinguish
notationally between $[\Dirac,f]$ and this operator.

\bigskip

Finally, we shall list one more important property of $\Dirac$: the
Dirac operator has {\em \cpt resolvent}; that is, $(\Dirac-\l)\inv$
is a \cpt operator, for all \linebreak ${\l \in
\r(\Dirac)=\bC/\s(\Dirac)}$.

\section{Spectral Triples}

As we stated earlier, and as we shall see below, much of the
structure of a \cpt Riemannian spin manifold can be expressed in
terms of its associated Dirac operator. This fact motivated Connes
to try and construct a \nc generalisation of spin manifolds based on
a type of generalised Dirac operator. In Connes' work the Dirac
operator is no longer an object associated to a manifold but rather
one of the data defining it.

\begin{defn} A {\em spectral triple} $(A,H,D)$ consists of a $*$-\alg
$A$ faithfully represented on a Hilbert space $H$, together with a
(possibly unbounded) self-adjoint operator $D$ on $H$ such that:
\begin{enumerate}
\item $\dom(D) \sseq H$ is a dense subset of $H$, and $a\,\dom(D)
\sseq \dom(D)$, for all $a \in A$;
\item the operator $[D,a]$ is bounded on $\dom(D)$, for all $a \in A$,
and so it extends to a unique bounded operator on $H$;
\item $(D-\l)^{-1}$ is a {\em compact operator}, for all  $\l \notin
\s(D)$.
\end{enumerate}
\end{defn}

(Since no confusion will arise, we shall not distinguish
notationally between an element of $A$ and its image in $B(H)$, nor
between $[D,a]$ and its unique extension. Moreover, when we speak of
$\ol{A}$, the closure of $A$, we mean the closure of its image in
$B(H)$ \wrt the operator norm.)

If $M$ is a \cpt Riemannian spin manifold, then, from our comments
above, $\cantrip$ is a spectral triple; we call it the {\em
canonical triple} associated to $M$.

The definition of a spectral triple is partly motivated by the
notion of a Fredholm module. We shall not enter into a discussion of
such involved topics here, instead, we refer the interested reader
to \cite{CON,VAR}.

\bigskip

Associated to every spectral triple $(A,H,D)$ is a distinguished
differential calculus. Let $\Om_u(A)$ denote the universal calculus
of $A$, and let $\pi$ denote the unique \alg homomorphism from
$\Om_u(A)$ to $B(H)$ for which
\[
\pi(a_0da_1 \cdots da_n)=a_0[D,a_1]\cdots [D,a_n].
\]
(The fact that $\pi$ is a homomorphism easily follows from the fact
that $d$ and $[D,\cdot]$ are derivations.) Now $\pi(\Om_u(A)) \simeq
\Om_u(A)/\ker(\pi)$ is canonically a graded \algn, and it would be
natural to define a differential $d$ on it by $d(\pi(\w))=\pi(d\w)$.
However, this is not always a well-defined mapping since there can
exist forms $\w \in \Om_u(A)$ \st $\w \in \ker(\pi)$ and $d\w \notin
\ker(\pi)$.

With a view to overcoming this problem, let us consider the sub\alg
\linebreak ${J=\ker(\pi)+d\ker(\pi)}$. If $a_0+da_1 \in J \cap
\Om^k_u(A)$ and $b \in \Om_u(A)$, then
\[
(a_0 + da_1)b = a_0b + (da_1)b = a_0b + d(a_1b) - (-1)^{k-1}a_1db.
\]
Since $a_0b$, $a_1b$, and $a_1db$ are elements of $\ker(\pi)$, we
have that $(a_0 + da_1)b \in J$. Similarly, $b(a_0+da_1) \in J$.
Hence, $J$ is an ideal of $\Om_u(A)$. Moreover, since
${dJ=d(\ker(\pi)) \sseq J}$, $J$ is a differential ideal. This means
that we can canonically give the quotient \alg
\[
\Om_D(A)=\pi(\Om_u(A))/\pi(J)=\pi(\Om_u(A))/\pi(d(J))
\]
the structure of a differential \algn; we call it the {\em
differential \alg of $D$-forms}. Furthermore, since
$\Om^0_D(A)=\pi(\Om_u^0(A))/\pi(J_0)$, and $J_0=\{0\}$, we have that
$\Om^0_D(A)=A$. Hence, $\Om_D(A)$ is a differential calculus over
$A$.

\bigskip

Now if $M$ is a  Riemannian spin manifold and $f_0d_uf_1 \cdots
d_uf_n \in \Om_u^n(\csm)$, then, since $[\Dirac,f]=c(df)$, for all
$f \in \csm$,
\begin{eqnarray*}
\pi(f_0d_u f_1 \cdots d_uf_n) & = & f_0[\Dirac,f_1] \cdots [\Dirac,f_n]\\
                              & = & f_0c(df_1) \cdots c(df_n).
\end{eqnarray*}
One can build upon this fact to show that $\Om_{\Dirac}^p(\csm)
\simeq \Om^p(M)$, for all $p \geq 0$. Moreover, one can show that
$(\Om_{\Dirac}^p(\csm),d)$ and $(\Om(M),d)$ are isomorphic as
differential \algn s. Thus, the differential calculus of a
Riemannian spin manifold is entirely encoded in its canonical
spectral triple.

\bigskip

\subsubsection{Rennie's Spin Manifold Theorem}

A natural question to ask is which spectral triples arise from spin
manifolds and which do not. With a view to answering this question
Connes introduced a refinement of the notion of a spectral triple
called a \ncgn. (The definition of a \ncg is rather lengthy, and we
shall not go into the details here. It is, however, interesting to
note that it contains \nc generalisations of orientability and
Poincar\'{e} duality.) Connes showed that every spectral triple
arising from a \cpt spin manifold is a \ncgn. He went on to claim
that all commutative \nc geometries (that is,  all \nc geometries
\st the $*$-\algs of their underlying spectral triples are
commutative) arise from \cpt spin manifolds \cite{CON2}; Rennie and
Varilly \cite{RENNIE} would later prove this. At present, there is
no extension of this work to the locally \cpt case.

\subsubsection{The Connes--Moscovici Index Theorem}

From what we have presented above, one might get the impression that
spectral triples and cyclic cohomology have little in common.
However, the two areas are intimately linked. In fact, in
\cite{CONMOS2} Connes and Moscovici used spectral triples to
construct new versions of the Chern--Connes maps. Their
reformulation is extremely attractive for practical calculations. In
a related development, they also formulated a very important \nc
generalisation of the Atiyah--Singer index theorem; it is called the
Connes--Moscovici index theorem.

\subsection{The Noncommutative Riemannian Integral}
\label{sec:4:ncint}

Let $M$ be an $n$-dimensional manifold and let $(U,\f)$ be one of
its coordinate \nbdn s. For $i=1,2,\ldots , n$, let the \fn $x_i$ be
some global extension of $\pi_i \circ \f \inv$, where $\pi_i$ is the
canonical projection onto the $i^{\rm th}$ coordinate of $\bR^n$. We
can define a matrix-valued \fn on $U$ by $g_U(p) =
[g_p(dx_i,dx_j)(p)]_{ij}$ (we shall assume that ${\sqrt{\det
(g_U(p))} \geq 0}$, for all $p \in U$, and for all coordinate \nbds
$U$). Now, it is straightforward to to show that there exists a
unique $d\mu \in \Om^n(M)$ \stn, for any coordinate \nbd $U$,
\[
d\mu(p)= \sqrt{\det (g_U(p))}dx_1 \wed dx_2 \wed \cdots \wed
dx_n,~~~~\text{for all $p \in U$}.
\]
We call $d\mu$ the {\em Riemannian volume form}. For any $f \in
\csm$, we call $\int f d\mu$ the {\em Riemannian integral} of $f$
over $M$.

\bigskip

A major achievement of Connes has been to re-express the Riemannian
integral of a spin manifold purely in terms of its Dirac operator.
More explicitly, Connes established the following equation:
\[
\frac{1}{c_n} \int f d\mu = \text{Tr}_\w (f |D|^{-n}),~~~~ ~~~f \in
\csm;
\]
where $c_n$ is a constant, and $\text{Tr}_\w$ denotes the Dixmier
trace. (A precise exposition of this equation is well outside the
scope of this presentation. Instead, we refer the interested reader
to \cite{CON,VAR,ALBERT}.) Using this formula Connes was able to
introduce a generalised \nc Riemannian integral for certain suitable
types of spectral triples.

One of Connes' principal uses for his \nc integral has been to
construct \nc action functionals. Using one of these, he
reconstructed the standard model of particle physics in terms of a
distinguished non-classical spectral triple. While Connes'
reformulation does not address any questions of renormalisation (see
Chapter $5$ for a discussion of renormalisation), it does have far
greater simplicity than the usual presentation. For further details
see \cite{CON,CONLOT}.

The \nc integral has also found use in the study of fractals, see
\cite{GUIDOFRAC}.

\subsection{Connes' State Space Metric}

Let $c:[a,b] \to M$ be a smooth curve in a connected Riemannian
manifold $M$, and let $\dot{c}(t)$ denote the tangent vector to $c$
at $t \in [a,b]$. We can define a nonnegative \cts \fn $s_c:[a,b]
\to \bR$ by
\[
s_c(t)=\sqrt{g_{c(t)}(\dot{c}(t),\dot{c}(t))},~~~~t \in [a,b].
\]
The nonnegative real number
\[
L(c)=\int_a^bs_c dt
\]
is called the {\em arc-length} of $c$. If $c$ is a piecewise smooth
curve, then its {\em arc length} is defined to be the sum of the arc
lengths of its components.

For any two points $p,q \in M$, let us use the symbol $\Om(p,q)$ to
denote the set of all piecewise smooth curves in $M$ from $p$ to
$q$. It is easily shown that $\Om(p,q)$ is non-empty. The
non-negative real number
\[
d(p,q)= \inf \{L(c):c \in \Om(p,q)\},
\]
is called the {\em Riemannian distance} from $p$ to $q$. It can be
shown that the mapping
\[
d:M \by M \to \bR,~~~~(p,q) \to d(p,q),
\]
is a metric on $M$; it is called the {\em Riemannian metric}. An
important result of Riemannian geometry is that the topology
determined by this metric coincides with the original topology of
$M$.

\bigskip

If $M$ is a spin manifold and $\Dirac$ is its associated Dirac
operator, then Connes \cite{CON} showed that
\begin{eqnarray}
\label{eqn2:Diracgeodesic} d(p,q)=\sup\{|f(p)-f(q)|:f \in C(M),
\|[\Dirac,f]\| \leq 1 \}.
\end{eqnarray}
Thus, the Riemannian distance can be recovered from the canonical
spectral triple of $M$.

Let $M$ be a Riemannian spin manifold and let $\cantrip$ be its
canonical spectral triple. We can define a metric on
$\Om(\ol{\csm})$ by
\begin{equation} \label{4:eqn:charmetric}
d(\f,\ps)=\sup\{|\f(f)-\ps(f)|:f \in C(M), \|[\Dirac,f]\| \leq 1\}.
\end{equation}
equation (\ref{eqn2:Diracgeodesic}) and Theorem
\ref{1:thm:cptsurjective} imply that $M$ and $\Om(\ol{\csm})$ are
identical as metric spaces. (We note that, since the representation
of $\csm$ on $L^2(S)$ is isometric, the closure of $\csm$ \wrt the
supremum norm is isometrically isomorphic to the closure of the
image of $\csm$ in $B(L^2(S))$ \wrt the operator norm. Thus, no
ambiguity arises when we speak of the closure of $\csm$.)

\bigskip

If $(A,H,D)$ is a spectral triple, then the metric in
(\ref{4:eqn:charmetric}) can easily be generalised to a metric on
$\Om(\ol{A})$. As we noted in Chapter $1$, however, $\Om(\ol{A})$ is
not guaranteed to be non-empty when $\ol{A}$ is \ncn. Thus, this
metric is of limited interest.

Recall that a positive linear functional of norm $1$ on a \calg $\A$
is called a state. We denote the set of all states on $\A$ by
$S(\A)$, and we call it the {\em state space} of $\A$. The state
space of a \calg is always non-empty. Now since every \pve element
of a \calg $\A$ is of the form $aa^*$, for some $a \in \A$, and
every character on $\A$ is Hermitian and multiplicative, it holds
that every character is positive. If we also recall that every
character is of norm $1$, then we see that $\Om(\A) \sseq S(\A)$.
State spaces are of great importance in the study of \calgn s and in
quantum mechanics.

Now the metric defined above in (\ref{4:eqn:charmetric}) can easily
be modified to define a metric on $S(\ol{A})$; it is called {\em
Connes' state space metric}. While Connes did not explore this
metric very much for the \nc case, Marc Rieffel used it as his
starting point when he was developing the theory of \cqmsn s. We
shall present the theory of \cqmsn s in detail in Chapter $6$.

\section{Dirac Operators and Quantum Groups}

As we saw in Chapter $3$, the boundary between cyclic (co)homology
and quantum groups is a very active area of research. The same is
true of the boundary between spectral triples and quantum groups. At
present there is no comprehensive theory linking Connes' calculi to
covariant differential calculi. While examples of spectral triples
have been constructed over the coordinate \alg of $SU_q(2)$
\cite{CHAKRA,CHAKRA2}, it is known that the basic examples of
calculi over $SU_q(2)$ can not be realized by spectral triples, see
\cite{SCHMUGDEN}.

With a view to better understand the relationship between the two
theories, a lot of effort has been put into constructing generalised
Dirac operators on quantum groups. One approach is to introduce
quantum analogues of Clifford and spinor bundles and to define a
Dirac operator in terms of them, see \cite{KRAM,SIT2} for example.
Here, however, we shall present the approach pursued by Kustermans,
Murphy, and Tuset in \cite{KUSMURTUS0}. They  generalised the
K\"{a}hler--Dirac operator $d+d^*$ by introducing a generalised
codifferential. This generalisation has lead to some very
interesting results, most notably a generalised version of the Hodge
decomposition. We shall  briefly discuss why their work lead them to
call for a generalisation of spectral triples.

\subsection{The Dirac and Laplace Operators}

Let $G$ be a \cqg and let $(A,\DEL)$ be the Hopf $ *\text{-\alg}$
underlying it. A left-covariant $*$-differential calculus over $G$
is a \mbox{$*$-differential} calculus $\Om$ over $A$ \st there
exists an \alg homomorphism ${\DEL_\Om:A \to A \oby \Om}$ that
satisfies $\DEL_\Om(a)=\DEL(a)$, for all $a \in A$; and  ${(\id_A
\oby d)\DEL_\Om = \DEL_\Om d}$. Now the multiplication of $A \oby
\Om$ respects the grading given by
 ${A \oby \Om = \bigoplus_{n=0}^\infty A \oby \Om^n}$. Hence,
$A \oby \Om$ can be considered as a graded \algn. Moreover, since
\[
\DEL_\Om(a_0da_1 \cdots da_n)=\DEL(a_0)[(1 \oby d)\DEL(a_1)] \cdots
[(1 \oby d)\DEL(a_n)] \in A \oby \Om^n,
\]
we have that $\DEL_\Om$ is a mapping of degree zero \wrt this
grading.

\bigbreak

We say that $\w \in \Om$ is  {\em left-invariant} if $\DEL(\w)=1
\oby \w$. Clearly, the set of all left-invariant elements forms a
sub\alg of $\Om$; we denote it by $\Om_\invt$. It is obvious that
$\Om_\invt^k=\Om^k \cap \Om_\invt$ is also a sub\algn. We say that a
differential calculus over a \cqg is {\em strongly
finite-dimensional} if $\Om_\invt$ is finite-dimensional as a linear
space. Using results from general Hopf \alg theory, it can be shown
that every strongly finite-dimensional calculus is
finite-dimensional as a graded \algn. It is easy to show that in the
classical case left-invariant elements correspond to left-invariant
differential forms, as defined in the next chapter. As is well known
and easily seen, the \alg of differential forms of every Lie group
is strongly finite-dimensional.

Now if $\w = \sum_{k=1}^n \w_k$ is a left-invariant element of a
strongly finite-dimensional calculus $\Om$, and each $\w_k \in
\Om^k$, then we have that $\sum_{k=1}^n 1 \oby \w_k = \sum_{k=1}^n
\DEL(\w_k)$. Since $\DEL_\Om$ is a mapping of degree $0$, it follows
that $\DEL(\w_k)=1 \oby \w_k$, for ${k=1, \ldots ,n}$. Therefore,
${\Om_\invt=\bigoplus_{k=0}^\infty \Om_\invt^k}$. Furthermore, the
fact that $(\id_A \oby d)\DEL_\Om = \DEL_\Om d$, implies that
$\Om_\invt$ is invariant under the action of $d$. Hence,
$(\Om_\invt,\del)$ is a graded differential \algn, where $\del$
denotes the restriction of $d$ to $\Om_\invt$.

Using results from general Hopf \alg theory again, it can be shown
that $A \oby \Om_\invt$ and $\Om$ are isomorphic as left
$A$-modules. An isomorphism is provided by the unique mapping that
sends $a \oby \w$ to $a \w$. Clearly, we also have that $A \oby
\Om_\invt^k \simeq \Om^k$.

\bigskip

Inner products on $\Om$ will play an important part in our work. We
call an inner product on $\Om$ {\em graded} if the subspaces $\Om_k$
are orthogonal \wrt it. We call an inner product $\<\cdot,\cdot\>$
{\em left-invariant} if $\<a\w_1,b\w_2\>=h(ba^*)\<\w_1,\w_2\>$, for
all $a,b \in A$, and $\w_1,\w_2 \in \Om$; where $h$ denotes the
restriction of the Haar integral of $G$ to $A$. We shall only
consider graded left-invariant inner products here. Clearly, the
inner product defined in equations (\ref{4:eqn:hodgeinnerproduct1})
and (\ref{4:eqn:hodgeinnerproduct2}) is graded and left invariant.

We shall now introduce some useful linear operators. If $\chi$ is a
linear functional on the Hopf \alg $A$, and $a \in A$, then we write
$\chi * a$ for $(\id \oby \chi)\DEL(a)$; and we denote by $E_\chi$
the linear operator defined by setting $E_\chi(a)=\chi * a$. We
define $\chi^*$, the {\em adjoint} of $\chi$, by setting
$\chi^*(a)=\ol{\chi(S(a)^*)}$; where $S$ is the antipode of $A$.
And, if $\w \in \Om_\invt$, then we denote by $M_\w$ the linear
operator on $\Om_\invt$ defined by setting $M_\w(\eta)=\w\eta$,
$\eta \in \Om_\invt$.

Now if $\Om$ is endowed with an inner product, then it can be shown
(see \cite{KUSMURTUS0} and references therein) that $E_\chi$ is
adjointable \wrt this inner product and that $E^*_\chi=E_{\chi^*}$.
It can also be shown that if $\w_1,\dots,\w_m$ is an orthonormal
basis for $\Om_\invt^1$, then there exist unique linear functionals
${\chi_1,\dots,\chi_m}$ on $A$ such that
\begin{equation} \label{eqn: d-chi equation}
da={\sum_{r=1}^m E_{\chi_r}(a)\w_r},
\end{equation}
for all $a \in A$, (see \cite{SCHMU} for details).

\begin{thm} \label{thm: d-adjointability theorem}
Let $(\Om,d)$ be a strongly finite-dimensional $*$-differential
calculus over a \cqg $G$, and let $(A,\DEL)$ be the Hopf $*$-\alg
underlying $G$. If $\Om$ is endowed with an inner product, then $d$
is adjointable, and $d^*$ is of degree $-1$. Indeed, if
$\w_1,\dots,\w_m$ and $\chi_1,\dots,\chi_m$ are as in
equation~(\ref{eqn: d-chi equation}), then
\begin{equation} \label{eqn: tensor d-equation}
d=\id_A\otimes \del+\sum_{j=1}^m E_{\chi_j}\otimes
M_{\w_j}\end{equation} and
\begin{equation} \label{eqn: tensor d*-equation}
d^*=\id_A\otimes \del^*+\sum_{j=1}^m E_{\chi_j^*}\otimes
M_{\w_j}^*.\end{equation}
\end{thm}

\demo Recall that we can identify $\Om$ and $A \oby \Om_{\invt}$ by
identifying $a\w$ with $a \oby \w$, for $a \in A$, $\w \in
\Om_{\invt}$. If $a\in A$, and $\w \in \Om_\invt$, then since
\[
{d(a\w)=(da)\w+ad\w}{=\sum_{j=1}^m E_{\chi_j}(a)\w_j\w}+ad\w,
\]
we have that $d={\id_A\otimes \del}+{\sum_{j=1}^m E_{\chi_j}\otimes
M_{\w_j}}$. As remarked above, the operators $E_{\chi_j}$ are
adjointable; since $\Om_\invt$ is finite-dimensional, $M_{\w_j}$ and
$\del$ are also adjointable. The adjointability of $d$ and the
formula for $d^*$ in the statement of the theorem follow
immediately. Since $d$ is of degree $1$ and the inner product is
graded, it is easily verified that $d^*(a)=0$, for all $a \in A$,
and if $\w_k \in \Om^k$ and $k \geq 1$, then $d^*(\w_k) \in
\Om^{k-1}$.

The operator $d^*$ is called the {\em codifferential} of $d$; the
sum $D=d+d^*$ is called the {\em Dirac operator}; and the square
$\nabla=(d+d^*)^2$ is called the {\em Laplacian}. Since $d^2=0$
implies that $d^{*2}=0$, we have that $\nabla=dd^*+d^*d$. We call
$\w \in \Om$ a {\em harmonic form} if $\nabla(\w)=0$; we denote the
linear space of harmonic forms by $\Om_\nabla$. Clearly, these
operators generalise the classical K\"ahler--Dirac and
Hodge--Laplacian operators (or, more correctly, it generalises its
restriction to $\Om_{\Pol}(M)$).

\subsection{The Hodge Decomposition}

Now that we have defined generalised Dirac and Laplacian operators,
we are ready to present the following important theorem. It contains
a generalisation of the Hodge decomposition, which is surprising,
since such a decomposition is known not to exist in the setting of
general noncommutative  geometry. (Recall that an operator $T$ on a
linear space $V$ is {\em diagonalisable} if $V$ admits a Hamel basis
consisting of eigenvectors of $T$.)

\begin{thm} [Hodge decomposition \cite{KUSMURTUS0}] Let $(\Om,d)$ be a
strongly  finite-\newline dimensional $*$-differential calculus over
a compact quantum group $G$ and suppose that we have a graded, left
invariant, inner product on $\Om$. Then

\begin{enumerate}

\item The Dirac operator $D$ and the Laplacian $\nabla$ are
diagonalisable;

\item The space $\Om$ admits the orthogonal decomposition
\[
\Om= \Om_\nabla \oplus d(\Om) \oplus d^*(\Om).
\]
\end{enumerate}
\end{thm}

\subsection{The Hodge Operator}

Let $(\Om,d)$ be a strongly finite-dimensional differential calculus
over a \cqg $G$, and assume that $(\Om,d)$ is endowed with a graded,
left-invariant, inner product. As we remarked earlier, the fact that
$\omd$ is strongly finite-dimensional implies that $\omd$ is
finite-dimensional; let us assume that it is of dimension $n$. It is
straightforward to show that $\Om_{\invt}^n$ contains a self-adjoint
element $\ta$ of norm $1$. If $\dim (\Om_\invt^n)=1$, and we shall
always that it is, then this element is unique up to a change of
sign. (The assumption that $\dim (\Om_\invt^n)=1$ is motivated by
the fact it holds in the classical case.) Indeed, since
$A\Om^n_\invt=\Om^n$, it holds that $\{a \ta : a \in A\} = \Om^n$.
Let us define a linear functional $\int$ on $\Om^n$ by setting $\int
\w = h(a)$, when $\w = a \ta$, for some $a \in A$. We call $\int$
the {\em integral} associated to  $\ta$. Now we say that a calculus
$(\Om',d')$ is {\em non-degenerate} if, whenever $k=0, \ldots ,n$
and $\w\in \Om'_k$, and $\w'\w=0$, for all $\w'\in \Om'_{n-k}$, then
necessarily $\w=0$. If we further assume that $(\Om,d)$ is
non-degenerate and that $d(\Om_\invt^{n-1})=\{0\}$, then it can be
shown that $\int$ is a closed twisted graded trace, for details see
\cite{KUSMURTUS}. (Again, these assumptions are motivated by the
classical case.)

A proof of the following important theorem can be found in
\cite{KUSMURTUS0}.

\begin{thm} \label{thm: existence of Hodge operator}
Suppose that $\omd$ is a non-degenerate. Then there exists a unique
left $A$-linear operator $L$ on $\Om$ such that
$L(\Om^k)=\Om^{n-k}$, for ${k=0,\dots,n}$, and such that $\int
\w^*L(\w')=\<\w',\w\>$, for all $\w,\w'\in \Om$. Moreover, $L$ is
bijective.
\end{thm}

The operator $L$ is  called the {\em Hodge operator} on $\Om$; we
see that it generalises the classical Hodge operator defined in
Section \ref{sec:4:geomdirac}. Moreover, since
$L(\Om_\invt^k)=\Om_\invt^{n-k}$, for $k=0, \cdots n$, we have
$\dim(\Om_\invt^k)=\dim(\Om_\invt^{n-k})$.

The following pleasing result for the codifferential generalises
equation (\ref{4:eqn:Hodgecodifferential}) to the quantum setting.
\begin{cor}
Suppose $\omd$ is non-degenerate and $d(\Om^\invt_{n-1})=\{0\}$.
Then, if $k=0, \cdots n$ and $\w \in \Om_k$, we have
\[
d^*\w=(-1)^k L^{-1}d L(\w).
\]
\end{cor}

\demo If $k=0$, then clearly $d^*\w=0$. Also, $d L(\w)=0$, since
$L(\w)$ is in $\Om_n$ and $d(\Om_n)=0$. Hence, $d^*\w=(-1)^kL^{-1}d
L(\w)$ in this case.

Suppose now that $k>0$. Clearly, $L^{-1}d L(\w)\in \Om_{k-1}$.
Hence, if ${\w'\in \Om_{k-1}}$, then
\[
\begin{tabular}{ll}
$\<(-1)^kL^{-1}d L(\w),\w'\>$ & $={(-1)^k\int\w^{'*}d L(\w)}= {\int
d(\w')^*L(\w)}$\\       & $=\<\w,d(\w')\>=\<d^*\w,\w'\>$;
\end{tabular}
\]
(note that we have used the fact that $d(\Om^{\invt}_{n-1})=0$
implies that $\int$ is closed). Thus, $d^*\w=(-1)^kL^{-1}d L(\w)$ in
this case also. \qed

\subsection{A Dirac Operator on Woronowicz's Calculus}

As we mentioned in Chapter $3$, the first ever example of a \nc
differential calculus was Woronowicz's $3$-dimensional calculus over
$SU_q(2)$. It is strongly finite-dimensional, and in
\cite{KUSMURTUS0} Kustermans, Murphy, and Tuset calculated the
eigenvalues of the Dirac and Laplace operators arising from a
canonical choice of inner product for this calculus. From their
investigations they feel that it is very unlikely that this Dirac
operator fits into the framework of Connes' spectral triples;
although they do not have a proof of this fact. This provides an
indication that generalisations of spectral triples will have to be
studied. Their Dirac operator does, however, have properties that
are close to those required to enable it to fit into A.~Jaffe's
\cite{JAFFE} extension of Connes' theory. It may be that Jaffe's
theory can, and should, be further developed to cover this example.

\bigskip

A very interesting recent development in this area is Connes and
Moscovici's introduction of a twisted version of the spectral triple
definition \cite{CONMOS3}. Let $(A,H,D)$ be a triple consisting of
an \alg $A$, a Hilbert space $H$ upon which $A$ is represented, and
a self-adjoint operator $D$ with compact resolvent. If $\s$ is an
\alg automorphism of $A$ \st $Da - \s(a)D$ is bounded, for all  $a
\in A$, then we say that $(A,H,D)$ is a {\em $\s$-spectral triple}.
Clearly, if $\s=\id_A$, then we recover the ordinary spectral triple
definition. Connes and Moscovici believe that, since the domain of
quantum groups is an arena where twisting frequently occurs,
$\s$-spectral triples could be useful in the study of quantum
groups. It is interesting to note that they cite Hadfield and
Kr\"ahmer's paper \cite{HADKRA2}. The obvious question to ask is
whether or not Kustermann, Murphy, and Tuset's Dirac operator fits
into this framework. At present no work has been done in this area.

\chapter{Fuzzy Physics}

In this chapter we shall give a brief account of some of the work
being done by the Dublin node of the European Union Operator
Algebras Network: the Dublin Institute for Advanced Studies (DIAS).
The members of the DIAS are pursuing research in an area called
fuzzy physics. Essentially, fuzzy physics is an application of \ncg
to the problem of $UV$-divergences in \qftn. We shall begin this
chapter by introducing \qftn, we shall then discuss the role \ncg
plays, and finally we shall discuss some of the DIAS work.

Unfortunately, because of the physical nature of the material, our
presentation in this chapter will sometimes be a little vague.
Moreover, we shall also assume some prior knowledge of physics.

\section{Quantum Field Theory}

As we saw in Chapter $4$, Dirac's relativistic equation of the
electron  introduces into physics the notion of antimatter. The key
property of  antimatter is that a particle and its antiparticle can
come together and annihilate one another, their combined mass being
converted into energy in accordance with Einstein's equation:
$E=mc^2$. Conversely, if sufficient energy is introduced into a
system, localised in a suitably small region, then there arises the
strong possibility that this energy might serve to create some
particle together with its antiparticle. This does not violate any
conservation laws since the conserved quantum numbers of a particle
and its antiparticle have opposite signs. For example, for an
electron-positron pair, the electron has an electric charge of $-1$
and the positron has an electric charge of $1$, thus, the addition
of the pair has no effect on the total charge of the system. All
this means that the number of particles in a system is always
variable.

However, ordinary  non-relativistic quantum mechanics does not allow
us to describe systems in which the number of particles is variable.
In order to cope with the introduction of anti-matter a new
framework was required. The theory that emerged was called {\em
\qftn}. There exist a number of different approaches to the subject.
Any reader interested in the details of \qft can find a good, well
referenced, non-technical introduction in \cite{PEN}.

\bigskip

The laws of classical mechanics were formulated by Isaac Newton in
$1680$'s. In the centuries that followed two important
reformulations emerged: one due to Joseph Louis Lagrange, and
another due to William Rowan Hamilton. The standard approach to
dynamics in quantum mechanics is based on Hamilton's reformulation.
In the $1940$'s the brilliant American physicist Richard Feynman
introduced an approach that is based on Lagrange's reformulation. It
is based on the idea that if a particle travels between two points,
then, in a certain `quantum mechanical sense', it must travel every
path between those two points. The basic structure in Feynman's
approach is the Feynman path integral. Feynman used his ideas to
formulate an approach to \qftn. At present, it is arguably the most
widely used framework.

Soon after its introduction, however, it became clear that Feynman's
approach suffers from serious conceptual and technical difficulties.
From a mathematical point of view the whole construction is horribly
ill-defined. From a practical point of view physical processes
occurring at arbitrarily small distances cause the theory to give
infinite answers to questions that should have a finite answers;
these are the so called {\em UV-divergences}. In the late $1940$'s,
the ground breaking work of Feynman, Tomonaga, Dyson, Swinger and
others on these divergence problems produced successful, if somewhat
ad hoc, methods for extracting information from the theory. Their
approach is generally known as {\em renormalisation}. When
renormalisation  was applied to systems involving the
electromagnetic field, the result was quantum electrodynamics (QED).
Astonishingly, the predictions of QED match experimental data to a
level of accuracy never seen before. QED has been described as the
`most successful scientific theory ever'. Later, a quantum field
theory that united the electromagnetic, strong nuclear, and weak
nuclear fields would emerge; it is called the standard model. It is
regarded as a high point in twentieth century physics. However,
despite all efforts, gravity (the last remaining fundamental field)
has not been incorporated into this framework.

\section{Noncommutative Regularization}

Despite the experimental successes of quantum field theory, it is
still wholly unsatisfactory from a mathematical point of view.
Moreover, renormalization (which is in essence a perturbation
theory) is not a `fundamental' process (even Feynman himself was of
this opinion). The fact that physics has failed to incorporate
gravity into the framework of the standard model seems to confirm
that fundamental changes need to be made. Study in this area is
known as {\em nonpertubative regularization}.

At present, the conventional approach to nonperturbative
regularization is {\em lattice field theory}: it works by replacing
a \cts manifold with a discrete lattice. In this context Feynman's
formulation can be well defined mathematically, and the divergences
that arise from processes occurring at arbitrarily small distances
disappear. There is, however, one feature that must be criticized:
lattices do not retain the symmetries of the original theory, except
in some rough sense. This is a very serious matter because questions
of symmetry lie at the heart of modern physics. A related feature is
that the topology and differential geometry of the underlying
manifold are only treated indirectly; the lattice points are
generally manipulated like a trivial topological set. There do exist
radical attempts to overcome these limitations using partially
ordered sets \cite{FRAMP}, but their potentials have yet to be fully
realised.

\subsection{The Fuzzy Sphere}

In the early $1990$'s John Madore began to propose that \ncg be used
in nonpertubative regularisation. (Madore was not the first to
consider such an idea; as early as the $1940$'s Synder made a
suggestion that space-time coordinates should be noncommutative.)
Madore showed that to every \cpt coadjoint orbit (a distinguished
type of manifold) one could associate a canonical sequence of \nc
\algn s, each of which retains the symmetries of the original space.
As $n$ goes to infinity, the \algn s approach the \alg of \cts \fns
of the space (in a certain loose sense). Since each of these \nc
\algs is finite-dimensional, it was hoped that a well-defined
version of \qft could be expressed in terms of them. Moreover, it
was also hoped that the resulting theory would not suffer from
$UV$-divergences.

Madore's construction is not an entirely mathematical exercise. His
approach is also physically motivated. We quote Hawkins
\cite{HAWKINS2}:

`The observation of structures at very small distances requires
radiation of very short wavelength and correspondingly large energy.
Attempting to observe a sufficiently small structure would thus
require such a high concentration of energy that a black hole would
be formed and no observation could be made.  If this is so, then
distances below about the Planck scale are unobservable  and thus
operationally meaningless. If short distances are meaningless, then
perhaps precise locations are as well. This suggests the possibility
of uncertainty relations between position and position, analogous to
the standard ones between position and momentum.  An uncertainty
relation between, say, $x$-position and $y$-position, would mean
that the $x$ and $y$ coordinates do not commute (see \cite{d-f-r}).'

\bigskip

We shall begin by outlining Madore's `fuzzification' of the
$2$-sphere. Let us recall that $S^2$, the {\em $2$-sphere}, is the
submanifold of $\bR^3$ consisting of all points $(x_1,x_2,x_3)$ for
which
\begin{equation}\label{5:eqn:sphereradius}
x_1^2+x_2^2+x_3^2=1.
\end{equation}
It is felt that if a nonpertubative regularized version of \qft
could be constructed on $S^2$, then a regularized version of \qft on
ordinary space-time would follow.

Let us define three functions $\hat{x}_i$,~ $i=1,2,3$, on $S^2$ by
setting $\hat{x}_i(x_1,x_2,x_3)=x_i$, for $(x_1,x_2,x_3) \in S^2$.
We shall denote by $\text{Cor}(S^2)$ the smallest sub\alg of
$C(S^2)$ containing the \fns $\hat{x}_i$. Clearly, $\text{Cor}(S^2)$
is dense in $C(M)$, and each element of $\text{Cor}(S^2)$ is of the
form
\[
z_0+\sum_{i=1}^3z_i\hat{x}_i+\sum_{1\leq i < j \leq
3}z_{ij}\hat{x}_i\hat{x}_j+\sum_{1\leq i < j < k \leq
3}z_{ijk}\hat{x}_i\hat{x}_j\hat{x}_k + \ldots .
\]

Madore's method uses truncations of the above sum to construct a
sequence of \nc \algs approximating $\text{Cor}(S^2)$.

If we discard all functions with nonconstant terms, then
$\text{Cor}(S^2)$ reduces to $\bC$; this is our first approximation
and we denote it by $S^2_1$.

If we add to $S^2_1$ all terms linear in $\hat{x}_i$, then we get a
four dimensional linear space; we shall denote it by $S^2_2$.
Consider the unique isomorphism $\phi_2:S^2_2 \to M_2(\bC)$ for
which $\phi_2(1)=1$ and
\[
\phi_2(\hat{x}_i)=\frac{1}{\sqrt{3}}\s_i;
\]
where $\s_i$ are the Pauli spin matrices. We can use $\phi_2$ to
give $S^2_2$ a \nc \alg structure in an obvious way. Note that
equation (\ref{5:eqn:sphereradius}) holds in this context.

We shall now construct $S^2_3$. To begin with, we include all the
elements of $S^2_2$, and all terms of the form
$z_{ij}\hat{x}_i\hat{x}_j$, for $i<j$. Since we no longer intend to
use the multiplication of $C(S^2)$, we shall also include terms of
the form $z_{ij}\hat{x}_i\hat{x}_j$, for $j<i$, where
${\hat{x}_i\hat{x}_j \neq \hat{x}_j\hat{x}_i}$. It will also be
assumed that equation (\ref{5:eqn:sphereradius}) holds. We can give
$S^2_3$ the structure of a \nc \alg (that is consistent with our
assumptions) in much the same way as we did for $S^2_2$; however, we
shall need to present some facts about the representation theory of
$\frak{su}(2)$ first.

Let $\frak{su}(2)$ be the Lie \alg of $SU(2)$ ($\frak{su}(2)$ is a
three dimensional \alg canonically associated to $SU(2)$, in a sense
that we shall explain below). As is well known, for every \pve
integer $n$, there is a unique $n$-dimensional representation $U_n$
of $\frak{su}(2)$. As is also well known, it is possible to choose a
self-adjoint basis $\{J_i\}_{i=1}^3$ of $U_n(\frak{su}(2))$  \st
$J_i^2=1$, and
\begin{equation} \label{5:eqn:comm}
[J_j,J_k]=2i\epsilon_{jkl}J_l,~~~~~ 1 \leq j,k,l \leq 3;
\end{equation}
where $\epsilon$ is the Levi--Civita permutation symbol. (For $n=2$,
the set of Pauli spin matrices is such a basis.) It can be shown
that every element of $M_n(\bC)$ is expressible as a polynomial in
$J_1,J_2,J_3$.

Let $\phi_3:S^2_3 \to M_3(\bC)$ be the unique linear space
isomorphism for which $\phi_3(1)=1$ and
\[
\phi_3(\hat{x_i})=\frac{1}{\sqrt{8}}J_i.
\]
Obviously, we can use $\phi_3$ to give $S^2_3$ the structure of a
\nc \algn. Note this new structure is consistent with our earlier
assumptions about the linear structure of $S^2_3$.

For $n > 3$, it is now clear how to define the linear space $S^2_n$.
As before, we assume that $\hat{x}_i\hat{x}_j \neq
\hat{x}_j\hat{x}_i$, and that equation (\ref{5:eqn:sphereradius})
holds. We define a \nc multiplication for $S^2_n$ using the unique
isomorphism $\phi_n:S^2_n \to M_n(\bC)$ for which $\phi_n(1)=1$ and
\[
\phi_n(\hat{x_i})=\frac{1}{\sqrt{n^2-1}}J_i.
\]
Note that this new structure is again consistent with the linear
space structure of $S^2_n$. We call $S^2_n$ the $n$-{\em fuzzy
sphere}.

\bigskip

The fuzzy spheres are intuitively thought of as `\nc lattice
approximations to $S^2$'. Unlike the lattice approximations,
however, each $M_n(\bC)$ retains the symmetries of the sphere: The
symmetry group of $S^2$ is clearly $SU(2)$. As is well known, there
exists an irreducible representation of $SU(2)$ on $\bC^n$, for each
$n > 0$. Thus, the primary shortcoming of the lattice approach is
rectified.

The dimension of the \alg $S_n^2$ is thought of as the number of
points in the \nc lattice approximation. (Note that in the
commutative case, the dimension of the function \alg of the lattice
is equal to the number of points in the lattice.) Thus, as $n \to
\infty$, we think of the \nc lattice as `approaching' the continuum.
Furthermore, it is easy to see that equation (\ref{5:eqn:comm})
implies that
\[
[\hat{x}_i,\hat{x}_j]=\frac{2i}{\sqrt{n^2-1}}\epsilon_{ijk}
\hat{x}_k, ~~~~~~ \hat{x_i},\hat{x_j}, \in S_n.
\]
This motivates the intuitive statement: `In the limit $n \to
\infty$, the \alg $S^2_n$ becomes commutative'. For these reasons,
physicists write that `$S_n^2$ converges to $\text{Cor}(S^2)$'.

\bigskip

It is worth our while to comment briefly on the mathematical form
that the fuzzy sphere takes. Its presentation  owes a lot to its
physical origins. In the quantum theory of rotations, $SU(2)$ and
$\frak{su}(2)$ play central roles. As a result, physicists are
usually well versed in their properties.

\subsection{Fuzzy Coadjoint Orbits}

As we have presented it above, Madore's construction of the fuzzy
sphere seems somewhat ad hoc. However, it is in fact a special case
of a much more general procedure known as Berezin--Toeplitz
quantization. This is a method used in quantum mechanics for
`quantising' a phase space $M$ that is a K\"{a}hler manifold (a
K\"{a}hler manifold is a special type of symplectic manifold). By
{\em quantising} we mean producing a sequence of finite-dimensional
Hilbert spaces $H_n$, and a sequence of maps~ $T_n:C(M) \to B(H_n)$.
The maps $T_n$ are called the {\em Toeplitz maps} for the
quantization. The dimension of the Hilbert spaces are governed by
the Riemann-Roch formula; see \cite{HAWKINS2} for details. Madore's
insight was that he could use Berezin--Toeplitz quantisation to
construct fuzzy versions of spaces.

\subsubsection{Coadjoint Orbits}

If the K\"ahler manifold in question is a coadjoint orbit, then
Berezin--Toeplitz quantisation takes a simpler form. Let us recall
some of the basic theory of coadjoint orbits. Let $G$ be a \cpt Lie
group (we shall assume that all the Lie groups we work with are
connected and semi-simple). The map
\[
\l_g:G \to G, ~~~~~~ h \mto gh,
\]
is clearly a diffeomorphism, for all $g \in G$. Moreover, $\l_g$
induces a map \linebreak ${\l^*_g:T_h(G) \to T_{gh}(G)}$ that is
defined by setting
\begin{equation} \label{2:eqn:multtang}
\l^*_g(v_h) = v_h \circ \l_g, ~~~~~~~~~ h \in G.
\end{equation}
If $X$ is a vector field \st $[\l^*_g(X)](h)=X(gh)$, for all $g,h
\in G$, then we say that $X$ is {\em left-invariant}. The
left-invariant vector fields on a Lie group form an \alg that is
closed under the Lie bracket; we call it the {\em Lie \alg} of $G$
and we denote it by $\fg$. As a little thought will verify, $\fg$ is
canonically isomorphic to the tangent plane at any point of $G$. It
is common practice to equate it with the tangent plane at $e$, where
$e$ is the identity of $G$. Thus, $\fg$ is a complex linear space
whose dimension is equal to the manifold dimension of
$G$. 
(It is straightforward to build on the definition of a
left-invariant field and define the notion of a left-invariant
differential form, as referred to in the previous chapter.)

Let us consider the {\em conjugation map}
\[
I_g:G \to G,~~~~~ h \mto ghg \inv, ~~~~~~~~  g \in G.
\]
In direct analogy with equation (\ref{2:eqn:multtang}), $I_g$
induces a map $I^*_g:T_h(G) \to T_{I_g(h)}(G)$, for all $h \in G$.
Since $I_g(e)=e$, $I^*_g$ is a linear mapping on $T_e(G)$. The fact
that we equated $\fg$ and $T_e(G)$ means that we can also regard
$I^*_g$ as a linear map on $\fg$. A routine calculation will now
verify that
\[
\text{Ad}:G \to \End(\fg), ~~~~~ g \to I^*_g,
\]
is a representation of $G$; we call it the {\em adjoint
representation} of $G$. The {\em coadjoint representation} of $G$ is
the mapping
\[
\text{Ad}^*:G \to \End(\fg^*), ~~~~~ g \mto \Ad^*_g;
\]
where $\Ad^*_g$ is the linear dual of  $\Ad_g$, that is, it is the
unique mapping for which
\[
[\Ad^*_g(\mu)](\z)=\mu(\Ad_g(\z)), ~~~~~~~~ \z \in \fg, \mu \in
\fg^*.
\]
We define the {\em coadjoint orbit} of $\mu \in \fg^*$ to be the
subset
\[
\O_\mu=\{\Ad^*_g(\mu):g \in G\}.
\]
Consider the subgroup $H=\{h \in G:\Ad^*(\mu)=\mu\}$. It can be
shown that $\O_\mu$ is homeomorphic to $G/H$, and that $G/H$ is a
homogeneous space for $G$ (recall that a {\em homogeneous space} is
a topological space with a transitive group action).

\subsubsection{Quantized Coadjoint Orbits}

Each coadjoint orbit $\O_\mu \simeq G/H$ comes naturally equipped
with a symplectic form. One can use this form to endow $\O_\mu$ with
the structure of a K\"ahler manifold. When Berezin--Toeplitz
quantization is applied to $\O_\mu$, each Hilbert space $H_n$
produced is the representation space of a unitary irreducible
representation of $G$. (Strictly speaking, we should only be
considering {\em integral} coadjoint orbits, but we shall not need
to concern ourselves with such details here. See \cite{HAWKINS1} for
further details.) Thus, each {\em fuzzy coadjoint orbit} $B(H_n)$
retains the homogeneous space symmetries of $\O_\mu$. As $n$ goes to
infinity, the fuzzy coadjoint orbits $B(H_n)$ converge to
$C(\O_\mu)$ in a loose intuitive sense that is analagous to the
fuzzy sphere case.

\bigskip

The $2$-sphere is a coadjoint orbit of $SU(2)$. When the
Berezin--Toeplitz quantisation is applied to it, one gets the
sequence of $n$-fuzzy spheres described earlier. The Toeplitz maps
$T_n$ produced correspond to the unique set of maps
$F_n:\text{Cor}(S^2) \to S_n^2$, for which $F_n(\hat{x_i}) \to
\hat{x_i}$.

Besides the fuzzy sphere, many other examples of fuzzy coadjoint
orbits have been explicitly described. For example, fuzzy complex
projective spaces \cite{fuzzyprof},  fuzzy complex Grassmannian
spaces \cite{fuzzygrass}, and fuzzy orbifolds \cite{ORBI} have
appeared in the papers of DIAS members.

\subsection{Fuzzy Physics}

Now that we have introduced the notion of a fuzzy space, we are
ready to (briefly) discuss fuzzy physics. In essence, the subject
consists of \nc generalisations of physical theories expressed in
terms of fuzzy space matrix \algn s. Its long term goal is to
construct a fuzzy version of the standard model. Towards this end,
a lot of work has been put into the construction of fuzzy vector
bundles, fuzzy Lagrangians, fuzzy Dirac operators, and fuzzy gauge
theories. While a good deal of progress has been made, it seems
that a fuzzy version of the full standard model is still  some way
off. A large portion of this work has focused on the fuzzy physics
of the fuzzy sphere; in particular the fuzzification of scalar
field theories on $S^2$.

An obvious question to ask is if these new theories give finite
answers where ordinary \qft fails to do so. Unfortunately, the
answer in general is no.  At present, it is not clear whether these
divergences are due to `incorrect' \nc generalisations; or whether
there are fundamental problems with the present form of fuzzy
physics. It should be noted, however, that there are a large number
of examples that {\em are} naturally regularised; a notable example
is \cite{PRES}.

An interesting (albeit troublesome) feature of fuzzy physics is
$UV$--$IR$ mixing. In physics it is usually possible to organize
physical phenomena according to the energy scale or distance scale.
The short-distance, ultraviolet ($UV$) physics does not directly
affect qualitative features of the long-distance, infrared physics
($IR$), and vice versa. However, in fuzzy physics interrelations
between $UV$ and $IR$ physics start to emerge. This occurrence is
known as {\em $UV$--$IR$-mixing}. Unfortunately, this mixing leads
to divergences. In \cite{DIAS} members of the DIAS showed how to
overcome this problem in the case of the fuzzy sphere.

In recent years DIAS members Denjoe O' Connor and Xavier Martin
have been using computer based numerical simulations to study the
fuzzy physics of the fuzzy sphere $S_n^2$. Their work studies the
behaviour of certain fuzzy models as $n \to \infty$. For details
see \cite{DIASNUM1,DIASNUM2}.

Recently, an application of the fuzzy sphere to the study of black
hole entropy has emerged. It has been proposed that the event
horizon of a black hole should be modeled by a fuzzy sphere. A paper
on the subject \cite{blackhole} has been written by DIAS member
Brian Dolan.

\bigskip

Before we finish our discussion of fuzzy physics, we should note
that there exist other ways of applying \ncg to \qftn. As one
notable example, we cite Seiberg and Witten's paper \cite{SEIB}.

\chapter{Compact Quantum Metric Spaces}

As we saw in the previous chapter, a number of physicists are now
working with \ncg in the hope that it will provide a means to solve
the problems of quantum field theory. In practice, their work
involves the construction of `field theories' on \algn s of $n \by
n$ matrices. As $n$ goes to infinity, the matrix \algn s converge
(in some intuitive sense) to the \alg of \cts \fns on a coadjoint
orbit; and the field theories converge to some classical field
theory. The prototypical example of this process is the
fuzzification of the two sphere.

After reading the fuzzy physics literature, Marc Rieffel suspected
that the matrix \alg convergence referred to therein, involved some
kind of continuous $C^*$-field structure. (A $C^*$-field is a type
of bundle construction where each fibre is a \calgn.) However, he
later changed his mind. In \cite{RIEF2} he wrote:

`There is much more in play than just the continuous-field aspect.
Almost always there are various lengths involved, and the writers
are often careful in their bookkeeping with these lengths as $n$
grows. This suggested to me that one is dealing here with metric
spaces in some quantum sense, and with the convergence of quantum
metric spaces'.

Rieffel set himself the formidable task of formulating a
mathematical framework in which he could precisely express these
ideas. Previous to this the only major example of metric
considerations in the literature of quantum mathematics was Connes'
state space metric, as discussed in Chapter $4$. Rieffel used this
as his initial point of reference. However, for reasons that we
shall explain later, Connes' metric proves inadequate for the task
at hand. So Rieffel decided to base his work on that of Kantorovich
instead.

In this chapter we shall begin by presenting \cpt quantum metric
spaces, quantum Gromov--Hausdorff convergence, and some important
examples of both. We shall then go on to establish what is arguably
the most important achievement of the theory to date, a
formalisation of the convergence of the fuzzy spheres to the sphere.
Finally, we shall present some recently proposed modifications to
Rieffel's definitions.

\bigskip

The theory of \cqmsn s is still a very young area. Connes' state
space metric, which can be considered its starting point, was only
introduced in $1989$. Rieffel's first papers on the subject emerged
in the late nineties, and the subject only began to take definite
shape after the year $2000$. In fact, some of the more recent
developments discussed in this chapter occurred in  $2005$ and
$2006$.

With regard to future work, it seems to be Rieffel's intention to
expand his theory so that it will be able to define a distance
between some appropriately defined \nc version of \vbdn s. This
would then enable one to discuss the convergence of  field theories
on fuzzy spaces as well as the convergence of the spaces themselves.
Two recent papers of Rieffel \cite{RIEFVB1,RIEFVB2} indicate that
work in this area is well under way.

We should also mention that Nik Weaver has formulated his own \nc
generalisation of metric spaces based on von Neumann \algs and $w^*$
derivations. Weaver's book \cite{WEAV} gives a very good
presentation of this theory. While the two subjects have a very
different feel to them, there does seem to be some common ground.
For a brief discussion of their relationship see \cite{RIEF4}.

\section{Compact Quantum Metric Spaces}

\subsection{Noncommutative Metrics and the State Space}

Recall that if $X$ is a \cpt Hausdorff space, then by Theorem
\ref{1:thm:cptsurjective} it is homeomorphic to $\Om(C(X))$ endowed
with the weak$^*$ topology. Thus, metrics that metrize the topology
of $X$ are in one-to-one correspondence with metrics that induce the
weak$^*$ topology on $\Om(C(X))$. As a result a natural candidate
for a \nc metric space would be a pair $(\A,\r)$, where $\A$ is a
\calg and $\r$ is a metric on $\Om(\A)$. However, as we showed
earlier, the character space of a \nc \calg is not guaranteed to be
non-empty. Thus, any generalisation based upon it is unlikely to be
very fruitful. A better proposal might be to examine metrics that
induce the weak$^*$ topology on some non-empty set of functionals
that is equal to $\Om(\A)$ in the commutative case and properly
contains  it in the \nc case.

Recall that a state $\mu$ on a \calg is a \pve linear functional of
norm $1$, and that $\Om(\A) \sseq S(\A)$. Unlike $\Om(\A)$, $S(\A)$
is always non-empty; it is well known that for every normal element
$a \in \A$, there exists a state $\mu$ \st $|\mu(a)|=\|a\|$. Using
the Banach--Alaoglu Theorem $S(\A)$ can easily be shown to be \cpt
with respect to the weak$^*$ topology. Furthermore, it is quite easy
to show that $S(\A)$ is a convex subset of $\A^*$, and that when
$\A$ is commutative, $\Om(\A)$ is equal to the set of extreme points
of $S(\A)$; or, as they are better known, the {\em pure-states} of
$\A$. The Krein--Milman Theorem implies that it is also non-empty in
the \nc case. Thus, the set of pure-states of a \calg seems to be
the set of functionals we are looking for. Unfortunately, however,
this set is quite badly-behaved \wrt our needs; and, even though the
generalisation is no longer strict, it turns out to be much more
profitable to use the entire state space of a \calg instead. As we
shall see, the compactness of $S(\A)$  makes it well suited to our
needs.

\subsection{Lipschitz Seminorms}

To restate, if $\A$ is the \calg of \cts \fns on some \cpt Hausdorff
space $X$, then there is a canonical correspondence between metrics
that metrize the topology of $X$ and metrics that induce the
weak$^*$ topology on $\Om(\A)$. We have no such correspondence for
metrics on $S(\A)$. It seems reasonable to assume that a study of
all metrics that induce the weak$^*$ topology on $S(\A)$ would be
too broad. However, it is not clear what subfamily of metrics we
should restrict our attention to.

Recall that for a spectral triple $(A,H,D)$ Connes defined a metric
on $S(\ol{A})$ by
\[
\r_D(\mu, \nu) = \sup\{|\mu(a)-\nu(a)|:\|[D,a]\| \leq 1\}.
\]
In general, it proves quite difficult to say when Connes' metric
induces the weak$^*$ topology on $S(\ol{A})$. This makes $\r_D$
unsuitable for our later work; not only because of Theorem
\ref{1:thm:cptsurjective}, but also because the weak$^*$ topology is
needed for Rieffel's definition of quantum Gromov--Hausdorff
distance to work smoothly.

Connes' work on \nc metrics was in part based upon Kantorovich's
work \cite{KANT} with Lipschitz seminorms. Thus, a different
generalisation of Kantorovich's work might identify a `suitable'
subfamily of state space metrics.

\subsubsection{Kantorovich's Work}

For a compact metric space $(X,\rho)$, the {\em Lipschitz seminorm}
is a commonly used seminorm on $\A = C(X)$ defined by
\[
L_{\rho}(f) = \sup\{|f(x)-f(y)|/\rho(x,y): x \ne y\}.
\]
Finiteness of the seminorm is not guaranteed, but the subset $\{f:
L_\r(f) < \infty \}$ is always dense in $C(X)$. (Note that in this
chapter we shall allow seminorms and metrics to take infinite
values.) Just as we recovered the geodesic distance using $\r_D$,
the well known result
\begin{equation} \label{5:eqn:recovermetric}
\rho(x,y) = \sup \{|f(x)-f(y)|/\rho(x,y): L_\r(f) \leq 1\},
\end{equation}
allows us to recover the metric using $L_\rho$. We now have a
reformulation of the metric space data in terms of the commutative
\calg $C(X)$ and the seminorm $L_\rho$. (If $(\csm,L^2(M,S),\Dirac)$
is the canonical spectral triple on a \cpt Riemannian spin manifold,
then it holds that $\|[\Dirac,f]\|=L_\rho(f)$, where $\rho$ is the
geodesic distance of $M$. Thus, in this case equation
(\ref{5:eqn:recovermetric}) and equation (\ref{eqn2:Diracgeodesic})
coincide. However, we should note that equation
(\ref{5:eqn:recovermetric}) is true for any metric space $X$, and
so, it is a much more general result.)

Anticipating Connes, Kantorovich  defined a metric on the state
space of $C(X)$ by
\begin{equation}\label{5:eqn:kantmetric}
\rho_L(\mu,\nu) = \sup\{|\mu(f)-\nu(f)|: L_{\rho}(f) \le 1\}.
\end{equation}
He showed that, amongst other properties, the topology induced by
$\rho_L$ on $S(C(X))$ always coincides with the weak$^*$ topology.

All of this suggested to Rieffel that a \nc metric space should
consist of a noncommutative unital $C^*$-algebra $\A$ endowed with a
suitable seminorm.  A metric $\rho_L$ could then be defined on the
state space by the obvious generalisation of
equation~(\ref{5:eqn:kantmetric}) to the \nc case.

\subsection{Order-unit Spaces}

One might now assume that we were ready to give a definition of a
\nc metric space. However, our proposed formulation must undergo a
simplification.

Let $\A$ be a \calg and let $L$ be a seminorm on $\A$. If $L$ is to
be considered as a generalised Lipschitz seminorm, then it seems
reasonable to assume that, just as in the classical case,
$L(a^*)=L(a)$, for all $a \in \A$. Let us define a possibly
infinite-valued metric on $S(\A)$ using the obvious generalisation
of equation (\ref{5:eqn:kantmetric}), and let $\mu,\nu \in S(\A)$,
and $\d > 0$ be given. Then there exists an $a \in \A$ \st $L(a)
\leq 1$ and
\[
\r_L(\mu,\nu) - \d \leq  |\mu(a)-\nu(a)|=\mu(a)-\nu(a).
\]
Define an element of $\A_{sa}$ by $b=(a+a^*)/2$ and note that since
$L(a^*)=L(a)$, $L(b) \leq 1$. Now it is well known that all states
are Hermitian (see \cite{MUR} for a proof). Therefore
\[
|\mu(b)-\nu(b)|= \Re(\mu(a)-\nu(a)) =  \mu(a)-\nu(a) \geq
\r_L(\mu,\nu)-\d.
\]
This, means that when we are calculating the value of $\r_L$ it
suffices to take the supremum over the self-adjoint elements of
$\A$. This fact seems to suggest that self-adjoint elements have a
distinguished role to play in any formulation of \nc metric spaces.
In fact, it prompted Rieffel to suggest that order-unit spaces
should play a part.

\begin{defn}
An {\em order-unit space} is a pair $(A,e)$, where $A$ is a real
partially-ordered linear space, and $e \in A$ is a distinguished
element called the {\em order-unit}, such that;

\begin{enumerate}
\item For each $a \in A$, there
is an $r \in {\mathbf R}$ such that $a \le re$ ({\em order-unit
property}).

\item If $a \in A$ and if $a \le
re$, for all $r \in {\mathbf R}$ with $r > 0$, then $a \le 0$ ({\em
Archimedean property}).
\end{enumerate}
\end{defn}

For any \calg $\A$, $\A_{sa}$ is an order-unit space \wrt the
partial order defined by setting $a \geq b$ if $(a-b) \geq 0$, $a,b
\in \A$. We shall take some care to establish this. For any $a \in
\A_{sa}$, note that $C^*(a,1)$, the smallest $C^*$-sub\alg of $\A$
containing $a$ and $1$, is a commutative unital \calgn. By the
Gelfand--Naimark Theorem $C^*(a,1) \cong C(X)$, for some \cpt
Hausdorff space $X$. Using this fact, it can easily be shown that
$\A_{sa}$ satisfies properties $1$ and $2$. Hence, it is indeed an
order-unit space. It is interesting to note that just as any \calg
can be concretely realised as a norm closed $*$-sub\alg of $B(H)$,
for some Hilbert space $H$, any order-unit space can be concretely
realised as a real linear subspace of the self-adjoint operators on
a Hilbert space.

The {\em standard norm} on an order-unit space is defined by
\[
\|a\| = \inf\{r \in {\mathbf R}: -re \le a \le re\}.
\]
By returning to the fact that $C^*(a,1) \cong C(X)$, for any $a \in
\A_{sa}$, we can establish that the order-unit norm on $\A_{sa}$
coincides with the restriction of the $C^*$-norm to $\A_{sa}$.

A {\em state} on an order-unit space $(A,e)$ is a bounded linear
functional $\mu$ on $A$ such that $\mu(e) = 1 = \|\mu\|$ (where
$\|\mu\|=\sup\{\mu(a):\|a\| \leq 1,~ a \in A\}$). It can be shown
that all states are automatically positive. We denote the set of all
states on $A$ by $S(A)$ and call it the {\em state-space} of $A$.
Just as for \calgsn, we can show that $S(A)$ is a weak$*$ closed
subset of the unit ball of $A^*$. Thus, by the Banach--Alaoglu
theorem, it is compact \wrt the weak$^*$ topology on $A^*$. Again,
just as in the \calg case, $S(A)$ is a convex subset of $A^*$.

It is well known \cite{MUR} that that a bounded linear functional
$\mu$ on a unital \calg is \pve \iff $\mu(1)=\|\mu\|$. Thus, if $\A$
is unital and $\mu \in \A^*$, then $\mu$ is a state \iff
$1=\mu(1)=\|\mu\|$. As a little thought will verify, this implies
that the restriction of a state on $\A$, to $\A_{sa}$, will be a
state in the order-unit sense, and conversely, that every state on
$\A_{sa}$ has a unique extension to a \calg state on $\A$. Thus,
there is a one-to-one correspondence between the elements of $S(\A)$
and $S(\A_{sa})$. If we endow both spaces with their respective
weak$^*$ topologies, then it is a simple exercise to show that they
are homeomorphic. As we saw above any seminorm $L$ on $\A$
satisfying $L(a^*)=L(a)$ will induce the same metric on $S(\A)$ as
its restriction will induce on $S(\A_{sa})$. Thus, two \calgs whose
collection of  self-adjoint elements are isomorphic as order-unit
spaces will, loosely speaking, produce the same `metric data'. When
we also take into account the greater technical flexibility afforded
by working with order-unit spaces, it seems that they are the
natural structure upon which to base a \nc metric theory.

\subsection{Compact Quantum Metric Spaces}

Now that we have settled on the category that we shall be working
in, we are ready to formulate a generalisation of the Lipschitz
seminorm.
\begin{defn}
Let $A$ be an order-unit space and let $L$ be a seminorm on $A$
taking finite values on a dense order-unit subspace of $A$. We say
that $L$ is a {\em Lip-norm} if the following conditions hold:
\begin{enumerate}
\item $L(a) = 0$ ~~~   \iff   $a \in \bR e$.\\
\item The topology on $S(A)$ induced by the metric
\[
\rho_L(\mu,\nu) = \sup\{|\mu(a)-\nu(a)|: L(a) \leq 1\},~~~~\mu,\nu
\in S(\A),
\]
coincides with the weak$^*$ topology.
\end{enumerate}
\end{defn}

The first requirement on $L$ is a direct generalisation of the fact
that $L_{\rho}(f)=0$, \iff $f$ is a constant \fnn, while the~ second
requirement has already been motivated. As might be expected,
Lip-norms do not directly generalise Lipschitz seminorms; that is,
if $X$ is a \cpt Hausdorff space, then there can exist Lip-norms on
$C(X)$ that are not the Lipschitz seminorm for any metric on $X$.

We now are ready to define our \nc generalisation of \cpt metric
spaces.

\begin{defn}
A {\em compact quantum metric space} is a pair $(A,L)$ where $A$ is
an order-unit space and $L$ is a Lip-norm on $A$.
\end{defn}

Rieffel uses the term quantum metric space instead of \nc metric
space because there is no multiplication in an order-unit space, and
so, there is no noncommutativity to speak of. He also cites the
central role that states play in quantum mechanics as a motivation.

The metric $\r_L$ only takes finite values. To see this assume that
$\r_L(\mu_0,\nu_0)=+\infty$, for some $\mu_0,\nu_0 \in S(A)$. The
proper subset $\{\mu:\r_L(\mu,\nu_0) < \infty\}$ is both open and
closed, which is impossible since the convexity of $S(A)$ implies
that it is connected \wrt the weak$^*$ topology. Since $S(A)$ is
compact, $\r_L$ is bounded, and so we can speak of the radius of the
metric space $S(A)$. We define the {\em radius} of a \cqms to be the
radius of its state space.

\bigbreak

In naturally occurring examples it can often be quite difficult to
verify directly whether or not a metric induces the weak$^*$
topology. Fortunately, however, Rieffel \cite{RIEF3} managed to
reformulate this property in simpler terms.

\begin{thm}\label{thm:5:finiteradius}
Let $L$ be a seminorm on an order-unit space $A$, such that ${L(a) =
0}$ \iff $a \in \bR e$, and let $B_1 = \{a \in A: L(a) \le 1\}$.
Then $\r_L$ induces the weak$^*$ topology on $S(A)$ exactly if
\begin{enumerate}
\item  $(A, L)$ has finite radius, and
\item  $B_1$ is totally bounded in $A$.
\end{enumerate}
\end{thm}

An analogous reformulation is not known to exist for Connes' state
space metric.

\bigskip

It is interesting to note that recently Frederic
Latr{\'e}moli{\'e}re \cite{LAT2}, a former doctoral student of
Rieffel, has shown how to extend the definition of a \cqms to one
that generalises locally compact metric spaces.

\subsection{Examples}

Let $(X,\r)$ be a \cpt metric space and let $L_\r$ be the Lipschitz
seminorm on $C(X)$. Let us also use $L_\rho$ to denote the
restriction of $L_\r$ to $C(X)_{sa}$, the space of real-valued \fns
on $X$. The pair $(C(X)_{sa},L_\r)$ is clearly a \cqmsn; we call it
the {\em classical \cqms associated} to $(X,\r)$. It is not too hard
to show that given $(C(X)_{sa},L_\r)$ one can reproduce $(X,\r)$,
and so no information is lost by focusing on the order-unit spaces.

\bigskip

What we would now like to see are some purely quantum examples. Most
of these are constructed using the actions of groups on \calgsn. Let
$\A$ be a \linebreak \calg and let $G$ be a compact group. An {\em
action of $G$ on $\A$ by automorphisms} is a homomorphism $\a: G \to
\Aut (A),~g \mto \a_g$, where $\Aut (A)$ is the space of
automorphisms of $\A$. An action is called {\em strongly \cts} if,
for each $a \in \A$, the mapping
\[
G \to \A,~~~g \mto \a_g(a)
\]
is \ctsn. A {\em length-function} $\ell$ on $G$ is a \fn that takes
values in $\bR^+$ \st $\ell(xy) \leq \ell(x) + \ell(y)$,
$\ell(x^{-1}) = \ell(x)$, and  $\ell(x) = 0$ \iff $x = e$. (We note
that every length \fn $\ell$ on a group gives a metric $\r$ that is
defined by $\r(x,y) = \ell(xy\inv)$. Moreover, since
\[
\r(xz,yz)=\ell(xz(yz)\inv)=\ell(xzz\inv y \inv)=\r(x,y),
\]
the metric is right invariant. Conversely, given a right-invariant
metric on $G$ we can define a length \fn by $\ell(x)=\r(x,e)$.)
Define a seminorm $L$ on $\A$ by
\begin{equation}\label{6:defn:groupseminorm}
L(a) = \sup\{\|\alpha_x(a)-a\|/\ell(x): x \ne e\}.
\end{equation}
(Note that $L(a^*)=L(a)$.) The restriction of $L$ to $\A_{sa}$ seems
like a reasonable candidate for a Lip-norm. If we are to have that
$L(a) = 0$ only when $a \in \bC 1$, then it is clear that $\alpha$
must be ergodic on $\A$; an action $\a$ of a group on a unital \alg
$A$ is called {\em ergodic} if $\a(a)=a$ only when  $a \in \bC1$. In
fact, in \cite{RIEF3} Rieffel showed that ergodicity is all we need.
\begin{thm} \label{thm:5:lengthcqms}
If $G$ is a compact group endowed with a \cts length \fn and
$\alpha$ is a strongly \cts ergodic action of $G$ on $\A$, then $L$,
as defined in equation (\ref{6:defn:groupseminorm}), restricted to
$\A_{sa}$, is a Lip-norm.
\end{thm}
This result is established by verifying the criteria of Theorem
\ref{thm:5:finiteradius}. (It easily fails if $G$ is not compact; as
yet no effort has been made to construct a noncompact version.) This
result provides us with a large stock of good examples and motivates
us to make the following definition: Let $\A$ be a unital \calg and
let $L$ be a seminorm on $\A$ satisfying $L(a^*)=L(a)$, for all $a
\in \A$. If the restriction of $L$ to $\A_{sa}$ is a Lip-norm, then
we call the pair $(\A,L)$ a {\em Lip-normed \calgn} and we call
$(\A_{sa},L)$ its associated \cqmsn. Note that if $X$ is a compact
metric space, then $(C(X),L_\r)$ is a Lip-normed \calgn.

Our search for \cqmsn s now turns into a search for ergodic actions
of compact groups on unital \calgsn. If $G$ is a \cpt group endowed
with a \cts length \fn $\ell$, and if $U$ is an irreducible unitary
representation of $G$ on a Hilbert space $H$, then we can define a
strongly \cts group action $\a:G \to B(H)$ by setting
$\a_g(B)=U_gBU^*_g$. We shall examine this example in greater detail
in Section $6.3$. The corresponding \cqms will be of great
importance to us.

\subsubsection{Quantum Tori}

One of the most important families of examples of spaces in \ncg is
the family of quantum tori. For $\hbar \in \bR$, the quantum torus
$C_\hbar(\bT^2)$ is a $C^*$-sub\alg constructed as follows: let $H$
be the Hilbert space $L^2(\bT^2)$, where $\bT^2 = \bR^2/2 \pi
\bZ^2$; and let $U$ and $V$ be the two bounded linear operators
defined on $H$ by setting
\begin{eqnarray*}
U\,f(x_1,x_2) &=& e^{ix_1}f(x_1,x_2-\frac{1}{2}\hbar),\\
V\,f(x_1,x_2) &=& e^{ix_2}f(x_1+\frac{1}{2}\hbar,x_2),
\end{eqnarray*}
for $f \in H, x_i \in \bT^2$. These are unitary operators and they
obey the commutation relation $UV=e^{i\hbar}VU$. We define the {\em
quantum torus}, for $\hbar$, to be the closed span in
$B(L^2(\bT^2))$ of the operators $U^mV^n$, for $m,n \in \bZ$.

It can be shown that when $\hbar=0$, $C_\hbar(\bT^2) \simeq C(T^2)$.
This is the motivation for the name quantum torus.

Quantum tori have a number of important applications, the most
notable being in Connes' study of the quantum Hall effect
\cite{CON}. They are also a central example in cyclic cohomology
theory.

In \cite{RIEF1} Rieffel defined a canonical strongly \cts ergodic
action of $\bT^2$ on $C_\hbar(\bT^2)$ (where $\bT^2$ is considered
as a compact group in the obvious way). Thus, by choosing a \cts
length \fn on $\bT^2$, which it is always possible to do, one can
give $C_\hbar(\bT^2)_{sa}$ the structure of a \cqmsn.

\bigskip

Other examples of \cqmsn s have been produced from Connes and
Landi's $\ta$-deformed spheres \cite{CONTHETA}. Lip-normed AF-\algs
have been produced using Bratteli's non-commutative spheres
\cite{BRAT}. Quite recently Li \cite{LIQG} has used Podl\'{e}s
definition of an action of a compact quantum group on a \calg as a
means to generate \cqmsn s. He has used this structure to good
effect in studying the types of convergence in quantum field theory
that motivated Rieffel.

\subsection{Spectral Triples}

At this stage it might be interesting for us to reflect on what
connection, if any, exists between \cqmsn s and  spectral triples.
We recall that for a \cpt Riemannian spin manifold $M$, the
Lipschitz seminorm $L_\r$ and the Dirac operator $\Dirac$ are
related by the equation $L_\r(f)=\|[\Dirac,f]\|$, for all ${f \in
\csm}$. Rieffel has made some progress towards establishing a
similar relation in the \nc case. For an arbitrary \cqms $(A,L)$ he
constructed a faithful representation of $A$ on a Hilbert space $H$
that preserves the order-unit structure, and a self-adjoint operator
$D$ on $H$ \st $L(a) = \|[D,a]\|$, for all $a \in A$. A major
shortcoming of his construction is that in general $D$ does not have
\cpt resolvent. If $(\A,L)$ is a Lip-normed \mbox{\calgn}, then the
representation of $\A_{sa}$ can be extended to a linear
representation of $\A$ on $H$ \st $L(a)=\|[D,a]\|$, for all $a \in
\A$. However, the representation is not always a $*$-\alg
homomorphism. While there exist examples for which these problems
does not arise, it is an open question as to what additional \cdns a
\cqms would have to satisfy in order to ensure that they did not
arise in general.

\section{Quantum Gromov--Hausdorff distance}

As explained in the introduction, Rieffel introduced \cqmsn s in the
hope that they could be used to formalise statements about matrix
\algs converging to the sphere. Now, that we have presented \cqmsn s
we shall move onto defining what it means for them to converge.

\subsection{Gromov--Hausdorff Distance}

When speaking of \cvg of ordinary metric spaces the most frequently
used formulism is that of Gromov--Hausdorff distance. It is a
generalisation of Hausdorff distance and it is most commonly used in
the study of \cpt Riemannian manifolds. When the manifold is a spin
manifold, the associated Dirac operator plays a prominent role
\cite{LOTT}. This hinted to Rieffel that he might be able to discuss
convergence of \cqmsn s in terms of a suitable `quantum version' of
Gromov--Hausdorff distance.

We shall begin by recalling the definition of Hausdorff distance.
Let $(X,\rho)$ be a compact metric space and let $Y$ be a subset of
$X$. For any positive real number $r$, define  $\N_r(Y)$, the {\em
open $r$-neighborhood} of $Y$, by
\[
\N_r(Y) = \{x \in X: \rho(x,y) < r, \text{ for some } y \in Y\}.
\]
We define $\dist^\r_H(Y,Z)$, the {\em Hausdorff distance} between
two closed subsets $Y$ and $Z$ of $X$, by setting
\[
\dist^\r_H (X,Y) = \inf\{r: Y \subseteq \N_r(Z) \text{ and } Z
\subseteq \N_r(Y)\}.
\]
The Hausdorff distance defines a metric on the family of all closed
subsets of $X$. The resulting metric space can be shown to be
compact, and it is complete if $X$ is complete.

Gromov generalised this metric to one that defines a distance
between any two compact metric spaces. Let $(X,\rho_X)$ and
$(Y,\rho_Y)$ be two compact metric spaces and let $X {\dot \cup} Y$
denote the disjoint union of $X$ and $Y$. Let $\M(\rho_X,\rho_Y)$
denote the set of all metrics on  $X {\dot \cup} Y$ that induce its
topology, and whose restrictions to $X$ and $Y$ are $\rho_X$ and
$\rho_Y$ respectively. We call the elements of $\M(\rho_X,\rho_Y)$
{\em admissable metrics}. An element $\r \in \M(\rho_X,\rho_Y)$ can
be produced as follows: if $x,x' \in X$, define
$\rho(x,x')=\rho_X(x,x')$; if $y,y' \in Y$, then define
$\rho(y,y')=\rho_Y(y,y')$; if $x \in X$, $y \in Y$, then, for some
fixed $x_0 \in X$, some fixed $y_0 \in Y$, and some fixed positive
number $L$, define $\rho(x,y)=\rho_X(x,x_0)+L +\rho(y,y_0)$. It is
routine to show that this defines a metric that induces the topology
of $X \dot \cup Y$.

Now, for each  metric $\rho$ in $\M(\rho_X,\rho_Y)$, it is clear
that $X {\dot \cup} Y$ is compact, and that $X$ and $Y$ are closed
subsets of $X \dot \cup Y$. Thus, the Hausdorff distance between
them is well defined. Their {\em Gromov--Hausdorff distance}
$\dist_{GH}(X,Y)$ is defined by setting
\[
\dist_{GH}(X,Y) = \inf\{\dist^\r_H(X,Y): \rho \in
\M(\rho_X,\rho_Y)\}.
\]
Gromov showed that if $\dist_{GH}(X,Y)=0$, then $X$ and $Y$ are
isometric as metric spaces. He went on to establish that if $\C\M$
denotes the family of all isometry classes of \cpt metric spaces,
then the pair $(\C\M,\dist_{GH})$ is a complete metric space. He
also established necessary and sufficient conditions for a subset of
$\C\M$ to be totally bounded.

It is interesting to consider the relationship between Hausdorff and
Gromov--Hausdorff distance. As a little thought will verify, if
$(X,\r)$ is a compact Hausdorff space, and $Y$ and $Z$ are closed
subsets of $X$, then
\[
\dist_{GH}(Y,Z) \leq \dist^\r_H(Y,Z).
\]

\subsection{Quantum Gromov--Hausdorff}

We shall now construct a generalised version of Gromov--Hausdorff
distance that will define a distance between two compact quantum
metric spaces $(A,L_A)$ and $(B,L_B)$. An obvious, if somewhat
crude, way to do this would be to take the Gromov--Hausdorff
distance between $S(A)$ and $S(B)$. However, a distance that
involved the Lip-norms of $(A,L_A)$ and $(B,L_B)$ more directly
would be more natural. Let us look to the classical case for some
intuition on how to do this. The space of \cts real-valued \fns on
$X \dot{\cup} Y$ can be identified with the order-unit space
$C(X;\bR) \oplus C(Y;\bR)$. Thus, for any metric $\r$ on $X
\dot{\cup} Y$, we have a corresponding Lipschitz seminorm $L_\r$ on
$C(X;\bR) \oplus C(Y;\bR)$. This prompts us to generalise metrics on
the disjoint union of two \cpt spaces by Lip-norms on direct sum of
two order-unit spaces. (The {\em direct sum} of two order-unit
spaces $A$ and $B$ is defined in the obvious way: take $A \oplus B$,
the direct sum of $A$ and $B$ as linear spaces and define
$(e_A,e_B)$ to be the order-unit, then define a partial order by
setting $(a,b) \leq (c,d)$ if $a \leq c$ and $b \leq d$. The
standard norm on the direct sum is easily seen to satisfy
$\|(a,b)\|=\max\{\|a\|,\|b\|\}$.)

We now need to generalise to the quantum case the notion of an
admissable metric. Let $(X,\rho)$ be an arbitrary compact metric
space, let $Y$ be a closed subset of $X$, and let $f$ be an element
of $C(X)$. Denote the restriction of $\rho$ to $Y$ by $\rho_Y$ and
denote the restriction of $f$ to $Y$ by $\pi(f)$. If $g \in
C(Y;\bR)$ and $f \in C(X;\bR)$ \st $\pi(f)=g$, then it is clear that
\[
L_{\r_Y}(g) \le L_\r(f).
\]
Consider the \fn
\begin{equation*} \label{5:eqn:minfn}
h(x)=\inf_{y \in Y}(g(y)+L_{\r_Y}(g)\r(x,y)) \in C(X;\bR).
\end{equation*}
A short computation will verify that $\pi(h)=g$ and that
$L_{\r}(h)=L_{\r_Y}(g)$. Thus,
\[ \label{eqn:5:lipquotient}
L_{\r_Y}(g)= \inf\{L_\r(f):\pi(f)=g\}.
\]
(In fact, it is not always possible to find an analogue of $h$ for
complex-valued functions. This provides another important reason for
our emphasis on real-valued functions.)

This motivates us to make the following definition: let $(A,L_A)$
and $(B,L_B)$ be two \cqmsn s. We call a  Lip-norm $L$ on $A \oplus
B$ {\em admissible} if
\[
L^q_A(a)=\inf\{L(a,b):b \in B\}
\]
for all $a \in A$, and
\[
L^q_B(b)=\inf\{L(a,b):a \in A\}
\]
for all $b \in B$. We denote the set of all admissible Lip-norms by
$\M(L_A,L_B)$.

In a short while we shall use admissable Lip-norms to define a
distance between the state spaces of $A$ and $B$. Firstly, however,
for sake of clarity, we shall spell out some details about the
relationship between $S(A)$, $S(B)$, and $S(A \oplus B)$. Denote the
canonical injection of $S(A)$ into $S(A \oplus B)$ by $i$; that is,
if $\f \in S(A)$, then $i(\f)(a \oplus b)=\f(a)$. It is easily seen
that $i$ is injective, and that it is \cts \wrt to the weak$^*$
topology. Thus, since $S(A)$ and $S(A \oplus B)$ are both \cpt
Hausdorff spaces, $S(A)$ is homeomorphic to $i(S(A))$. In this sense
we shall consider $S(A)$ to be a closed subset of $S(A \oplus B)$.
Similarly, we shall consider $S(B)$ to be a closed subset of $S(A
\oplus B)$.

It can be shown \cite{RIEF1} that if $L$ is an admissible Lip norm
on $A\oplus B$, then the restriction of $\r_L$ to $S(A)$ is equal to
$\r_{L_A}$, and the restriction of $\r_L$ to $S(B)$ is equal to
$\r_{L_B}$. (The proof relies upon the fact that the metric on each
space induces the weak$^*$ topology.) This pleasing result allows us
to define a quantum version of Gromov--Hausdorff distance.

\begin{defn}
Let $(A,L_A)$ and $(B,L_B)$ be two compact quantum metric spaces.
The {\em quantum Gromov--Hausdorff distance} between them is
\[
\dist_q(A,B) = \inf\{\dist_H^{\rho_L}(S(A),S(B)): L \in {\mathcal
M}(L_A,L_B)\}.
\]
\end{defn}

Quantum Gromov--Hausdorff distance is clearly symmetric, that is,
for two \cqmsn s $(A,L_A)$ and $(B,L_B)$, $d_q(A,B)=d_q(B,A)$. In
\cite{RIEF1} Rieffel showed that if $(A,L_A)$, $(B,L_B)$, and
$(C,L_C)$ are compact quantum metric spaces, then
\[
\dist_q(A,C) \le \dist_q(A,B) + \dist_q(B,C).
\]
He also showed that if $\dist_q(A,B)=0$, then, \wrt an appropriately
defined notion of isometry based on Lip-norms, $(A,L_A)$ is
isometric to $(B,L_B)$. Thus, if we denote the family of isometry
classes of \cqmsn s by $\C \Q \M$, then the pair $(\C \Q \M ,
\dist_q)$ is a metric space. In fact, Rieffel went on to show that
it is a complete metric space, and that analogues of Gromov's
results on the total boundedness of subsets of $\C\M$ also hold.

\bigskip

When the definition of quantum Gromov--Hausdorff distance is applied
to compact metric spaces it does not in general agree with
Gromov--Hausdorff distance. The basic reason why the two definitions
fail to agree is not too difficult to understand. For ordinary
Gromov--Hausdorff distance one is looking, loosely speaking, at the
distance between the pure states of $C(X)$ and the pure states of
$C(Y)$. In the case of \qghd one is looking at the distance between
the states of $C(X)$ and the states of $C(Y)$. It turns out that the
\qghd between two \cpt metric spaces is always less than the
Gromov--Hausdorff distance. Loosely speaking, this is because it is
`more difficult' to find a pure state that is close to a pure state
than it is to find a state that is close to a state.

The set of pure states is badly behaved in the \nc case and it is
not clear how one would develop a useful theory that would define a
distance between the pure states of two \calgn s. In fact, Rieffel
is unsure as to whether the non-equivalence of the two definitions
should be viewed as a defect or as a `quantum feature'. For a more
detailed discussion of the relationship between Gromov--Hausdorff
distance and its quantum version see \cite{RIEF1}.

\bigskip

If $(\A,L_A)$ and $(\B,L_B)$ are two Lip-normed $C^*$-algebras, then
it is unfortunate, but true, that the quantum Gromov--Hausdorff
distance between $(\A_{sa},L_\A)$ and $(\B_{sa},L_\B)$ can be zero
even when $\A$ and $\B$ are not isomorphic as \calgn s. (As we noted
earlier, this cannot happen in the commutative case.) Both David
Kerr and Hanfeng Li have worked towards addressing this shortcoming
of the theory and we shall review their work in the last section of
this chapter.

\subsection{Examples}
In general it proves quite difficult to find the \qghd between two
\cqmsn s. Usually the best one can do is to establish an upper bound
for it (lower bounds are also quite difficult to find). An analogous
situation holds for classical Gromov--Hausdorff distance.

The first major example of \qgh \cvg that Rieffel established
involved a sequence of quantum tori \cite{RIEF1}}. (The quantum tori
being considered as \cqmsn s in the sense explained earlier.) It was
shown that if $\{\hbar(n)\}_{n}$ is a sequence of real numbers
converging to a real number $\hbar$, then the corresponding quantum
tori $C_{\hbar(n)}(T^2)$ converge to $C_{\hbar}(T^2)$ \wrt \qghdn.
In other words, the mapping from $\bR$ to $\C \Q \M$ given by $\hbar
\mto C_{\hbar}(T^2)$ is \cts \wrt the canonical topology of $\bR$
and the topology induced on $\C \Q \M$ by quantum Gromov�-Hausdorff
distance.

Another interesting example, that builds on Rieffel's work, comes
from Latr{\'e}moli{\'e}re \cite{LAT}. He has recently shown that
\wrt quantum Gromov--Hausdorff distance any quantum torus can be
approximated by a sequence of finite-dimensional $C^*$-algebras. He
loosely terms these finite-dimensional \calgs {\em fuzzy tori}. His
motivation for establishing such a result came again from various
statements in quantum field theory. According to Rieffel there is a
wealth of other examples in the physics literature that could be
given formal description using the langauge of \cqmsn s.

In the next section we shall present what is arguably the most
famous example of Gromov--Hausdorff convergence. It involves a
sequence of canonically constructed \cqmsn s converging to a
classical \cqms associated to the sphere. It is interesting because
it gives rigorous expression to our earlier discussion of fuzzy
spheres converging to $S^2$ and it demonstrates very well the
interplay between theoretical physics and mathematics that is so
prevalent in \ncgn.

\section{Matrix Algebras Converging to the Sphere}

As we saw in the previous chapter, the two sphere is a coadjoint
orbit of $SU(2)$. For sake of convenience and generality, most of
the discussion in this section will be in terms of a general
coadjoint orbit $\O_{\mu}$. (Strictly speaking we should only be
considering {\em integral} coadjoint orbits. However, just as in the
previous chapter, we are going to be a little careless about this.)
It is only as we near the end of the exposition that we shall return
to the special case of the two sphere.

\bigskip

To show that $\O_{\mu}$ is the limit of a sequence of matrix \algn
s, we must first find an `appropriate' way to give it the structure
of a \cqmsn. Recall that all coadjoint orbits of a Lie group $G$ are
of the form $G/H$, for some subgroup $H$. We shall use this fact to
endow $\O_{\mu}$ with a \cqms structure. Let $\ell$ be a \cts length
\fn on $G$. As we noted before, $\ell$ induces a metric on $G$ that
is defined by $\r(g,h)=\ell(gh \inv)$. If we assume that
\begin{equation}\label{6:eqn:transinvmetric}
\ell(xgx\inv)=\ell(x),
\end{equation}
then
\begin{equation} \r(xg,xh)=\r(g,h), ~~~~~\text{ for all } x,g,h \in G,
\end{equation}
then $\r$ in turn  induces a metric $\r_{\pi}$ on $G/H$ that is
defined by setting
\[
\r_{\pi}([x],[y])=\inf\{\r(x,y):x \in [x],y \in [y]\}.
\]
We shall use $L_\A$ to denote the Lipschitz seminorm that $\r_{\pi}$
induces on \linebreak ${\A_{sa}\sseq \A =C(G/H)}$.

We note that it is always possible to define a \cts length \fn on a
\cpt Lie group $G$ that satisfies equation
(\ref{6:eqn:transinvmetric}). In fact, `most' canonical metrics on
coadjoint orbits arise in this way. For example, the usual round
metric on the $2$-sphere is of this for

\bigskip

In its standard form $L_\A$ is somewhat awkward to work with.
Fortunately, however, there exists a more convenient formulation.
Let $\l$ be an action of $G$ on $\A$ defined by setting
\begin{equation}\label{defn:5:lamdaaction}
(\l_hf)[g]=f([h \inv g]), ~~~~~ f \in \A,~g,h \in G.
\end{equation}
A series of straightforward calculations will show that
\[
L_{\A}(f)=\sup_{g \neq e}\{\|\l_g(f)-f\|_{\infty}/\ell(g)\}.
\]
(This means that $L_{\A}$ is the Lip norm on $\A$ arising from the
ergodic action $\l$.)

\subsubsection{Compact Quantum Metric Spaces from Group
Representations}

We now need a suitable way to endow the fuzzy space matrix \algs
with a \cqms structure. Earlier in the chapter, we saw that one
could endow the \alg of operators on a Hilbert space with a \cqms
structure using the action of a \cpt group. Since every fuzzy
coadjoint orbit comes naturally endowed with a \cpt group action,
this seems like a very suitable formulation. Let us present what is
involved in detail: let $\a$ be a strongly \cts ergodic action of a
compact group $G$ on a unital \calg $\A$, and let $\ell$ be a \cts
length \fn defined on $G$. Then Theorem \ref{thm:5:lengthcqms}
states that the seminorm $L_\a$, defined by setting
\begin{equation}\label{6:eqn:groupactionlipnoprm}
L_{\a}(a)=\sup\{\|\a_g(a)-a\|/\ell(g):g\neq e_G\},
\end{equation}
is a Lip-norm on $\A_{sa}$. Now, let $U$ be an irreducible unitary
representation of $G$ on a Hilbert space $H$. (We note that every
irreducible representation of a \cpt group is finite-dimensional, as
is well known.) We can define an action $\a$ of $G$ on $\B=B(H)$ by
setting
\[
\a_g(T) = U_gTU_g^*, ~~~~~ g \in G,~T \in \B.
\]
Let us show that this action is strongly \cts and ergodic: if
$\a_g(T)=T$, for all $g \in G$, then $U_gT=TU_g$, for all $g \in G$.
Thus, if $\l$ is an eigenvalue of $T$, and $v$ is an element of the
corresponding eigenspace $E_{\l}$, then ${TU_gv=U_gTv=\l U_gv}$.
This means that $U_gv \in E_{\l}$, and so $E_{\l}$ is invariant
under $U$. To avoid a contradiction we conclude that $\a$ is
ergodic. To see that $\a$ is strongly \cts take a net ${g_\l}$ in
$G$ that converges to $g$ and note that since $U_{{g_\l}} \to U_g$,
and
\[
\|U_g T U_g^* - U_{g_\l} T U_{g_\l}^*\| \leq \|U_gT\|\|U_g^* -
U_{g_\l}^*\|+\|U_g-U_{g_\l}\|\|T U_{g_\l}^*\|,
\]
we must have that $\a_{g_\l}(T) \to \a_{g}(T)$. Hence, the seminorm
defined by equation (\ref{6:eqn:groupactionlipnoprm}) is a Lip-norm,
and $(\B_{sa},L_{\a})$ is a compact quantum metric space.

\subsection{The Berezin Covariant Transform}

Let $\O_\mu \simeq G/H$ and $(\A_\ssa,L_{\A})$ be as above, and let
$(\B_{sa},L_\B)$ be the \cqms associated the to $H_n$, the
$n$-fuzzification of $\O_\mu$, for some $n>0$. Most of the rest of
this section will be spent trying to find an upper bound for the
\qghd between $(\A_\ssa,L_{\A})$ and $(\B_{sa},L_\B)$. If we
calculated $\dist_H^{\rho_L}(S(\A_\ssa),S(\B_\ssa))$, for some
admissible Lip-norm $L$ on $\A_\ssa \oplus \B_\ssa$, then this would
give us such an upper bound.

\bigskip

Before we try to do this, however, we need to introduce an important
reformulation of $\O_\mu$. Let $U_n$ be the representation of
$\O_\mu$ on $H_n$, and let $\xi$ be a {\em highest-weight vector} of
$U_n$ (see \cite{SIM} for details on highest weight vectors). We
define $P$, the projection operator {\em corresponding} to $\xi$, by
setting ${Px = \la x,\xi \ra \xi}$, for $x \in H_n$. We define $R$,
the {\em stabilizer} of $P$, by setting
\[
R=\{g \in G~|~\a_g(P)=P\}.
\]
It can be shown \cite{RIEF2} that $R$ is equal to $H$, and so
$\O_\mu \simeq G/R$.

\bigskip

Inspired by previous work on Gromov--Hausdorff distance in
\cite{RIEF1} Rieffel made the following guess at a Lip norm $L$ on
$\A_\ssa \oplus \B_\ssa$:
\[
L(f,T)=L_{\A}(f) \vee L_{\B}(T) \vee N(f,T),~~~~~\g \in \bC;
\]
where $a \vee b$ denotes the maximum of $a$ and $b$, and $N$ is a
seminorm on $\A_\ssa \oplus \B_\ssa$ that satisfies $N(1_A,1_B)=0$,
among a number of other natural conditions; see \cite{RIEF1} for
details. By verifying the criteria of Theorem
(\ref{thm:5:finiteradius}), Rieffel \cite{RIEF1} showed that $L$
induces the weak$^*$ topology on $S(\A_{sa} \oplus \B_{sa})$. Thus,
since it is clear that $L(1_\A,1_\B)=0$, $L$ must be a Lip norm on
$\A_{sa} \oplus \B_{sa}$. For $L$ to be of use to us, however, it
must be admissable; that is, its quotient seminorms on $\A_{sa}$ and
$\B_{sa}$ must be $L_{\A}$ and $L_{\B}$ respectively. We shall begin
by establishing that the quotient of $L$ on $\B_{sa}$ is equal to
$L_\B$. This will require us to construct a specific form for $N$
using the Berezin covariant transform.

\subsubsection{The Berezin Covariant Transform}

For $T \in \B$, and $\t$ the trace, the {\em Berezin covariant
symbol} of $T$, \wrt $P$, is the \cts mapping
\[
\s_T:G \to \bC,~~~g \to \t(T\a_g(P)).
\]
The mapping
\[
\s:\B \to C(G),~~~~~T \mto \s_T,
\]
is called the {\em Berezin covariant transform}. We can easily see
that $\s_T(gh)=\s_T(g)$, for all $h \in R$. Thus, the function
\[
\s_T : G/R \to \bC,~~~[g] \mto \t(T \a_g(P))
\]
is a well-defined element of $\A=C(G/R)$; and $\s:\B \to \A,~~T \mto
\s_T$ is a well defined mapping. It has a number of useful
properties. Firstly, $\s_1=1$, as can be seen from
\[
\s_1([g])=\t(1\a_g(P))=\t(U_gPU_g^*)=\t(P)=1.
\]
If $T$ is a positive element of $\B$, then, since $\a_g(P)$ is
clearly positive, ${\t(T\a_g(P)) \geq 0}$. Hence, $\s$ is a positive
operator. If $T \in \B_{sa}$, then $\s_T \in \A_{sa}$. If $\l$ is
the action of $G$ on $G/R$, as defined in equation
(\ref{defn:5:lamdaaction}), then it is easily seen that $\s$ is
\mbox{$\l$-$\a$-equivariant}, that is,  $\l_g\s_T=\s_{\a_g(T)}$, for
all $g \in G,\, T \in \B$. Finally, we have that ${\|\s_T\|_{\infty}
\leq \|T\|}$, for all ${T \in \B}$. To see why this is so, note that
since $\a_g(P)$ is a rank-one projection, it is of the form
$\a_g(P)x=\<x,e_0\>e_0$, where $e_0$ is some norm-one element in
image of $\a_g(P)$. Now, if $\{e_i\}_{i=0}^n$ is an orthonormal
basis of $H$ containing $e_0$, then, for all $g \in G$,
\[
\left|\s_T[g]\right| =\left| \sum_{i=0}^n\<T\a_g(P)e_i,e_i\>\right|
= \<Te_0,e_0\>\leq \|T\|\|e_0\|^2=\|T\|,
\]
and the desired result follows.

\subsubsection{The Quotient of $L$ on $\B_{sa}$}

We now define $N(a,b)=\g \inv \|f-\s_T\|_{\infty}$, for some
constant $\g$. It is clear from the definition of $L$ that
$L^q_{\B}(T)=\inf\{L(f,T):f \in \A_{sa}\} \leq L_{\B}(T)$. Thus, to
establish equality between $L^q_\B$ and $L_\B$, it would suffice to
show that, for every $T \in \B_{sa}$, there exists an $f_T \in
\A_{sa}$, \st $L(f_T,T) = L_\B(T)$. This is where the Berezin
covariant transform comes into play. For any given $T \in \B_{sa}$,
try $f_T=\s_T$. Since $\s$ is $\l$-$\a$-equivariant, we have that
\begin{eqnarray*}
L_\A(\s_T) & =\sup_{g \neq
e}\{\|\l_g(\s_T)-\s_T\|_{\infty}/\ell(g)\} = \sup_{g \neq
e}\{\|\s_{(\a_g(T)-T)}\|_{\infty}/\ell(g)\}\\
 & \leq \sup_{g \neq e}\{\|\a_g(T)-T\|/\ell(g)\}= L_{\B}(T),
\end{eqnarray*}
Thus,
\[
L(\s_T,T)=L_A(\s_T)\vee L_{B}(T) \vee \g \inv
\|\s_T-\s_T\|_{\infty}= L_{B}(T);
\]
and so $L_{\B}^q(T)=L_{\B}$, for all choices of $\g$.

The quotient of $L$ on $\A_{sa}$ is not as easy to calculate. In
fact, we shall only be able to prove that it is equal to $L_A$ for a
certain adequately large values of $\g$. In order to calculate this
value we shall need to introduce a suitably defined adjoint of $\s$
called the Berezin contravariant transform.

\subsection{The Berezin Contravariant Transform}

In this section we shall make extensive use of the notion of
averaging an operator over a \cpt group. Therefore, before we
begin any presentation of the Berezin contravariant transform, it
would be wise to recall what it means to `average an operator over
a group'.

\subsubsection{Averaging Operators over Compact groups}

Let $G,U,H$ and $\a$ be as above. The compactness of $G$ implies
that the \cts mapping
\begin{eqnarray}\label{5:map:haarmeasurable}
g \mto \< y,\a_g(T) x\>
\end{eqnarray}
is Haar integrable, for all $x,y \in H$. (Note that as usual we
shall only consider the normalised Haar measure.) Thus, the mapping
\[
y \mto \int_G \<y,\a_g(T) x\> dg,
\]
is a well defined element of the \cts dual of $H$. Linearity of
the functional is obvious, and it is easily seen to be bounded. By
the Riesz representation theorem there exists a unique ${z \in H}$
\st
\begin{eqnarray}\label{5:map:avggrp}
{\<y,z\> =\int_{G}\< y,\a_g(T)x \> dg}.
\end{eqnarray}
We shall denote $z$ by $\int_{G}\a_g(T)\,x\,dg$. Consider the
operator
\[
\int_{G}\a_g(T)dg: H \to H,~~~ x \mto \int_{G}\a_g(T) \, x \, dg.
\]
It is easy to establish that it is linear and bounded, with norm
less than or equal to $\|T\|$. We call $\int_G \a_g(T) dg$ the {\em
average} of $T$ over $G$, and we denote it by $\wt{T}$. An important
point, that is easily verified, is that if $T \geq 0$, then
$\wt{T}=0$ \iff $T=0$.

\bigskip

A little thought will verify that the map defined in
(\ref{5:map:haarmeasurable}) can be replaced by the map ${g \mto
\<y, A\a_g(T)x\>}$, for any $A \in B(H)$, and that a well defined
meaning can then be ascribed to $\int_G A \a_g(T)dg$ as an element
of $B(H)$. Similarly, a well defined meaning can be ascribed to
$\int_G \a_g(T)Adg$. Let us note that since
\begin{eqnarray*}
\<\int_G A \a_g(T)\, x dg, y\> & = & \int_G \<\a_g(T)\,x,A^*y\>dg = \<\int_G \a_g(T)\,xdg ,A^*y\>\\
                               & = & \<A\int_G \a_g(T)\,x dg,y\>,
\end{eqnarray*}
for all $x,y \in H$, it holds that $\int A\a_g(T)dg =
A\int_G\a_g(T)dg$. It is also easily seen that $\int \a_g(T)Adg=\int
\a_g(T)dgA$.

An important consequence of these two results is that ${U_h \wt{T} =
\wt{T} U_h}$, for all $h \in G$. This can be seen from
\[
U_h \wt{T}=\int U_hU_gTU^*_gdg =\int U_{hg}TU^*_{hg}U^*_{h
\inv}dg=\wt{T}U_h.
\]
Thus, if $T$ is non-zero, then the ergodicity of $\a$ implies that
$\wt{T}=\l 1$, for some $\l \in \bC$. With a view to finding a value
for $\l$, consider the trace of $\int_G S_g dg$, where
${S_g=A\a_g(T)B}$, for some ${A,B \in B(H)}$. If $\{e_i\}_{i=1}^n$
is an orthonormal basis of $H$, then
\begin{eqnarray*}
\t (\int_{G} S_g dg) & = & \sum_{i=1}^n \< \int_G S_g e_idg,e_i\>=\sum_{i=1}^n \int_G \<S_ge_i,e_i\>dg  \\
                     & = & \int_G \sum_{i=1}^n \<S_ge_i,e_i\>dg =\int_G \t(S_g) dg.
\end{eqnarray*}
This implies that
\[
\t(\wt{T})=\t(\int_G U_gTU^*_g dg) = \int_G\t(U_gTU^*_g)dg=\int_G
\t(T)dg=\t(T).
\]
Since $\wt{T}=\l 1$, we also have that $\t(\wt{T})=\t(\l 1) = \l
n$, where $n$ is the dimension of $H$. Thus,
\begin{equation}\label{5:eqn:intidentity}
\wt{T} = \frac{\t(T)}{n}1.
\end{equation}


\subsubsection{The Berezin Contravariant Transform}


Endow $\B$ with the {\em Hilbert--Schmidt} inner product, which is
defined by setting
\[
\<T,S\>_{HS}=\frac{1}{n}\t(TS^*).
\]

Let $\mu$ be the Haar measure on $G$, and let $\pi$ be the canonical
projection from $G$ to $G/R$. We denote by $L^2(G/R)$ the linear
space of equivalence classes of Borel measurable \fns on $G/R$ that
are square integrable \wrt the measure $\wt{\mu}=\mu \circ \pi
\inv$. We endow $L^2(G/R)$  with its standard inner product, as
defined in Section \ref{1:sect:NCMT}. Since $G/R$ is compact,
$C(G/R) \sseq L^2(G/R)$. Thus, $\s$ can be viewed as a linear
mapping from the Hilbert space $\B$ to the Hilbert space $L^2(G/R)$.
This means that there exists an operator $\breve{\s}:L^2(G/R) \to
\B$ \st
\[
\<\s_T,f\>_{L^2}=\<T,\breve{\s}_f\>_{HS},
\]
for all $f \in L^2(G/R),~T \in B(H)$. We call $\breve{\s}$ the {\em
Berezin contravariant mapping}, and we call  $\breve{\s}_f$ the {\em
Berezin contravariant symbol} of $f$. We shall only consider the
restriction of $\breve{\s}$ to $\A$, which we denote by the same
symbol. This mapping is often viewed as a `quantization' operator
since it brings \fns to operators. It is related to the Toeplitz
maps discussed in the previous chapter.

\bigskip

Using the results that we established above for the average of an
operator over a group, we shall find a more explicit formulation of
$\breve{\s}$. To begin with, we note that we can regard $C(G/R)$ as
a subset of $C(G)$ (in the sense that there exists a canonical
embedding of $C(G/R)$ into $C(G)$; namely the mapping $f \mto
\wt{f}=f \circ \pi$). It proves profitable to do so since
\[
\int_G \wt{f} d\mu = \int_{G/R} f d\wt{\mu};
\]
as can be verified by a routine investigation. In what follows we
shall tacitly assume this observation, and we shall not
distinguish notationally between $f$ and $\wt{f}$.

For any $f \in \A$,~$T \in \B$, we have that
\begin{eqnarray*}
\frac{1}{n}\t (\breve{\s}_fT^*) & =& \<{\breve \s}_f,T\>_{HS} =
\<f,\s_T\>_{L^2} = \int_G f(g)\ol{(\s_T(g))}dg \\
& =& \int_G f(g)\tau(\a_g(P)T^*)dg = \tau \left(\int_G
f(g)\a_g(P)dg\;T^* \right).
\end{eqnarray*}
Since this is true for all $T$, it must hold that
\begin{equation} \label{6:eqn:brexbreve}
\breve{\s}_f= n \int_G f(g)\a_g(P)dg.
\end{equation}
Thus, since
\begin{eqnarray*}
\a_h (\breve{\s}_f) & = &n U_h \int_G f(g) \a_g(P)dg U_h^* = n
\int_G f(g)
\a_{hg}(P)dg \\
                    & = &n \int_G f(h \inv g) \a_g(P) dg = \breve{\s}_{(\l_h f)},\\
\end{eqnarray*}
it also holds that $\breve{\s}$ is $\a$-$\l$-equivariant. Following
similar lines of argument we can also show that $\breve{\s}$ is
norm-decreasing.

\subsubsection{The Quotient of $L$ on $\A_{sa}$}


Now, that we have constructed the Berezin contravariant mapping, we
are ready to approach the question of the quotient of $L$ on
$\A_{sa}$. For convenience sake we shall recall here that
\[
L(f,T)=L_{\A}(f)\vee L_{\B}(T) \vee \g \inv \|f-\s_T\|_{\infty}.
\]

An immediate consequence of the definition of $L$ is that
\[
L_\A(f) \leq L_\A^q(f)=\inf\{L(f,T):T \in \B\}.
\]
Thus, to establish equality between $L^q_\A$ and $L_\A$, it would
suffice to show that, for each $f \in \A$, there exists a $T_f \in
B(H)$ \st $L(f,T_f)\leq L_\A(f)$. Recalling the use we made of the
Berezin covariant symbol when we were examining the quotient of
$L$ on $B$, it seems reasonable to try $T_f=\breve{\s}_f$. Since
$\breve{\s}$ is norm-decreasing and $\a$-$\l$-equivariant, it
holds that
\begin{eqnarray*}
L_\B(\breve{\s}_f) & =\sup_{g \neq e}\|\a_g(\breve{\s}_f)-\breve{\s}_f\|/\ell(g) =\sup_{g \neq e}\|\breve{\s}_{(\l_g(f)-f)}\|/\ell(g)\\
                   & \leq \sup_{g \neq e}\|\l_g(f)-f\|/\ell(g)=L_{\A}(f).
\end{eqnarray*}
Thus, $L_{\B}(\breve{\s}_f) \leq L_\A(f)$, for all values of $\g$.
This means that if we could find a value for $\g$ \st
\begin{equation} \label{5:eqn:gamma}
L_{\A}(f) \geq \g \inv \|f - \s(\breve{\s}_f)\|_{\infty},
\end{equation}
then for the corresponding $L$, it would hold that
$L^q_{\A}=L_{\A}$. We shall spend the remainder of this section
trying to find such a value.

\bigskip

The map $f \mto \s(\breve{\s}_f)$ is called the {\em Berezin
transform}. We can derive a more explicit formulation of it as
follows:
\begin{eqnarray*}
(\s({\breve \s}_f))[h] &= & \tau({\breve \s}_f\a_h(P)) = \tau\left(n
\int_G f(g)\a_g(P)dg\;\a_h(P)\right) \\
&=& n \int_{G} f(g)\tau(\a_g(P)\a_h(P))dg = n \int_{G}
f(g)\tau(P\a_{g^{-1}h}(P))dg.
\end{eqnarray*}
For any rank-one projection $P$ on $H$, we shall find it useful to
introduce a \fn $k_P \in C(G/R)$ defined by
\[
k_P[g] = n \t(P\a_g(P)).
\]
(Note that $ k_P[g]=n \s_P[g]$; we use a distinct symbol for $k_P$
for sake of presentation.) Our formula for the Berezin transform now
becomes
\[
(\s({\breve \s}_f))[h] =  \int_G f(g)k_P(g \inv h) dg.
\]
The \fn $k_P$ has some pleasing properties: for any norm-one vector
$e_0$ contained in the image of $P$,
\begin{equation}\label{6:eqn:kppoasitive}
k_P([g])=n|\<U_ge_0,e_0\>|^2 \geq 0.
\end{equation}
(This is easily established by choosing a specific orthonormal basis
for $H$ that contains $e_0$.) Thus, $k_p$ is a positive \fnn. Using
equation (\ref{6:eqn:kppoasitive}) we can easily show that $k_P([g
\inv])=k_P([g])$. Finally, we also have that
\[
\int_{G}k_P(g)dg= \int_{G}n \t(P\a_g(P))dg = \t \left(
Pn\int_{G}\a_g(P)dg \right)=\t(P\breve{\s}_1)=1.
\]
We shall tacitly make use of these observations below.

\bigskip

We are now in a position to find a value for $\gamma$ for which
equation (\ref{5:eqn:gamma}) will be satisfied. To begin with, let
us note that
\begin{eqnarray*}
|f([h])-\s(\breve{\s}_f)([h])|& =  & \left| \int_{G}(f([h])-f(g))k_P(g \inv h)dg \right| \\
                              & \leq & \int_{G} \left|f([h])-f(g)\right|k_P(g \inv h) dg.
\end{eqnarray*}
Now, if $f$ is an element of the dense order-unit subspace of
$\A_{sa}$ on which the Lipschitz seminorm takes finite values, then
 $|f[h]-f[g]| \leq L_\A(f) \r_{\pi}([h],[g])$, for all
$g,h \in G$. Therefore,
\begin{eqnarray*}
|f([h])-\s(\breve{\s}_f)([h])| & \leq &  L_\A(f)\int_G \r_\pi([h],[g])k_P(g \inv h)dg\\
                               &   =  &  L_\A(f)\int_G \r_\pi([h],[g]) k_P(h \inv g)dg.\\
\end{eqnarray*}
Since we required the length \fn $\ell$ on $G$ to satisfy
$\ell(xg,xh) = \ell(g,h)$, for all $x,g,h \in G$, we have that
$\r_\pi([xg],[xh])=\r_\pi([g],[h])$. Consequently,
\begin{eqnarray*}
|f([h])-\s(\breve{\s}_f)([h])|  &  \leq  & L_\A(f)\int_G \r_\pi([h],[hg]) k_P(g)dg\\
                                &  =   & L_\A(f)\int_G \r_\pi([e],[g]) k_P(g)dg.\\
\end{eqnarray*}
Thus, if we choose
\[
\g= n\int_{G}\r_\pi([e],[g])\s_P(g)dg,
\]
then $\|f-\s(\breve{\s}_f)\|_{\infty} \leq \g L_\A(f)$, for all $f
\in \A_{sa}$. This gives us the following proposition.

\begin{prop}
For $\g$ chosen as above, the seminorm on $\A_{sa} \oplus \B_{sa}$
defined by setting
\[
L(f,T)=L_\A(f) \vee L_\B(T) \vee \g \inv \|f-\s_T\|_{\infty},
\]
has $L_\A$ as its quotient norm on $\A_{sa}$.
\end{prop}

\subsection{Estimating the QGH Distance}

Now, that we have shown that $L$ is an admissible Lip-norm on
$\A_{sa} \oplus \B_{sa}$, we shall try and estimate the quantum
Gromov--Hausdorff distance between $(\A_{sa},L_\A)$ and
$(\B_{sa},L_\B)$.

\begin{prop} \label{5:prop:gammanbdofB}
If $\g$ is a constant chosen \st the quotient of $L$ on $\A_{sa}$ is
$L_\A$, then $S(\A_{sa})$ is in the $\g$-neighborhood of
$S(\B_{sa})$ for $\r_L$.
\end{prop}

\demo For each~ $\mu \in S(\A_{sa})$~ we~ must~ produce a $\nu \in
S(\B_{sa})$, such~ that {${\r_L(\mu,\nu) \leq \g}$}. Let us try
$\nu=\mu \circ \s$. (Note that $\s$ is unital and positive, and
therefore $\nu$ is indeed contained in $S(\B_{sa})$). Recall that
\[
\r_L(\mu,\nu)=\sup \{|\mu(f,T)-\nu(f,T)|:(f,T)\in \A_{sa} \oplus
\B_{sa},\, L(f,T) \leq 1\}.
\]
For $(f,T) \in \A_{sa} \oplus \B_{sa}$,
\[
\begin{tabular}{ll}
$|\mu(f,T)-\nu(f,T)|$ & $=|\mu(f)-\nu(T)|=|\mu(f)-\mu(\s_T)|$\\
                      & $=|\mu(f-\s_T)| \leq \|\mu\|\|f-\s_T\|_{\infty}$\\
                      & $=\|f-\s_T\|_{\infty}.$\\
\end{tabular}
\]
Since $L(f,t) \leq 1$ and $L(f,t)=L(f) \vee L(T) \vee \g\inv
\|f-\s_T\|_\infty,$ we have ${\|f-\s_T\|_\infty \leq \g}$. It
follows that $\r_L(\mu,\nu) \leq \g$, and so $S(\A_{sa})$ is
contained in the $\g$-\nbd of $S(\B_{sa})$. \qed

Consequently, to put a suitably small upper bound on the quantum
Gromov--Hausdorff distance between $(\A_{sa},L_\A)$ and
$(\B_{sa},L_\B)$, it only remains to show that $S(\B_{sa})$ is
contained in a suitably small \nbd of $S(\A_{sa})$. That is, for
each $\nu \in S(\B_{sa})$, we must find a $\mu \in S(\A_{sa})$ \st
$\dist_q(\nu,\mu)$ is suitably small. Mimicking the proof of
Proposition \ref{5:prop:gammanbdofB}, we propose $\mu = \nu \circ
\breve{\s}$. If ${(f,T) \in \A_{sa} \oplus \B_{sa}}$, and $L(f,T)
\leq 1$, then $L_B(T)\leq 1$ and {$\|f-\s_T\| \leq \g$}. Thus,
\[
\begin{tabular}{ll}
$|\mu (f,T)-\nu(f,T)|$ &
                                         $=|\nu(\breve{\s}_f)-\nu(T)| \leq
\|\nu\|\|\breve{\s}_f-T\|$\\
                                       & $\leq
\|\breve{\s}_f-\breve{\s}(\s_T)\|+\|\breve{\s}(\s_T)-T\|$\\
                                       & $\leq \|f-\s_T\|_{\infty}
                                       +\|\breve{\s}(\s_T)-T\|$\\
                                       & $\leq \g +
                                       \|\breve{\s}(\s_T)-T\|$.
\end{tabular}
\]
This means that any bound that we can obtain on
$\|\breve{\s}(\s_T)-T\|$, for $L_B(T) \leq 1$, will give us a bound
on the quantum Gromov--Hausdorff distance between $(\A_{sa},L_\A)$
and $(\B_{sa},L_\B)$.

\subsection{Matrix Algebras Converging to the Sphere}

Let us summarise what we have established: If $(\A_\ssa,L_\A)$ is
the \cqms associated to $\O_\mu$, and $(\B_\ssa,L_\B)$ is the \cqms
associated to the $n$-fuzzification of $\O_\mu$, then the \qghd
between them is less than or equal to $ \g +\|\breve{\s}(\s_T)-T\|$;
where $T \in \B_\ssa$ \st $L_B(T) \leq 1$, and $\g =
n\int_{G}\r_\pi([e],g)\s_P(g)dg$.

Rieffel went on to show that $\g$ and $\|\breve{\s}(\s_T)-T\|$ are
dependent on $n$, and that as $n \to \infty$, $\g$ and
$\|\breve{\s}(\s_T)-T\|$ become arbitrarily small. Thus, \wrt
quantum Gromov--Hausdorff distance, the sequence of fuzzy
coadjoint orbits converges to the $\O_\mu$. (We shall not outline
the proof because it quite lengthy and would require the
introduction of an excessive amount of Lie group theory; for
details see \cite{RIEF2}.) A precise meaning has now been given to
statements involving the convergence of fuzzy spaces to a
coadjoint orbit. Moreover, since $S^2$ is a coadjoint orbit of
$SU(2)$, a precise meaning has also been given to statements
involving the convergence  of fuzzy spaces to the $2$-sphere.

\section{Matricial Gromov--Hausdorff Distance}

As we discussed earlier, a shortcoming of quantum Gromov--Hausdorff
distance is that two Lip-normed $C^*$-algebras can have distance
zero yet their \calgs may not be isomorphic. Following a suggestion
of Rieffel, David Kerr \cite{KERR} began to investigate the
possibility of defining a modified version of quantum
Gromov--Hausdorff using `matrix-valued states'. In this context the
notion of positivity gives way to the notion of complete positivity.
We shall now review this concept.

\subsubsection{Complete Positivity and Operator Systems}

Let $\A$ be a \calgn, and let $M_n(\A)$ denote the \alg of all  $n
\times n$ matrices with entries in $\A$. We can define an involution
on $M_n(\A)$ by setting $[a_{ij}]^*=[a_{ji}^*]$. If $\f$ is a
mapping from $\A$ to another \calg $\B$, then we define its {\em
$n$-inflation} to be the mapping
\[
\f_n:M_n(\A) \to M_n(\B),~~~[a_{ij}] \to [\f(a_{ij})].
\]
Note that if $\f$ is a $*$-\alg homomorphism, then its $n$-inflation
is also a  ${*\text{-\alg}}$ homomorphism, for all $n$.

Let $H$ be a Hilbert space, and consider the mapping
\[
\ps:M_n(B(H)) \to B(H^n),~~~u \mto \ps(u);
\]
where $H^n$ is the orthogonal $n$-sum of $H$, and $\ps(u)$ is the
operator defined by setting
\[
\ps(u)(x_1,x_2, \ldots,x_n)=(\sum_{j=1}^nu_{1j}(x_j), \ldots ,
\sum_{j=1}^n u_{nj}(x_j)).
\]
It is straightforward to show that $\ps$ is a $*$-\alg isomorphism.
This means that we can define a norm $\|\cdot\|$ on $M_n(B(H))$ that
makes it a \calgn, by setting $\|u\|=\|\ps(u)\|$. The following
useful inequalities are easily established:
\begin{equation}\label{5:eqn:cmatrix}
\|u_{ij}\| \leq \|u\|,~~~~ i,j=1, \ldots ,n.
\end{equation}

Let $\pi$ be a faithful representation of $\A$ in $B(H)$, for some
Hilbert space $H$. (As discussed in Chapter $1$, such a $\pi$ and
$H$ can be always be produced using the GNS construction.) Let
$\pi_p$ be the $p$-inflation of $\pi$, and let  $u(n)$ be a sequence
in $\pi_p(M_p(\A))$ that converges to some $u \in M_p(B(H))$. By
(\ref{5:eqn:cmatrix}), we have that $u_{ij}(n)$ converges to
$u_{ij}$, for each $i,j=1, \ldots ,n$. Since $\pi(\A)$ is complete,
each $u_{ij}$ is contained in $\pi(\A)$. Thus, $\pi(M_p(\A))$ is
closed in $M_p(B(H))$, implying that $\pi(M_p(\A))$ is a \calgn.
This enables us to define a norm on $M_p(\A)$, that makes it a
\calgn, by setting
\[
\|v\|=\|\pi(v)\|,~~~~~v \in M_p(\A).
\]
As discussed in Chapter $1$, it is the unique norm on $M_p(\A)$ that
does so.

\bigskip

Let $\f$ be a map from $\A$ to another \calg $\B$. Then $\f$ is
called {\em completely-positive} if its $n$-inflation $\f_n:M_n(\A)
\to M_n(\B)$ is positive, for all $n$. Not all \pve mappings are
completely-positive, the standard example of a positive,
non-completely-positive, mapping is
\[
\f:M_2(\bC) \to M_2(\bC),~~~~ A \to A^T.
\]

Complete-positivity is usually studied in the more general setting
of operator systems. A {\em (concrete) operator system} is a unital
self-adjoint closed linear subspace of a unital \calgn. (We use the
adjective {\em concrete} here because there exists a more general
definition of an operator system \cite{EFROS}; each such structure
can, however, be represented as a self-adjoint linear subspace of a
\calgn. Kerr works with the concrete definition, and so it is the
formulation that we shall use here.) Note that for any operator
system $X$, its subset of self-adjoint elements, which we denote by
$X_{\ssa}$, is an order-unit space.

If $\A$ is a \calgn, and $X \sseq \A$ is an operator system, then we
denote by $M_n(X)$ the subset of $M_n(\A)$ whose elements are the
matrices with entries in $X$. It is clear that $M_n(X)$ is a unital
self-adjoint linear subspace of $M_n(\A)$, and so it is also an
operator system. We define the set of {\em positive} elements of $X$
to be $X \cap \A_+$, and we define the set of $n$-\pve and
completely-\pve maps between two operator systems in exactly the
same way as for \mbox{\calgsn}. Finally, we denote by $UCP_n(X)$,
the set of unital completely-positive maps from  $X$ to $M_n(\bC)$;
or, more explicitly, the set of maps $\f:X \to M_n(\bC)$ whose
inflation
\[
\f_p:M_p(X) \to M_p(M_n(\bC)) \simeq M_{p \times n}(\bC)
\]
is \pve and unital, for all $p$. Since each $\f \in UCP_1(X)$ is
clearly a state, the following lemma tells us that $UCP_1(X)$ is
equal to the state space of $X$; for a proof see \cite{CONW}.
\begin{lem}
If $\f$ is a state on $X$, then $\f$ is a unital completely-positive
map.
\end{lem}

\subsubsection{Lip-Normed Operator Systems}

We shall now define the analogue for operator systems of  \cqmsn s:
let $(X,L)$  be a pair consisting of an operator system $X$, and a
Lip-norm $L$ defined on $X_{sa}$; if $D_1(L)=\{x \in D(L):  L(x)
\leq 1\}$ is closed in $X_{sa}$, then we say that $(X,L)$ is a {\em
Lip-normed operator system}. (The technical requirement that
$D_1(L)$ be closed will not be of great importance to us here; it is
included for the sake of accuracy.) Clearly, every Lip-normed unital
\calg $(\A,L)$ (for which $D_1(L)$ is closed) is a Lip-normed
operator system.

\bigskip

Mimicking the manner in which we defined a metric on the state space
of an order-unit space using a Lip-norm, we define a metric
$\r_{L,n}$ on $UCP_n(X)$, for each $n$, by setting
\[
\rho_{L,n} (\varphi , \psi ) = \sup\{ \| \varphi (a) - \psi (a) \|:
L(a) \leq 1\},
\]
for $\varphi , \psi\in UCP_n (X)$. (Note that by $\| \cdot \|$ we
mean the unique norm on $M_n(\bC)$ that makes it a \calgn.)

In the order-unit case we required that the metric $\r_L$ induce the
weak$^*$ topology on $S(A)$. For $UCP_n(X)$ the natural analogue of
the weak$^*$ topology is the {\em point-norm topology}; it is
defined to be the weakest topology \wrt which the family of \fns
$\{\wh{x}:x \in X\}$ is \ctsn, where $\wh{x}(\f)=\f(x)$. Just as the
state space is \cpt \wrt the weak$^*$ topology, each $UCP_n(X)$ is
point-norm \cptn. It would be natural to require that each
$\r_{L,n}$ induce the point-norm topology on $UCP_n(X)$. However,
Kerr established that this is a consequence of the fact that $L$ is
a Lip-norm on $X_{sa}$. Therefore, there is no need to impose such a
condition.

\subsubsection{Complete Gromov--Hausdorff Distance}

The definition of quantum Gromov--Hausdorff distance involves the
direct sum of two order-unit spaces, and the embedding of their
state spaces into the state space of their direct sum. We shall now
translate this process to the operator system setting. Let $(X,L_X)$
and $(Y,L_Y)$ be Lip-normed operator systems, and let $\A$ and $\B$
be two \calgs containing $X$ and $Y$ respectively. The {\em direct
sum} of $\A$ and $\B$ is defined to be their direct sum as normed
\algs endowed with the pointwise-defined addition and involution; it
is denoted by $\A \oplus \B$. Clearly, $\A \oplus \B$ is a \calgn.
The direct sum of $X$ and $Y$ as normed linear spaces, endowed with
the pointwise-defined multiplication and involution, is clearly a
unital self-adjoint linear subspace of $\A \oplus \B$. Hence, it is
an operator space. We shall denote it by $X \oplus Y$.

Now, $UCP_n(X)$ can be embedded into $UPC_n(X \oplus Y)$ in an
obvious manner, and it is is easily seen that this embedding is \cts
when we put the point-norm topology on both spaces. Since $UCP_n(X)$
and $UPC_n(X \oplus Y)$ are both \cpt Hausdorff spaces, the image of
$UCP_n(X)$ in $UCP_n(X \oplus Y)$, which we shall equate with
$UCP_n(X)$, is closed. Obviously, an entirely analogous situation
holds for $UCP_n(Y)$.

We can speak of admissable Lip-norms on $X \oplus Y$, since it is
just the sum of two order-unit subspaces. We denote the set of
admissable Lip-norms on $X \oplus Y$ by $\M(L_X,L_Y)$. Kerr showed
that if $L \in \M(L_X,L_Y)$, then the restriction of $\r_{L,n}$ to
$UCP_n(X)$ is equal to $\r_{L_X,n}$, and the restriction of
$\r_{L,n}$ to $UCP_n(Y)$ is equal to $\r_{L_Y,n}$. This is a direct
and pleasant generalisation of what happens in the order-unit case.
We can now imitate the definition of quantum Gromov--Hausdorff
distance. For each natural number $n$, we define $\dist^n(X,Y)$ the
{\em $n$-distance} between $(X,L_X)$ and $(Y,L_Y)$ by setting
\[
\dist^n(X,Y) = \inf\{\r_H^{\rho_{L,n}}(UCP_n(X),UCP_n(Y)): L \in
\M(L_X,L_Y)\}.
\]
We then define $\dist_c(X,Y)$ the {\em complete quantum
Gromov--Hausdorff distance} by setting
\[
\dist_c(X,Y) = \sup_{n \in \bN} \{\dist^n(X,Y)\}.
\]

If we bear in mind that $UCP_1(X)=S(X)$, then a little careful
reflection will verify that the quantum Gromov--Hausdorff distance
between two Lip-normed \mbox{$C^*$-algebras} is equal to their
$1$-distance. Thus, the complete distance is always greater than or
equal to the quantum Gromov--Hausdorff distance.

Clearly, the complete distance is symmetric in its arguments. Kerr
showed that it also satisfies the triangle inequality and is
positive definite on the family of appropriately defined
equivalences classes of Lip-normed operator systems. Hence, it is
well defined as a metric. A direct consequence of the proof of
positive definiteness is that two \calgs have complete distance zero
\iff they are $*$-isomorphic. Thus, Kerr's definition overcomes the
shortcoming of Rieffel's definition.

Kerr also showed that the continuity of quantum Gromov--Hausdorff
distance for non-commutative tori, as described earlier, carries
over to the complete quantum Gromov--Hausdorff distance case; as
does the convergence of matrix algebras to coadjoint orbits
described in Section $6.3$. In \cite{KERRLI} it was shown that the
family of equivalences classes of Lip-normed operator systems
endowed with $\dist_c$ is a complete metric space.

\subsubsection{Operator Gromov--Hausdorff Distance}

Hangfeng Li, a doctoral student of Rieffel, devised another strategy
for quantizing Gromov--Hausdorff distance that operates entirely at
the algebraic level. It also overcomes the shortcoming of Rieffel's
distance addressed above. His versatile approach was implemented in
both the order-unit and $C^*$-algebraic contexts under the
terminology {\em order-unit}, and {\em $C^*$-algebraic quantum
Gromov--Hausdorff distance} respectively \cite{LI1,LI2}. It affords
many technical advantages.

In a recent paper \cite{KERRLI} Kerr, working jointly with Li,
established an analogue for Lip-normed operator systems of Li's
distance. The pair then proved that this new distance is in fact
equal to the complete \qghdn. This consolidation of complete
Gromov--Hausdorff distance has motivated Kerr and Li to propose that
it be renamed {\em operator Gromov--Hausdorff distance}.

\subsubsection{Completeness and Lip-Ultraproducts}

Since it is primarily \calgsn, as opposed to operator systems, that
we are interested in, it would be pleasing if the subfamily of
Lip-normed \calgs were closed in the family of Lip-normed operator
systems. In \cite{KERRLI} Kerr and Li produced sufficient conditions
for a sequence of \calgs to converge to a \calgn. However, in a
recent paper Daniele Guido and Tommaso Isola \cite{GUIDO} (both
members of the European Union Operator Algebras Network) constructed
a Cauchy sequence of \calgs that converges, \wrt \qghdn, to an
operator system that is not a \mbox{\calgn}. Hence, the space of
Lip-normed \calgs is not complete \wrt $\dist_c$. Guido and Isola's
work is based upon their newly defined notion of a
`Lip-ultraproduct'. For natural reasons they propose that it be
viewed as the quantum analogue of the ultralimit of a sequence of
compact metric spaces (for details on ultralimits see
\cite{BRIDHAE}). Their work has lead them to define a new metric on
the space of Lip-normed \calgs \wrt which it is complete.
Consequently, they propose it as a more natural way to define
distance.

\end{document}